\input amstex
\documentstyle{amsppt}
\input miniltx.tex
\input graphicx.sty
\input epsf.tex
\def\stydate{May 10, 2002}

\chardef\tempcat\catcode`\@
\ifx\undefined\amstexloaded\input amstex \else\catcode`\@\tempcat\fi
\expandafter\ifx\csname amsppt.sty\endcsname\relax\input amsppt.sty \fi
\let\tempcat\undefined

\immediate\write16{This is LABEL.DEF by A.Degtyarev <\stydate>}
\expandafter\ifx\csname label.def\endcsname\relax\else
  \message{[already loaded]}\endinput\fi
\expandafter\edef\csname label.def\endcsname{%
  \catcode`\noexpand\@\the\catcode`\@\edef\noexpand\styname{LABEL.DEF}%
  \def\expandafter\noexpand\csname label.def\endcsname{\stydate}%
    \toks0{}\toks2{}}
\catcode`\@11
\def\labelmesg@ {LABEL.DEF: }
{\edef\temp{\the\everyjob\W@{\labelmesg@<\stydate>}}
\global\everyjob\expandafter{\temp}}

\def\@car#1#2\@nil{#1}
\def\@cdr#1#2\@nil{#2}
\def\eat@bs{\expandafter\eat@\string}
\def\eat@ii#1#2{}
\def\eat@iii#1#2#3{}
\def\eat@iv#1#2#3#4{}
\def\@DO#1#2\@{\expandafter#1\csname\eat@bs#2\endcsname}
\def\@N#1\@{\csname\eat@bs#1\endcsname}
\def\@Nx{\@DO\noexpand}
\def\@Name#1\@{\if\@undefined#1\@\else\@N#1\@\fi}
\def\@Ndef{\@DO\def}
\def\@Ngdef{\global\@Ndef}
\def\@Nedef{\@DO\edef}
\def\@Nxdef{\global\@Nedef}
\def\@Nlet{\@DO\let}
\def\@undefined#1\@{\@DO\ifx#1\@\relax\true@\else\false@\fi}
\def\@@addto#1#2{{\toks@\expandafter{#1#2}\xdef#1{\the\toks@}}}
\def\@@addparm#1#2{{\toks@\expandafter{#1{##1}#2}%
    \edef#1{\gdef\noexpand#1####1{\the\toks@}}#1}}
\def\make@letter{\edef\t@mpcat{\catcode`\@\the\catcode`\@}\catcode`\@11 }
\def\donext@{\expandafter\egroup\next@}
\def\x@notempty#1{\expandafter\notempty\expandafter{#1}}
\def\lc@def#1#2{\edef#1{#2}%
    \lowercase\expandafter{\expandafter\edef\expandafter#1\expandafter{#1}}}
\newif\iffound@
\def\find@#1\in#2{\found@false
    \DNii@{\ifx\next\@nil\let\next\eat@\else\let\next\nextiv@\fi\next}%
    \edef\nextiii@{#1}\def\nextiv@##1,{%
    \edef\next{##1}\ifx\nextiii@\next\found@true\fi\FN@\nextii@}%
    \expandafter\nextiv@#2,\@nil}
{\let\head\relax\let\specialhead\relax\let\subhead\relax
\let\subsubhead\relax\let\proclaim\relax
\gdef\let@relax{\let\head\relax\let\specialhead\relax\let\subhead\relax
    \let\subsubhead\relax\let\proclaim\relax}}
\newskip\@savsk
\let\@ignorespaces\ignorespaces
\def\@ignorespacesp{\ifhmode
  \ifdim\lastskip>\z@\else\penalty\@M\hskip-1sp%
        \penalty\@M\hskip1sp \fi\fi\@ignorespaces}
\def\ignorespaces{\protect\@ignorespacesp}
\def\@bsphack{\relax\ifmmode\else\@savsk\lastskip
  \ifhmode\edef\@sf{\spacefactor\the\spacefactor}\fi\fi}
\def\@esphack{\relax
  \ifx\penalty@\penalty\else\penalty\@M\fi   
  \ifmmode\else\ifhmode\@sf{}\ifdim\@savsk>\z@\@ignorespacesp\fi\fi\fi}
\let\@frills@\identity@
\let\@txtopt@\identyty@
\newif\if@star
\newif\if@write\@writetrue
\def\@numopt@{\if@star\expandafter\eat@\fi}
\def\checkstar@#1{\DN@{\@writetrue
  \ifx\next*\DN@####1{\@startrue\checkstar@@{#1}}%
      \else\DN@{\@starfalse#1}\fi\next@}\FN@\next@}
\def\checkstar@@#1{\DN@{%
  \ifx\next*\DN@####1{\@writefalse#1}%
      \else\DN@{\@writetrue#1}\fi\next@}\FN@\next@}
\def\checkfrills@#1{\DN@{%
  \ifx\next\nofrills\DN@####1{#1}\def\@frills@####1{####1\nofrills}%
      \else\DN@{#1}\let\@frills@\identity@\fi\next@}\FN@\next@}
\def\checkbrack@#1{\DN@{%
    \ifx\next[\DN@[####1]{\def\@txtopt@########1{####1}#1}%
    \else\DN@{\let\@txtopt@\identity@#1}\fi\next@}\FN@\next@}
\def\check@therstyle#1#2{\bgroup\DN@{#1}\ifx\@txtopt@\identity@\else
        \DNii@##1\@therstyle{}\def\@therstyle{\DN@{#2}\nextii@}%
    \expandafter\expandafter\expandafter\nextii@\@txtopt@\@therstyle.\@therstyle
    \fi\donext@}

\newread\@inputcheck
\def\@input#1{\openin\@inputcheck #1 \ifeof\@inputcheck \W@
  {No file `#1'.}\else\closein\@inputcheck \relax\input #1 \fi}

\def\loadstyle#1{\edef\next{#1}%
    \DN@##1.##2\@nil{\if\notempty{##2}\else\def\next{##1.sty}\fi}%
    \expandafter\next@\next.\@nil\lc@def\next@\next
    \expandafter\ifx\csname\next@\endcsname\relax\input\next\fi}

\let\pagebody@\pagebody
\let\pagetop@\empty
\let\pagebot@\empty
\let\@Xend\empty
\def\pagebody{\pagetop@\pagebody@\pagebot@\@Xend}
\let\@Xclose\empty

\newwrite\@Xmain
\newwrite\@Xsub
\def\W@X{\write\@Xout}
\def\make@Xmain{\global\let\@Xout\@Xmain\global\let\end\endmain@
  \xdef\@Xname{\jobname}\xdef\@inputname{\jobname}}
\begingroup
\catcode`\(\the\catcode`\{\catcode`\{12
\catcode`\)\the\catcode`\}\catcode`\}12
\gdef\W@count#1((\lc@def\@tempa(#1)%
    \def\\##1(\W@X(\global##1\the##1))%
    \edef\@tempa(\W@X(%
        \string\expandafter\gdef\string\csname\space\@tempa\string\endcsname{)%
        \\\pageno\\\cnt@toc\\\cnt@idx\\\cnt@glo\\\footmarkcount@
        \@Xclose\W@X(}))\expandafter)\@tempa)
\endgroup
\def\readaux{\bgroup\checkbrack@\readaux@}
\let\begin\readaux
\def\readaux@{%
    \W@{>>> \labelmesg@ Run this file twice to get x-references right}%
    \global\everypar{}%
    {\def\\{\global\let}%
        \def\/##1##2{\gdef##1{\wrn@command##1##2}}%
        \disablepreambule@cs}%
    \make@Xmain{\make@letter\setboxz@h{\@input{\@txtopt@{\@Xname.aux}}%
            \lc@def\@tempa\jobname\@Name\open@\@tempa\@}}%
  \immediate\openout\@Xout\@Xname.aux%
    \immediate\W@X{\relax}\egroup}
\everypar{\global\everypar{}\readaux}
{\toks@\expandafter{\topmatter}
\global\edef\topmatter{\noexpand\readaux\the\toks@}}
\let\@@end@@\end

\def\@Xclose@{{\def\@Xend{\ifnum\insertpenalties=\z@
        \W@count{close@\@Xname}\closeout\@Xout\fi}%
    \vfill\supereject}}
\def\endmain@{\@Xclose@
    \W@{>>> \labelmesg@ Run this file twice to get x-references right}%
    \@@end@@}
\def\disablepreambule@cs{\\\disablepreambule@cs\relax}

\def\include#1{\bgroup
  \ifx\@Xout\@Xsub\DN@{\errmessage
        {\labelmesg@ Only one level of \string\include\space is supported}}%
    \else\edef\@tempb{#1}\clearpage
      \DN@##1 {\if\notempty{##1}\edef\@tempb{##1}\DN@####1\eat@ {}\fi\next@}%
    \DNii@##1.{\edef\@tempa{##1}\DN@####1\eat@.{}\next@}%
        \expandafter\next@\@tempb\eat@{} \eat@{} %
    \expandafter\nextii@\@tempb.\eat@.%
        \relaxnext@
      \if\x@notempty\@tempa
          \edef\nextii@{\write\@Xmain{%
            \noexpand\string\noexpand\@input{\@tempa.aux}}}\nextii@
        \ifx\undefined\@includelist\found@true\else
                    \find@\@tempa\in\@includelist\fi
            \iffound@\ifx\undefined\@noincllist\found@false\else
                    \find@\@tempb\in\@noincllist\fi\else\found@true\fi
            \iffound@\lc@def\@tempa\@tempa
                \if\@undefined\close@\@tempa\@\else\edef\next@{\@Nx\close@\@tempa\@}\fi
            \else\xdef\@Xname{\@tempa}\xdef\@inputname{\@tempb}%
                \W@count{open@\@Xname}\global\let\@Xout\@Xsub
            \openout\@Xout\@tempa.aux \W@X{\relax}%
            \DN@{\let\end\endinput\@input\@inputname
                    \@Xclose@\make@Xmain}\fi\fi\fi
  \donext@}
\def\includeonly#1{\edef\@includelist{#1}}
\def\noinclude#1{\edef\@noincllist{#1}}

\def\arabicnum#1{\number#1}

\def\Romannum#1{\expandafter\uppercase\expandafter{\romannumeral#1}}
\def\alphnum#1{\ifcase#1\or a\or b\or c\or d\else\@ialph{#1}\fi}
\def\@ialph#1{\ifcase#1\or \or \or \or \or e\or f\or g\or h\or i\or j\or
    k\or l\or m\or n\or o\or p\or q\or r\or s\or t\or u\or v\or w\or x\or y\or
    z\else\fi}
\def\Alphnum#1{\ifcase#1\or A\or B\or C\or D\else\@Ialph{#1}\fi}
\def\@Ialph#1{\ifcase#1\or \or \or \or \or E\or F\or G\or H\or I\or J\or
    K\or L\or M\or N\or O\or P\or Q\or R\or S\or T\or U\or V\or W\or X\or Y\or
    Z\else\fi}

\def\ST@P{step}
\def\ST@LE{style}
\def\N@M{no}
\def\F@NT{font@}
\outer\def\newcounter{\checkbrack@{\expandafter\newcounter@\@txtopt@{{}}}}
{\let\newcount\relax
\gdef\newcounter@#1#2#3{{%
    \toks@@\expandafter{\csname\eat@bs#2\N@M\endcsname}%
    \DN@{\alloc@0\count\countdef\insc@unt}%
    \ifx\@txtopt@\identity@\expandafter\next@\the\toks@@
        \else\if\notempty{#1}\global\@Nlet#2\N@M\@#1\fi\fi
    \@Nxdef\the\eat@bs#2\@{\if\@undefined\the\eat@bs#3\@\else
            \@Nx\the\eat@bs#3\@.\fi\noexpand\arabicnum\the\toks@@}%
  \@Nxdef#2\ST@P\@{}%
  \if\@undefined#3\ST@P\@\else
    \edef\next@{\noexpand\@@addto\@Nx#3\ST@P\@{%
             \global\@Nx#2\N@M\@\z@\@Nx#2\ST@P\@}}\next@\fi
    \expandafter\@@addto\expandafter\@Xclose\expandafter
        {\expandafter\\\the\toks@@}}}}
\outer\def\copycounter#1#2{%
    \@Nxdef#1\N@M\@{\@Nx#2\N@M\@}%
    \@Nxdef#1\ST@P\@{\@Nx#2\ST@P\@}%
    \@Nxdef\the\eat@bs#1\@{\@Nx\the\eat@bs#2\@}}
\outer\def\everystep{\checkstar@\everystep@}
\def\everystep@#1{\if@star\let\next@\gdef\else\let\next@\@@addto\fi
    \@DO\next@#1\ST@P\@}
\def\counterstyle#1{\@Ngdef\the\eat@bs#1\@}
\def\advancecounter#1#2{\@N#1\ST@P\@\global\advance\@N#1\N@M\@#2}
\def\setcounter#1#2{\@N#1\ST@P\@\global\@N#1\N@M\@#2}
\def\counter#1{\refstepcounter#1\printcounter#1}
\def\printcounter#1{\@N\the\eat@bs#1\@}
\def\refcounter#1{\xdef\@lastmark{\printcounter#1}}
\def\stepcounter#1{\advancecounter#1\@ne}
\def\refstepcounter#1{\stepcounter#1\refcounter#1}
\def\savecounter#1{\@Nedef#1@sav\@{\global\@N#1\N@M\@\the\@N#1\N@M\@}}
\def\restorecounter#1{\@Name#1@sav\@}

\def\warning#1#2{\W@{Warning: #1 on input line #2}}
\def\warning@#1{\warning{#1}{\the\inputlineno}}
\def\wrn@@Protect#1#2{\warning@{\string\Protect\string#1\space ignored}}
\def\wrn@@label#1#2{\warning{label `#1' multiply defined}{#2}}
\def\wrn@@ref#1#2{\warning@{label `#1' undefined}}
\def\wrn@@cite#1#2{\warning@{citation `#1' undefined}}
\def\wrn@@command#1#2{\warning@{Preamble command \string#1\space ignored}#2}
\def\wrn@@option#1#2{\warning@{Option \string#1\string#2\space is not supported}}
\def\wrn@@reference#1#2{\W@{Reference `#1' on input line \the\inputlineno}}
\def\wrn@@citation#1#2{\W@{Citation `#1' on input line \the\inputlineno}}
\let\wrn@reference\eat@ii
\let\wrn@citation\eat@ii
\def\nowarning#1{\if\@undefined\wrn@\eat@bs#1\@\wrn@option\nowarning#1\else
        \@Nlet\wrn@\eat@bs#1\@\eat@ii\fi}
\def\printwarning#1{\if\@undefined\wrn@@\eat@bs#1\@\wrn@option\printwarning#1\else
        \@Nlet\wrn@\eat@bs#1\expandafter\@\csname wrn@@\eat@bs#1\endcsname\fi}
\printwarning\Protect
\printwarning\label
\printwarning\ref
\printwarning\cite
\printwarning\command
\printwarning\option

{\catcode`\#=12\gdef\@lH{#}}
\def\@@HREF#1{}
\def\@HREF#1#2{\@@HREF{a #1}{\let\@@HREF\eat@#2}\@@HREF{/a}}
\def\@@Hf#1{file:#1} \let\@Hf\@@Hf
\def\@@Hl#1{\if\notempty{#1}\@lH#1\fi} \let\@Hl\@@Hl
\def\@@Hname#1{\@HREF{name="#1"}{}} \let\@Hname\@@Hname
\def\@@Href#1{\@HREF{href="#1"}} \let\@Href\@@Href
\ifx\undefined\pdfoutput
  \csname newcount\endcsname\pdfoutput
\else
  \def\pdflinkattr{attr{/C [0 0.9 0.9]}}
  \let\pdflinkbegin\empty
  \let\pdflinkend\empty
  \def\@pdfHf#1{file {#1}}
  \def\@pdfHl#1{name {#1}}
  \def\@pdfHname#1{\pdfdest name{#1}xyz\relax}
  \def\@pdfHref#1#2{\pdfstartlink \pdflinkattr goto #1\relax
    \pdflinkbegin#2\pdflinkend\pdfendlink}
  \def\@ifpdf#1#2{\ifnum\pdfoutput>\z@\expandafter#1\else\expandafter#2\fi}
  \def\@Hf{\@ifpdf\@pdfHf\@@Hf}
  \def\@Hl{\@ifpdf\@pdfHl\@@Hl}
  \def\@Hname{\@ifpdf\@pdfHname\@@Hname}
  \def\@Href{\@ifpdf\@pdfHref\@@Href}
\fi
\def\@Hr#1#2{\if\notempty{#1}\@Hf{#1}\fi\@Hl{#2}}
\def\@localHref#1{\@Href{\@Hr{}{#1}}}
\def\@countlast#1{\@N#1last\@}
\def\@@countref#1#2{\global\advance#2\@ne
  \@Nxdef#2last\@{\the#2}\@tocHname{#1\@countlast#2}}
\def\@countref#1{\@DO\@@countref#1@HR\@#1}

\def\Href@@#1{\@N\Href@-#1\@}
\def\Href@#1#2{\@N\Href@-#1\@{\@Hl{@#1-#2}}}
\def\Hname@#1{\@N\Hname@-#1\@}
\def\Hlast@#1{\@N\Hlast@-#1\@}
\def\cntref@#1{\global\@DO\advance\cnt@#1\@\@ne
  \@Nxdef\Hlast@-#1\@{\@DO\the\cnt@#1\@}\Hname@{#1}{@#1-\Hlast@{#1}}}
\def\HyperRefs#1{\global\@Nlet\Hlast@-#1\@\empty
  \global\@Nlet\Hname@-#1\@\@Hname
  \global\@Nlet\Href@-#1\@\@Href}
\def\NoHyperRefs#1{\global\@Nlet\Hlast@-#1\@\empty
  \global\@Nlet\Hname@-#1\@\eat@
  \global\@Nlet\Href@-#1\@\eat@}

\HyperRefs{label}
{\catcode`\-11
\gdef\@labelref#1{\Hname@-label{r@-#1}}
\gdef\@xHref#1{\Href@-label{\@Hl{r@-#1}}}
}
\HyperRefs{toc}
\def\@HR#1{\if\notempty{#1}\string\@HR{\Hlast@{toc}}{#1}\else{}\fi}



\def\bftext{\ifmmode\fam\bffam\else\bf\fi}
\let\@lastmark\empty
\let\@lastlabel\empty
\def\lastmark{\@lastmark}
\let\lastlabel\empty
\let\everylabel\relax
\let\everylabel@\eat@
\let\everyref\relax
\def\newlabel{\bgroup\everylabel\newlabel@}
\def\newlabel@#1#2#3{\if\@undefined\r@-#1\@\else\wrn@label{#1}{#3}\fi
  {\let\protect\noexpand\@Nxdef\r@-#1\@{#2}}\egroup}
\def\w@ref{\bgroup\everyref\w@@ref}
\def\w@@ref#1#2#3#4{%
  \if\@undefined\r@-#1\@{\bftext??}#2{#1}{}\else%
   \@xHref{#1}{\@DO{\expandafter\expandafter#3}\r@-#1\@\@nil}\fi
  #4{#1}{}\egroup}
\def\@@@xref#1{\w@ref{#1}\wrn@ref\@car\wrn@reference}
\def\@xref#1{\rom{\@@@xref{#1}}}
\let\xref\@xref
\def\pageref#1{\w@ref{#1}\wrn@ref\@cdr\wrn@reference}
\def\thepage{\ifnum\pageno<\z@\romannumeral-\pageno\else\number\pageno\fi}
\def\label@{\@bsphack\bgroup\everylabel\label@@}
\def\label@@#1#2{\everylabel@{{#1}{#2}}%
  \@labelref{#2}%
  \let\thepage\relax
  \def\protect{\noexpand\noexpand\noexpand}%
  \edef\@tempa{\edef\noexpand\@lastlabel{#1}%
    \W@X{\string\newlabel{#2}{{\@lastmark}{\thepage}}{\the\inputlineno}}}%
  \expandafter\egroup\@tempa\@esphack}
\def\label#1{\label@{#1}{#1}}
\def\fn@P@{\relaxnext@
    \DN@{\ifx[\next\DN@[####1]{}\else
        \ifx"\next\DN@"####1"{}\else\DN@{}\fi\fi\next@}%
    \FN@\next@}
\def\eat@fn#1{\ifx#1[\expandafter\eat@br\else
  \ifx#1"\expandafter\expandafter\expandafter\eat@qu\fi\fi}
\def\eat@br#1]#2{}
\def\eat@qu#1"#2{}
{\catcode`\~\active\lccode`\~`\@
\lowercase{\global\let\@@P@~\gdef~{\protect\@@P@}}}
\def\Protect@@#1{\def#1{\protect#1}}
\def\disable@special{\let\W@X@\eat@iii\let\label\eat@
    \def\footnotemark{\protect\fn@P@}%
  \let\footnotetext\eat@fn\let\footnote\eat@fn
    \let\refcounter\eat@\let\savecounter\eat@\let\restorecounter\eat@
    \let\advancecounter\eat@ii\let\setcounter\eat@ii
  \let\ifvmode\iffalse\Protect@@\@@@xref\Protect@@\pageref\Protect@@\nofrills
    \Protect@@\\\Protect@@~}
\let\notoctext\identity@
\def\W@X@#1#2#3{\@bsphack{\disable@special\let\notoctext\eat@
    \def\chapter{\protect\chapter@toc}\let\thepage\relax
    \def\protect{\noexpand\noexpand\noexpand}#1%
  \edef\next@{\if\@undefined#2\@\else\write#2{#3}\fi}\expandafter}\next@
    \@esphack}
\newcount\cnt@toc
\def\writeauxline#1#2#3{\W@X@{\cntref@{toc}\let\tocref\@HR}
  \@Xout{\string\@Xline{#1}{#2}{#3}{\thepage}}}
{\let\newwrite\relax
\gdef\@openin#1{\make@letter\@input{\jobname.#1}\t@mpcat}
\gdef\@openout#1{\global\expandafter\newwrite\csname tf@-#1\endcsname
   \immediate\openout\@N\tf@-#1\@\jobname.#1\relax}}
\def\@@openout#1{\@openout{#1}%
  \@@addto\readaux@{\immediate\closeout\@N\tf@-#1\@}}
\def\auxlinedef#1{\@Ndef\do@-#1\@}
\def\@Xline#1{\if\@undefined\do@-#1\@\expandafter\eat@iii\else
    \@DO\expandafter\do@-#1\@\fi}
\def\beginW@{\bgroup\def\do##1{\catcode`##112 }\dospecials\do\@\do\"
    \catcode`\{\@ne\catcode`\}\tw@\immediate\write\@N}
\def\endW@toc#1#2#3{{\string\tocline{#1}{#2\string\page{#3}}}\egroup}
\def\do@tocline#1{%
    \if\@undefined\tf@-#1\@\expandafter\eat@iii\else
        \beginW@\tf@-#1\@\expandafter\endW@toc\fi
}
\auxlinedef{toc}{\do@tocline{toc}}

\let\protect\empty
\def\Protect#1{\if\@undefined#1@P@\@\PROTECT#1\else\wrn@Protect#1\empty\fi}
\def\PROTECT#1{\@Nlet#1@P@\@#1\edef#1{\noexpand\protect\@Nx#1@P@\@}}
\def\pdef#1{\edef#1{\noexpand\protect\@Nx#1@P@\@}\@Ndef#1@P@\@}

\Protect\operatorname
\Protect\operatornamewithlimits
\Protect\qopname@
\Protect\qopnamewl@
\Protect\text
\Protect\topsmash
\Protect\botsmash
\Protect\smash
\Protect\widetilde
\Protect\widehat
\Protect\thetag
\Protect\therosteritem
\Protect\Cal
\Protect\Bbb
\Protect\bold
\Protect\slanted
\Protect\roman
\Protect\italic
\Protect\boldkey
\Protect\boldsymbol
\Protect\frak
\Protect\goth
\Protect\dots
\Protect\cong
\Protect\lbrace \let\{\lbrace
\Protect\rbrace \let\}\rbrace
\let\root@P@@\root \def\root@P@#1{\root@P@@#1\of}
\def\root#1\of{\protect\root@P@{#1}}

\def\frills{\ignorespaces\@txtopt@}
\def\frillsnotempty#1{\x@notempty{\@txtopt@{#1}}}
\def\numberline{\@numopt@}
\newif\if@theorem
\let\@therstyle\eat@
\def\@headtext@#1#2{{\disable@special\let\protect\noexpand
    \def\chapter{\protect\chapter@rh}%
    \edef\next@{\noexpand\@frills@\noexpand#1{#2}}\expandafter}\next@}
\let\AmSrighthead@\rightheadtext
\def\rightheadtext{\checkfrills@{\@headtext@\AmSrighthead@}}
\let\AmSlefthead@\leftheadtext
\def\leftheadtext{\checkfrills@{\@headtext@\AmSlefthead@}}
\def\@head@@#1#2#3#4#5{\@Name\pre\eat@bs#1\@\if@theorem\else
    \@frills@{\csname\expandafter\eat@iv\string#4\endcsname}\relax
        \ifx\protect\empty\@N#1\F@NT\@\fi\fi
    \@N#1\ST@LE\@{\counter#3}{#5}%
  \if@write\writeauxline{toc}{\eat@bs#1}{#2{\counter#3}\@HR{#5}}\fi
    \if@theorem\else\expandafter#4\fi
    \ifx#4\endhead\ifx\@txtopt@\identity@\else
        \headmark{\@N#1\ST@LE\@{\counter#3}{\frills\empty}}\fi\fi
    \@Name\post\eat@bs#1\@\ignorespaces}
\ifx\undefined\endhead\Invalid@\endhead\fi
\def\@head@#1{\checkstar@{\checkfrills@{\checkbrack@{\@head@@#1}}}}
\def\@thm@@#1#2#3{\@Name\pre\eat@bs#1\@
    \@frills@{\csname\expandafter\eat@iv\string#3\endcsname}
    {\@theoremtrue\check@therstyle{\@N#1\ST@LE\@}\frills
            {\counter#2}\@theoremfalse}%
    \@DO\envir@stack\end\eat@bs#1\@
    \@N#1\F@NT\@\@Name\post\eat@bs#1\@\ignorespaces}
\def\@thm@#1{\checkstar@{\checkfrills@{\checkbrack@{\@thm@@#1}}}}
\def\@capt@@#1#2#3#4#5\endcaption{\bgroup
    \edef\@tempb{\global\footmarkcount@\the\footmarkcount@
    \global\@N#2\N@M\@\the\@N#2\N@M\@}%
    \def\shortcaption##1{\global\def\sh@rtt@xt####1{##1}}\let\sh@rtt@xt\identity@
    \DN@{#4{\@tempb\@N#1\ST@LE\@{\counter#2}}}%
    \if\notempty{#5}\DNii@{\next@\@N#1\F@NT\@}\else\let\nextii@\next@\fi
    \nextii@#5\endcaption
  \if@write\writeauxline{#3}{\eat@bs#1}{{} \@HR{\@N#1\ST@LE\@{\counter#2}%
    \if\notempty{#5}.\enspace\fi\sh@rtt@xt{#5}}}\fi
  \global\let\sh@rtt@xt\undefined\egroup}
\def\@capt@#1{\checkstar@{\checkfrills@{\checkbrack@{\@capt@@#1}}}}
\let\captiontextfont@\empty

\ifx\undefined\subsubheadfont@\def\subsubheadfont@{\it}\fi
\ifx\undefined\proclaimfont\def\proclaimfont{\sl}\fi
\ifx\undefined\proclaimfont@\let\proclaimfont@\proclaimfont\fi
\def\proclaimfont{\proclaimfont@}
\ifx\undefined\definitionfont@\def\AmSdeffont@{\rm}
    \else\let\AmSdeffont@\definitionfont@\fi
\ifx\undefined\remarkfont@\def\remarkfont@{\rm}\fi

\def\newfont@def#1#2{\if\@undefined#1\F@NT\@
    \@Nxdef#1\F@NT\@{\@Nx.\expandafter\eat@iv\string#2\F@NT\@}\fi}
\def\newhead@#1#2#3#4{{%
    \gdef#1{\@therstyle\@therstyle\@head@{#1#2#3#4}}\newfont@def#1#4%
    \if\@undefined#1\ST@LE\@\@Ngdef#1\ST@LE\@{\headstyle}\fi
    \if\@undefined#2\@\gdef#2{\headtocstyle}\fi
  \@@addto\moretocdefs@{\\#1#1#4}}}
\outer\def\newhead#1{\checkbrack@{\expandafter\newhead@\expandafter
    #1\@txtopt@\headtocstyle}}
\outer\def\newtheorem#1#2#3#4{{%
    \gdef#2{\@thm@{#2#3#4}}\newfont@def#2#4%
    \@Nxdef\end\eat@bs#2\@{\noexpand\revert@envir
        \@Nx\end\eat@bs#2\@\noexpand#4}%
  \if\@undefined#2\ST@LE\@\@Ngdef#2\ST@LE\@{\proclaimstyle{#1}}\fi}}%
\outer\def\newcaption#1#2#3#4#5{{\let#2\relax
  \edef\@tempa{\gdef#2####1\@Nx\end\eat@bs#2\@}%
    \@tempa{\@capt@{#2#3{#4}#5}##1\endcaption}\newfont@def#2\endcaptiontext%
  \if\@undefined#2\ST@LE\@\@Ngdef#2\ST@LE\@{\captionstyle{#1}}\fi
  \@@addto\moretocdefs@{\\#2#2\endcaption}\newtoc{#4}}}
{
\outer\gdef\newtoc#1{{%
    \@DO\ifx\do@-#1\@\relax
    \global\auxlinedef{#1}{\do@tocline{#1}}{}%
    \@@addto\tocsections@{\make@toc{#1}{}}\fi}}}

\toks@\expandafter{\itembox@}
\toks@@{\bgroup\let\therosteritem\identity@\let\rm\empty
  \let\@Href\eat@\let\@Hname\eat@
  \edef\next@{\edef\noexpand\@lastmark{\therosteritem@}}\donext@}
\edef\itembox@{\the\toks@@\the\toks@}
\def\firstitem@false{\let\iffirstitem@\iffalse
    \global\let\lastlabel\@lastlabel}

\let\rosteritemrefform\therosteritem
\let\rosteritemrefseparator\empty
\def\rosteritemref#1{\hbox{\rosteritemrefform{\@@@xref{#1}}}}
\def\local#1{\label@\@lastlabel{\lastlabel-i#1}}

\def\xRef@P@{\gdef\lastlabel}
\def\xRef#1{\@xref{#1}\protect\xRef@P@{#1}}

\def\iref@P@{\gdef\lastref}
\def\itemref#1#2{\rosteritemref{#1-i#2}\protect\iref@P@{#1}}
\def\iref#1{\@xref{#1}\rosteritemrefseparator\itemref{#1}}

\def\eqtag{\tag\counter\equation}
\def\eqref#1{\thetag{\@@@xref{#1}}}
\def\tagform@#1{\ifmmode\hbox{\rm\else\rom{\fi
        (\ignorespaces#1\unskip)\iftrue}\else}\fi}

\let\AmSfnote@\makefootnote@
\def\makefootnote@#1{\bgroup\let\footmarkform@\identity@
  \edef\next@{\edef\noexpand\@lastmark{#1}}\donext@\AmSfnote@{#1}}

\def\clearpage{\ifnum\insertpenalties>0\line{}\fi\vfill\supereject}

\def\proof{\checkfrills@{\checkbrack@{%
    \check@therstyle{\@frills@{\demo}{\frills{Proof}}{}}
        {\frills{}\envir@stack\endremark\envir@stack\enddemo}%
  \envir@stack\endproof\ignorespaces}}}
\def\everyendproof{\qed}
\def\endproof{\nofrillscheck{\frills@{\everyendproof}\revert@envir\endproof\enddemo}}

\let\AmSref\ref
\let\AmSrefstyle\refstyle
\let\plaincite\cite
\def\citei@#1,{\citeii@#1\eat@,}
\def\citeii@#1\eat@{\w@ref{#1}\wrn@cite\@car\wrn@citation}
\def\mcite@#1;{\plaincite{\citei@#1\eat@,\unskip}\mcite@i}
\def\mcite@i#1;{\DN@{#1}\ifx\next@\endmcite@
  \else, \plaincite{\citei@#1\eat@,\unskip}\expandafter\mcite@i\fi}
\def\endmcite@{\endmcite@}
\def\cite#1{\mcite@#1;\endmcite@;}
\PROTECT\cite
\def\refstyle#1{\AmSrefstyle{#1}\uppercase{%
    \ifx#1A\relax \def\@ref@##1{\AmSref\xdef\@lastmark{##1}\key##1}%
    \else\ifx#1C\relax \def\@ref@##1{\AmSref\no\counter\refno}%
        \else\def\@ref@{\AmSref}\fi\fi}}
\refstyle A
\newcounter\refno\null
\newif\ifRefs
\gdef\Refs{\checkstar@{\checkbrack@{\csname AmSRefs\endcsname
  \nofrills{\frills{References}%
  \if@write\writeauxline{toc}{vartocline}{\@HR{\frills{References}}}\fi}%
  \def\ref{\@ref@}\Refstrue\ignorespaces}}}
\let\ref\xref

\newif\iftoc
\pdef\tocbreak{\iftoc\hfil\break\fi}
\def\tocsections@{\make@toc{toc}{}}
\let\moretocdefs@\empty
\def\newtocline@#1#2#3{%
  \edef#1{\def\@Nx#2line\@####1{\@Nx.\expandafter\eat@iv
        \string#3\@####1\noexpand#3}}%
  \@Nedef\no\eat@bs#1\@{\let\@Nx#2line\@\noexpand\eat@}%
    \@N\no\eat@bs#1\@}
\def\MakeToc#1{\@@openout{#1}}
\def\newtocline#1#2#3{\Err@{\Invalid@@\string\newtocline}}
\def\make@toc#1#2{\penaltyandskip@{-200}\aboveheadskip
    \if\notempty{#2}
        \centerline{\headfont@\ignorespaces#2\unskip}\nobreak
    \vskip\belowheadskip \fi
    \@openin{#1}\relax
    \vskip\z@}
\def\contents{\readaux\checkfrills@{\checkbrack@{\@contents@}}}
\def\@contents@{\toc@{\frills{Contents}}\envir@stack\endcontents%
    \def\nopagenumbers{\let\page\eat@}\let\newtocline\newtocline@\toctrue
  \def\@HR{\Href@{toc}}%
  \def\tocline##1{\csname##1line\endcsname}
  \edef\caption##1\endcaption{\expandafter\noexpand
    \csname head\endcsname##1\noexpand\endhead}%
    \ifmonograph@\def\vartoclineline{\Chapterline}%
        \else\def\vartoclineline##1{\sectionline{{} ##1}}\fi
  \let\\\newtocline@\moretocdefs@
    \ifx\@frills@\identity@\def\\##1##2##3{##1}\moretocdefs@
        \else\let\tocsections@\relax\fi
    \def\\{\unskip\space\ignorespaces}\let\maketoc\make@toc}
\def\endcontents{\tocsections@\vskip-\lastskip\revert@envir\endcontents
    \endtoc}

\if\@undefined\selectf@nt\@\let\selectf@nt\identity@\fi
\def\Err@math#1{\Err@{Use \string#1\space only in text}}
\def\textonlyfont@#1#2{%
    \def#1{\RIfM@\Err@math#1\else\edef\f@ntsh@pe{\string#1}\selectf@nt#2\fi}%
    \PROTECT#1}
\tenpoint

\def\newshapeswitch#1#2{\gdef#1{\selectsh@pe#1#2}\PROTECT#1}
\def\shapeswitch#1#2#3{\@Ngdef#1\string#2\@{#3}}
\shapeswitch\rm\bf\bf
\shapeswitch\rm\tt\tt
\shapeswitch\rm\smc\smc
\newshapeswitch\em\it
\shapeswitch\em\it\rm
\shapeswitch\em\sl\rm
\def\selectsh@pe#1#2{\relax\if\@undefined#1\f@ntsh@pe\@#2\else
    \@N#1\f@ntsh@pe\@\fi}

\def\@itcorr@{\leavevmode
    \edef\prevskip@{\ifdim\lastskip=\z@ \else\hskip\the\lastskip\relax\fi}\unskip
    \edef\prevpenalty@{\ifnum\lastpenalty=\z@ \else
        \penalty\the\lastpenalty\relax\fi}\unpenalty
    \/\prevpenalty@\prevskip@}
\def\rom@P@#1{\@itcorr@{\selectsh@pe\rm\rm#1}}
\def\rom{\protect\rom@P@}
\def\Rom@P@#1{\@itcorr@{\rm#1}}
\def\Rom{\protect\Rom@P@}
{\catcode`\-11
\HyperRefs{idx}
\HyperRefs{glo}
\newcount\cnt@idx \global\cnt@idx=10000
\newcount\cnt@glo \global\cnt@glo=10000
\gdef\writeindex#1{\W@X@{\cntref@{idx}}\tf@-idx
 {\string\indexentry{#1}{\Hlast@{idx}}{\thepage}}}
\gdef\writeglossary#1{\W@X@{\cntref@{glo}}\tf@-glo
 {\string\glossaryentry{#1}{\Hlast@{glo}}{\thepage}}}
}
\def\emph#1{\@itcorr@\bgroup\em\ignorespaces#1\unskip\egroup
  \DN@{\DN@{}\ifx\next.\else\ifx\next,\else\DN@{\/}\fi\fi\next@}\FN@\next@}
\def\makequoteactive{\catcode`\"\active}
{\makequoteactive\gdef"{\FN@\quote@}
\gdef\quote@{\ifx"\next\DN@"##1""{\quoteii{##1}}\else\DN@##1"{\quotei{##1}}\fi\next@}}
\let\quotei\eat@
\let\quoteii\eat@
\def\MakeIndex{\@openout{idx}}
\def\MakeGlossary{\@openout{glo}}

\def\endofpar#1{\ifmmode\ifinner\endofpar@{#1}\else\eqno{#1}\fi
    \else\leavevmode\endofpar@{#1}\fi}
\def\endofpar@#1{\unskip\penalty\z@\null\hfil\hbox{#1}\hfilneg\penalty\@M}

\newdimen\normalparindent\normalparindent\parindent
\def\firstparindent#1{\everypar\expandafter{\the\everypar
  \global\parindent\normalparindent\global\everypar{}}\parindent#1\relax}

\@@addto\disablepreambule@cs{%
    \\\readaux\relax
    \\\begin\relax
    \\\readaux@\relax
    \\\@openout\eat@
    \\\@@openout\eat@
    \/\Monograph\empty
    \/\MakeIndex\empty
    \/\MakeGlossary\empty
    \/\MakeToc\eat@
    \/\HyperRefs\eat@
    \/\NoHyperRefs\eat@
}

\csname label.def\endcsname


\def\punct#1#2{\if\notempty{#2}#1\fi}
\def\sppunct{\punct{.\enspace}}
\def\varpunct#1#2{\if\frillsnotempty{#2}#1\fi}

\def\headstyle#1#2{\numberline{#1\sppunct{#2}}\ignorespaces#2\unskip}
\def\headtocstyle#1#2{\numberline{#1\punct.{#2}}\space #2}

\def\specialtocstyle#1#2{#2}
\newcounter\section\null
\newcounter\subsection\section
\newcounter\subsubsection\subsection
\newhead\specialsection[\specialtocstyle]\null\endspecialhead
\newhead\section\section\endhead
\newhead\subsection\subsection\endsubhead
\newhead\subsubsection\subsubsection\endsubsubhead
\def\firstappendix{\global\sectionno0 %
  \counterstyle\section{\Alphnum\sectionno}%
    \global\let\firstappendix\empty}

\def\appendixtocstyle#1#2{\space\numberline{Appendix #1\sppunct{#2}}#2}
\newhead\appendix[\appendixtocstyle]\section\endhead

\let\endAmSdef\enddefinition
\def\proclaimstyle#1#2{\numberline{#2\varpunct{.\enspace}{#1}}\frills{#1}}
\copycounter\thm\subsubsection
\theorem\thm\endproclaim
\proposition\thm\endproclaim
\lemma\thm\endproclaim
\corollary\thm\endproclaim
\definition\thm\endAmSdef
\example\thm\endAmSdef
\examples\thm\endAmSdef

\def\captionstyle#1#2{\frills{#1}\numberline{\varpunct{ }{#1}#2}}
\newcounter\figure\null
\newcounter\table\null
\newcaption{Figure}\figure\figure{lof}\botcaption
\newcaption{Table}\table\table{lot}\topcaption

\copycounter\equation\subsubsection

\expandafter\ifx\csname label.def\endcsname\relax\input label.def \fi
\def\stydate{June 26, 2000}
\def\styname{DEBUG.DEF}
\immediate\write16{This is \styname\space by A.Degtyarev <\stydate>}
\expandafter\ifx\csname debug.def\endcsname\relax\else
  \message{[already loaded]}\endinput\fi
\expandafter\edef\csname debug.def\endcsname{%
  \catcode`\noexpand\@\the\catcode`\@\edef\noexpand\styname{\styname}
  \def\expandafter\noexpand\csname debug.def\endcsname{\stydate}}
\catcode`\@=11
{\edef\temp{\the\everyjob\W@{\styname: <\stydate>}}
\global\everyjob\expandafter{\temp}}

\def\n@te#1#2{\leavevmode\vadjust{%
 {\setbox\z@\hbox to\z@{\strut\eightpoint\let\quotei\filename#1}%
  \setbox\z@\hbox{\raise\dp\strutbox\box\z@}\ht\z@=\z@\dp\z@=\z@%
  #2\box\z@}}}
\def\leftnote#1{\n@te{\hss#1\quad}{}}
\def\rightnote#1{\n@te{\quad\kern-\leftskip#1\hss}{\moveright\hsize}}
\def\?{\FN@\qumark}
\def\qumark{\ifx\next"\DN@"##1"{\leftnote{\rm##1}}\else
 \DN@{\leftnote{\rm??}}\fi{\rm??}\next@}
\def\filename#1{\hbox{\tt #1}}
\def\mnote@@#1{\rightnote{\vtop{%
 \ifcat\noexpand"\noexpand~\def"##1"{\filename{##1}}\fi
 \hsize2.0in \baselineskip7\p@\parindent\z@
 \tolerance\@M\spaceskip2.6\p@ plus10\p@ minus.9\p@\rm#1}}}
\def\mnote#1{\@bsphack\mnote@{#1}\@esphack}

\def\nonotes{\let\mnote@\eat@}
\def\printnotes{\let\mnote@\mnote@@}
\printnotes

\def\PrintLabels{%
 \gdef\printlabel@##1##2{\ifvmode\else\leftnote{\eighttt##2}\fi}}
\def\NoLabels{\global\let\printlabel@\eat@ii}
\NoLabels
\@@addparm\everylabel@{\printlabel@#1}

\def\PrintFiles{\gdef\outputmark@{\line{\hfill\smash{\raise1cm\vbox{%
  \hbox to\z@{\kern1cm\tenrm\the\month/\the\day/\the\year\hss}%
  \hbox to\z@{\kern1cm\tentt\jobname.tex\hss}%
  \ifx\filecomment\undefined\else\hbox to\z@{\kern1cm\tenrm\filecomment\hss}\fi}}}}}
\def\NoFiles{\global\let\outputmark@\empty}
\def\PageMark{\gdef\outputmark@}
\let\outputmark@\empty
\@@addto\pagetop@\outputmark@

\def\tracerefs{\def\wrn@reference##1##2{\W@X{\string\reference{##1}}}}
\def\tracecites{\printwarning\nocite
  \def\wrn@citation##1##2{\@Nxdef\cite@-##1\@{##1}\W@X{\string\citation{##1}}}}
\def\wrn@@nocite#1#2{\ifRefs\wrn@@@nocite{#2}\fi}
\def\wrn@@@nocite#1{\if\@undefined\cite@-#1\@\warning@{citation `#1' [\@lastmark] not used}\fi}
\let\wrn@nocite\eat@ii
\let\reference\eat@
\let\citation\eat@
\@@addparm\everylabel@{\wrn@nocite#1}

\@@addto\disablepreambule@cs{%
    \/\PrintFiles\empty
    \/\NoFiles\empty
    \/\PageMark\empty
}

\csname debug.def\endcsname
\input xypic

\MakeToc{toc}

\let\textsection=\S

\def\CC{\Cal C}
\def\conj{\operatorname{conj}}
\def\C{{\Bbb C}}
\def\R{{\Bbb R}}
\def\Z{{\Bbb Z}}

\def\P{{\Bbb P}}
\def\Rp#1{\R\roman P^{#1}}

\def\Im{\mathop{\roman{Im}}\nolimits}

\def\Ker{\mathop{\roman{Ker}}\nolimits}
\def\ind{\mathop{\roman{ind}}\nolimits}

\def\Cl{\mathop{\roman{Cl}}\nolimits}

\def\la{\langle}
\def\ra{\rangle}
\def\id{\mathop{\roman{id}}\nolimits}

\def\Arf{\mathop{\roman{Arf}}\nolimits}
\def\emptyset{\varnothing}
\def\oo{\varnothing}

\def\rank{\operatorname{rank}}

\def\D{\Delta}
\def\G{\Gamma}
\def\S{\Sigma}

\def\til{\widetilde}

\def\b{\beta}
\def\g{\gamma}

\def\l{\lambda}

\def\s{\sigma}
\def\ss{\frak s}

\def\Quadr{\mathop{\roman{Quadr}}\nolimits}
\def\sm{\smallsetminus}

\def\Fix{\operatorname{Fix}}
\def\Sym{\mathop{\roman{Sym}}\nolimits}
\def\Sing{\mathop{\roman{Sing}}\nolimits}
\def\mon{G_{\roman{mon}}}
\def\Map{\mathop{\roman{Map}}}

\def\trans{\mathop{\roman{\tau}}\nolimits}
\def\hP{\widehat P^2}
\def\hspec{\widehat\spec}
\def\res{\roman{res}}

\def\Gres{\G_{\roman{res}}}
\def\GII{\G_{I\!I}}
\def\GIII{V_{I\!I\!I}}
\def\spec{S}

\def\spect{\widetilde{\spec}}

\def\cs{c_{\spect}}
\def\TT{\Cal T}

\let\t=\trans

\let\ge\geqslant
\let\le\leqslant
\let\+\sqcup
\let\dsum\+


\NoBlackBoxes
\rightheadtext{Deformation classification of real cubic 3-folds with a line}

\topmatter
\title
Deformation classification of real non-singular cubic threefolds with a marked line
\endtitle
\author
S.~Finashin, V.~Kharlamov
\endauthor
\address Middle East Technical University,
Department of Mathematics\endgraf Damlupinar Bul. 1, Ankara 06800 Turkey
\endaddress
\address
Universit\'{e} de Strasbourg et IRMA (CNRS)\endgraf 7 rue Ren\'{e}-Descartes, 67084 Strasbourg Cedex, France
\endaddress

\abstract
We prove that the space of pairs $(X,l)$ formed by a real non-singular cubic hypersurface  $X\subset P^4$
with a real line $l\subset X$ has 18 connected components and
 give for them several quite explicit interpretations.
The first one relates these components to the orbits of the monodromy action
on the set of connected components of the Fano surface $F_\R(X)$ formed by real lines on $X$.
For another interpretation we associate with each of the 18 components a well defined real
deformation class of real non-singular plane quintic curves and show that this deformation class together with the
real deformation class of $X$ characterizes completely the component.
\endabstract

\contents
\endcontents

\endtopmatter
\hskip1.6in La derni\`ere chose qu'on trouve en faisant un ouvrage,

\hskip1.7in est de savoir celle qu'il faut mettre la premi\`ere.
\vskip3mm
\hskip2.5in  Blaise Pascal, {\it Pens\'ees }(1670), 19.
 \vskip6mm

\document

\section{Introduction}

\subsection{Sketch of principal results}
The study of lines on projective hypersurfaces had started from discovery of 27 complex lines by Cayley and Salmon. It followed by an analysis of reality of these lines performed
by  Schl\"afli, which in its turn motivated Klein's interpretation of 5 species of cubic surfaces introduced by Schl\"afli as a classification of non-singular real cubic surfaces $X$ up to deformation
(here and further on, by deformation we mean variations preserving non-singularity).
A major step in further analysis of the real aspects of the theory of cubic surfaces was made by B.~Segre. Among other results, he determined, for each of the 5 deformation classes, the group
of substitutions induced by the monodromy action on the set of real lines $l\subset X$, which implies, in particular,
the deformation classification of the  pairs $(X,l)$, where $X$ is a real non-singular cubic surface and $l\subset X$ is a real line.
For a modern view and contemporary further spectacular development, as well as a correction of one of Segre's results on monodromy,
we send the reader to  D.~Allcock, J.A.~Carlson, and D.~Toledo paper \cite{ACT}.

The real aspects of the theory of cubic projective hypersurfaces of higher dimensions are
less studied. Recent advance in the case of real cubic threefolds  $X\subset P^4$
is due to V.~Krasnov who gave their
deformation classification
(see  \cite{Kr-rigid})
along with a description of
how the topology of the real locus $F_\R(X)$ of the {\it Fano surface of $X$}
formed by the lines on $X$
depends on the deformation class of $X$ (see \cite{Kr-Fano}).

Here, we make the next step and perform
the deformation classification of the pairs $(X,l)$, where $X$ is a real non-singular cubic hypersurface in $P^4$  and $l\subset X$ is a real line.

Our interest to this problem has grown, indeed,
from a wish to disclose a puzzling coincidence of two graphs: the
adjacency graph of the deformation classes of
real non-singular cubic threefolds and
a similar graph for real non-singular plane quintic curves, see Figure  \ref{2graphs}.
The vertices of these graphs represent the deformation classes of real non-singular cubic threefolds $X$ at the top graph and those of real non-singular plane quintics at the bottom one. For greater clarity, the vertices are shown as circles containing inside an indication of the topological type of
the locus of real points
in the case of cubics and a schematic picture of position of
the locus of real points
in the real plane in the case of quintics.
Two vertices are connected by an edge if the corresponding deformation classes are adjacent through a wall, that is if these classes can be joined by a continuous family of
real cubic threefolds or, respectively,  real plane quintics, such that all but one element of the family is
non-singular, the singular one (cubic or quintic) has only one singular point, and this singular point is a simple double point.
\midinsert
\epsfbox{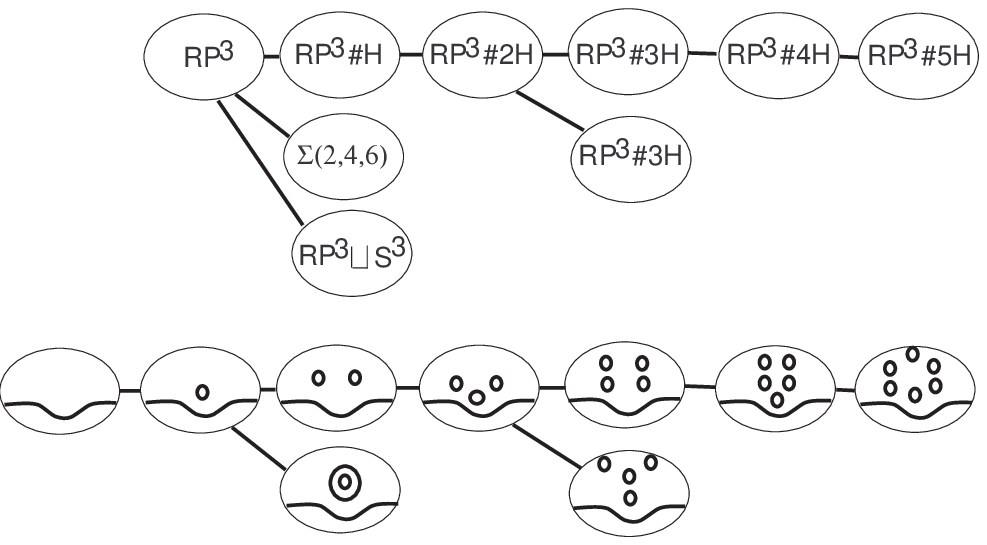}
\figure\label{2graphs}{
$H$ states for a handle,
$H=S^1\times S^2$, and $\Sigma(2,4,6)$ denotes a Seifert manifold with 3 multiple fibers of multiplicities 2, 4, 6;
the double occurrence of the type ${\Rp3}
\# 3H$ means that the real non-singular cubics with
this topological type of the real locus
form two distinct deformation classes.}
 \endfigure
\endinsert

Due to the classical correspondence between cubic threefolds and plane quintics that associates to
a cubic threefold $X$ a plane quintic $S$ (we call $S$
the {\it spectral quintic}) arising as the discriminant locus of the projection $X\to P^2$
with a chosen line $l\subset X$ as the center, it was natural to expect
 that such a coincidence should not be accidental.
However, two
 difficulties obstruct an immediate application of this correspondence. The first difficulty is in
degenerating of $S$ for some choices of $l$ even on a non-singular $X$. The second one is due to presence of several
connected components in the real locus $F_\R(X)$ of the Fano surface of $X$. For most of the deformation classes of $X$, this leads
to dependence of the deformation class of $S$ on the choice of a connected component of $F_\R(X)$ that contains $l$.

The first difficulty turns out to be inessential: we show that despite such degenerations
the deformation class $[S]$ of $S$ is not changing under deformations of pairs $(X,l)$. Furthermore, we  prove (see Theorem \ref{main}, parts (1-2))
that deformation classes of pairs $(X,l)$ are completely determined by
{\it matchings} $([X],[S])$ between the deformation classes $[X]$ and $[S]$.
The second difficulty happens to be a more subtle issue.

We present the final deformation classification in a form of a complete list of  possible matchings
in terms of simple topological invariants of $X$ and $S$, such as Smith discrepancy and Klein type,
see Theorem \ref{main} (part (3))
and Corollary  \ref{18classes}.
Accordingly, the number of deformation classes happens to be 18.

In fact, the set of matchings, and hence the set of deformation classes, turns
out to be in 1--1 correspondence with the set of orbits of
the monodromy action induced by real deformations
on the sets of connected components
of the real loci of Fano surfaces. In particular, each real Fano surface has a unique connected component whose Euler characteristic is odd,
this component is preserved by the monodromy,
and, as we show, it is the choice of $l$ on these components that
establishes a natural isomorphism between the mentioned above two graphs (see Theorem \ref{N-correspondence}). The monodromy on
the set of  toric components is less trivial, and we give a full explicit description of these orbits in  Theorem \ref{T-orbits}.

For example,
if $X_\R= \Rp3\# 5H$ (the case of so-called maximal cubic threefolds), then
$F_\R(X)$
is formed by a non-orientable component $N_5=\#_5\Rp2$
(of Euler characteristic $-3$)
and 15 torus-components, and, according to Theorem \ref{T-orbits}, the
monodromy action on these tori has two orbits: 6 torus components in one and 9 components in another orbit;
namely, a choice
of a line $l\subset X$ on a torus gives an $(M-2)$-quintic, which is of Klein type I for tori from the first orbit,
and of type II for the second orbit.

\subsection{Main tools and constructions}
The key tools involved traditionally in studying complex cubic threefolds $X$
are Fano surfaces, Jacobian and Prym  varieties, conic bundles that project
$X$ from lines $l\subset X$ to $\P^2$,
and discriminant loci
of these projections (spectral quintics) equipped with odd theta-characteristics.
Over the reals they are not well studied (if at all),
and our first task was to develop them accordingly
disclosing
an interaction between geometrical and topological properties of these tools
over the reals.

In particular, it is the interplay
between the arithmetics of the action of the complex conjugation
on the homology of the cubic threefolds and that of the action on the homology of the associated plane quintics that helped us
to give a complete effective description of the relations between Klein types and Smith discrepancies of a
 real cubic threefold $X$ and its spectral quintic $S$ (see  Lemma \ref{Klein-type-relation}  and Corollary \ref{discrepancy-relation}),
 and as a consequence, allowed us to obtain a strong lower bound on the number of matchings.
 This description, which plays a crucial role in our paper, is achieved by developing a technique that discloses the above relations
 in arithmetic terms via the homological representative of a theta-characteristic (Subsections \ref{HC}--\ref{Klein-types})
 and, alternatively, in geometric terms via the mutual
 position of the theta-conic and the ovals of the quintic (Subsections \ref{OB}--\ref{surjectivity of SMC}).

In showing that it is a sharp bound,
it is the geometric methods that prevail. For example, to
determine the orbits of the monodromy action on the set of real components of Fano surfaces of maximal cubic threefolds $X$ we
make use of their degeneration
to remarkable  6-nodal cubic threefolds $X_0$ studied by C.~Segre (see Subsection \ref{Proof-T-orbits-for-M-cubics}).
Such a degeneration establishes
a 1--1 correspondence between the 15 torus components of $F_\R(X)$ we are interested in, the 15 lines connecting pairwise the 6 nodes of $X_0$,
and the 15 real lines on some real $(M-1)$-cubic surface $Y$ that can be associated to $F(X_0)$.
To conclude, we make appeal to the fundamental bipartition of these 15 lines in two groups (due to B.~Segre): 9 hyperbolic and 6 elliptic lines.

In addition to extending classical tools to the real setting, we develop also certain tools looking a somewhat unconventional,
like, for example, the {\it conormal projection} of Fano surface $F(X)$ to $\P^2$ (see Section \ref{conormal-section}).
The latter is crucial in out study of what we call exotic cubics (the case $X_\R=\Sigma(2,4,6)$).
In this regard, let us note also that, as it is underlined in \ref{conormal-remark},  this new tool is not only
 crucial in our study of real cubic threefolds but also allows us to
 disclose new geometric properties of real plane quintics.

\subsection{Structure of the paper} Our proof of the main results splits into two separate parts. We start with a lower bound
for the number of deformation classes in terms of matchings (Sections 2 - 5), and then respond to the question of how many deformation classes correspond
to each of the matchings (Sections 6 - 9). The techniques used in different parts are of different nature, and we develop them as we go along.

In Section 2 we recall some basic constructions associated to a cubic threefold $X$ with a chosen line $l$,
like those of the conic bundle, the spectral quintic $S\subset P^2$ and the
theta-conic $\Theta$.
We review their principal properties (which are mostly well-known)
and adapt them to our real setting.

In Section 3, we interpret
theta-conics as odd non-degenerate theta-characteristics, relate them to spectral coverings, and elaborate the {\it real spectral correspondence}
which is a real version of classical White's theorem that ensures a recovery of
a pair $(X,l)$ from a pair $(S,\Theta)$.
Moreover, we prove the invariance of the deformation class of $S$
under deformations of the pair $(X,l)$.

In Section 4 we review
the deformation classification of real cubic threefolds and real plane quintics, and in particular, discuss
the topological invariants involved in these classifications (Smith discrepancy and Klein's type).
We review also basic facts
on the Fano surfaces $F(X)$ and, in particular, recall the description of their real loci $F_\R(X)$.
Finally, we introduce the key notion of {\it spectral matching} and formulate our main results:
Theorem  \ref{main}, that describes the deformation classification of pairs $(X,l)$ in terms of matchings, Theorem \ref{N-correspondence}
relating the Smith discrepancy and Klein's type under the spectral matching
correspondence, and Theorem \ref{T-orbits} describing the orbits of the monodromy action on the set of torus components of Fano surfaces $F_\R(X)$.

In Section 5 we compare
the Smith discrepancies and Klein's types of $X$ and $S$
and prove part (3) of Theorem  \ref{main}. This gives
a lower bound to the number of deformation classes of pairs $(X,l)$.

In Section 6 we treat a special kind of real cubic threefolds. These cubics, that we call {\it exotic} are
exceptional in several respects.
For instance, topology of $X_\R$ is more sophisticated for them
than for the other types of real non-singular cubic threefolds, namely, it
is diffeomorphic to a
Seifert manifold $\Sigma(2,4,6)$,
and the Fano surface $F_\R(X)$ has two non-orientable real components
instead of one (see Section \ref{SMC}). We end Section 6 by proving Theorem \ref{N-correspondence} for exotic threefolds.
This is done through the conormal projection
of $F(X)$ to the projective plane.

In Section 7 we analyze real nodal cubic threefolds $X_0$ and study their
{\it quadrocubics}, that is the curves of bidegree $(2,3)$ in $P^3$ associated to each of the nodes of $X_0$.
In particular, we express in terms of quadrocubics
some topological and monodromy properties of cubic threefolds $X$ obtained by perturbation of $X_0$.
We relate quadrocubics to the spectral curves $S_0$ of $X_0$ and analyze how
the latter are perturbed into
the spectral curves $S$ of $X$.

In Section 8 we study the Fano surfaces $F(X_0)$ of
real nodal cubic threefolds $X_0$ and relate them to the symmetric square of quadrocubics.
Our key example is the case of real 6-nodal Segre cubic threefolds $X_0$, in which we find
the monodromy action on the nodes of $X_0$.
This allows us to find the monodromy of torus components of $F_\R(X)$, where $X$ is
the maximal cubic threefold obtained by perturbation of $X_0$.

In Section 9 we complete the proof of the main results: Theorems  \ref{main}, \ref{N-correspondence}, and \ref{T-orbits}.

Section 10 contains a few additional remarks and applications: on the spectral correspondence in terms of the moduli spaces
(in a setting a bit more general that we use in our paper), on the permutation groups realized by the monodromy action on the
set of real components of the Fano surfaces,
and on some
topological and geometrical information that one can obtain about real non-singular quintics using methods developed.

\subsection{Conventions}\label{conventions}
In this paper algebraic varieties by default are complex, for example, $P^n$ stands for the complex projective space.
Such a variety $X\subset P^n$ is {\it real} if it is invariant under the complex conjugation in $P^n$.
In the case of abstract complex varieties $X$ not embedded in
$P^n$ a {\it real structure} on $X$ is usually defined as an anti-holomorphic
involution $\conj\:X\to X$, but in this work we deal mostly
with subvarieties of $P^n$ or
other Grassmannians, and the
involution $\conj$ on $X$ is then induced from there.

We denote by $X_\R$ the fixed-point-set of the real structure, $X_\R=\Fix(\conj)$, and call it the {\it real locus of $X$}.
Speaking on non-singular real varieties, we mean that the whole $X$ (not only $X_\R$) has no singular points.

We prefer  to use classical rather than scheme-theoretic viewpoint and terminology,
allowing, for example,
to use, whenever possible, the same notation for polynomials
and for varieties that they represent, as well as the same notation both for the variety and its complex point set.

By a {\it nodal singularity }(a node) we mean
a singularity of type $A_1$, that is a simple
double point. A {\it nodal variety} (curve, surface, threefold) is
a one whose singular locus contains only nodes.

Recall that the space $\Cal C^{d,n}$ formed by hypersurfaces $Z\subset P^{n+1}$ of degree $d$ is a projective space
of dimension $\binom{d+n+1}{n+1}$
where the singular $Z$ form a so-called
{\it discriminant hypersurface} $\Delta^{d,n}\subset \Cal C^{d,n}$.
The real locus $\Cal C^{d,n}_\R$ is constituted by real hypersurfaces $Z$.
The connected components of the complement $\Cal C^{d,n}_\R\sm\Delta^{d,n}_\R$ are called {\it real deformation classes}:
they can be viewed as {\it chambers}
separated by {\it walls}, that is the connected components of the space
of one-nodal real hypersurfaces (the principal stratum of $\Delta^{d,n}_\R$).
The mutual position of the chambers can be characterized by
the  {\it adjacency graph} $\G_{d,n}$ whose vertex set $V_{d,n}$ is the set of
chambers and an edge between two vertices means that the corresponding chambers lie on the opposite sides of some wall
(note that $\G_{d,n}$ indicates only existence of walls and not their number).
We have already presented such graphs $\G_{3,3}$ and $\G_{5,1}$ on Figure 1
(and in this paper we are not concerned about any other values of $d$ and $n$).

In accordance with the above definition, by a {\it real deformation}
of non-singular hypersurfaces we mean simply a
continuous
path $[0,1]\to\Cal C^{d,n}_\R\sm\Delta^{d,n}_\R$, $t\mapsto Z_t$.
By a {\it real equisingular deformation} of a nodal hypersurface $Z$
we mean a
continuous
path in the stratum of $\Delta^{d,n}_\R$ parameterizing hypersurfaces $Z$ with a fixed number of nodes.

A {\it nodal degeneration} of a non-singular hypersurface
$Z$ is a continuous family
$Z_t$, $0\le t\le 1$, $Z=Z_1$, of
varieties which are non-singular for $t>0$ and
nodal for $t=0$.
The same family $Z_t$ can be called also  a {\it real perturbation} of a real nodal variety $Z_0$.

Curves of degree $6$ in $P^3$ which are
complete intersections of a cubic hypersurface with a quadric will be called {\it quadrocubics}.
By a conic we mean a plane curve of degree 2.

Whenever an equivalence relation is clear from the context, we use notation $[X]$ for the equivalence class of $X$.

If the most of the content of a section, or a subsection, is well known, we put at the beginning of it a reference to sources containing the corresponding results
and provide proofs only for those ones for which we did not find an appropriate reference.

\subsection{Acknowledgements}
This work was done during visits of the first author to
the University of Strasbourg and joint visits
to  the Max Planck Institute in Bonn
and the  Mathematisches Forschungsinstitut Oberwolfach.
We thank all these institutions for providing excellent working conditions.


\section{Spectral quintics and theta-conics of cubic threefold}\label{Spec_generalities}

\subsection{Conic bundles: generalities}\label{conic-bundles}
By a {\it conic bundle} we mean
a nonsingular irreducible projective variety $X$ equipped with a regular map $\pi\:X\to P^2$
such that each scheme theoretic fiber is isomorphic to a plane conic.

\proposition\label{nodality}\cite{Beau}
Let $\pi\:X\to P^2$ be a conic bundle. Then:
\roster
\item There exists a rank 3 vector bundle $E$ over $P^2$ and an integer $k$ such that the conic bundle is defined in the associated projective
 fibration
$P (E)$ by a quadratic form $h\in H^0(P^2, Sym^2(E^*)\otimes \Cal O(k))$.
The equation $\det h =0$ defines either an empty set or a curve in $P^2$.
\item The curve $S\subset P^2$ defined by $\det h =0$ has at most nodes as singular points.
\item  At each point $s\in \Sing S$ there exist local coordinates $u,v$ such that
over a neighborhood $U$ of $s$ the conic bundle is isomorphic to a conic bundle given in $U\times P^2$ by equation $ux^2+vy^2+z^2=0$;
at each point $s\in S\sm \Sing S$ there exist local coordinates $u,v$ such that over a neighborhood $U$ of $s$ the local model in $U\times P^2$ is given by equation $ux^2+y^2+z^2=0$;
while for $s\in P^2\sm S$ there exist local models given by equation $x^2+y^2+z^2=0$.
In particular, the fibers $X_s=\pi^{-1}(s)\subset P(E_s)$ are non-singular conics for $s\in P^2\sm S$, reduced singular conics for $s\in S$,  and
double lines for $s\in \Sing S$. \quad \quad \qed
\endroster
\endproposition

The curve $S$ described in this Proposition is called the {\it discriminant curve} of the conic bundle $\pi\:X\to P^2$.

Each conic bundle $\pi : X\to  P^2$
gives rise to a ramified double covering
$p_S\:\tilde S\to S$, where for each $s\in S$ the points of $p_S^{-1}(s)\subset\tilde S\subset  \operatorname{Gr}(E,1)$ are the lines
contained in the fibers $X_s\subset P(E_s)$.
This covering is unramified over $S\setminus \Sing S$, and it has a simple ramification over each branch at each node.

\subsection{The conic bundle for a line on a cubic threefold}\label{conic-bundles-for-cubics}
Consider a non-singular
cubic threefold $X\subset P^4$ and  a line $l\subset X$.
The net of planes $\{P_s\}_{s\in P^2}$ passing through $l$
induces a conic bundle $\pi_l\:X_l\to P^2$ whose fibers are {\it residual conics}
 $r_s$ traced by $X$ in the net of planes,
namely, $P_s\cap X=l+r_s$ (as a divisor). Here $X_l$ is obtained by blowing up $X$ along $l$.
The discriminant curve $\spec=\spec_l$ of this conic bundle will be called  the {\it spectral curve of $(X,l)$}.

Let us choose projective coordinates $(x,y,z,u,v)$ in $P^4$ so that $l$ is defined by equations $x=y=z=0$. Then, the equation of $X$ can be written as
$$
u^2 L_{11}+2uv L_{12}+v^2L_{22}+2uQ_{1}+2vQ_{2}+C=0,\eqtag\label{canon.equat}
$$
where $L_{ij}, Q_i, C$ stand for homogeneous polynomials in $x,y,z$ of degrees 1,2, and 3 respectively. Next,
the {\it residual conic} $r_s$, for each point $s=[x\!:\!y\!:\!z]\in P^2$, is
defined by the symmetric {\it fundamental matrix}
$$
A_s=\pmatrix
L_{11}&L_{12}&Q_{1}\\
L_{12}&L_{22}&Q_{2}\\
Q_{1}&Q_{2}&C
\endpmatrix
\eqtag\label{spec-matrix}
$$
and the spectral curve $\spec$ is defined by a degree 5 equation
$$
\det A_s=0.
\eqtag\label{spec-curve}
$$

The conic defined by equation
$$
\det\pmatrix
L_{11} & L_{12} \\
L_{12} & L_{22}
\endpmatrix
=0, \eqtag\label{theta-conic}
$$
will be called the {\it theta-conic} and denoted by $\Theta$.

Note that for the definition of the fundamental matrix non-singularity of the cubic $X$ is not required, but for certain degenerations of $X$
the conic $\Theta$ or/and the quintic $\spec$ may be not well-defined (if their equations identically vanish).
For non-singular or one-nodal $X$ this cannot happen:
for $\spec$ it follows, for instance, from \cite{Clemens-Griffiths, Lemma 8.1}.
For $\Theta$, it is straightforward; for nonsingular $X$ it follows, for example, from the next Lemma.

\lemma\label{l-nonsing}
The cubic threefold $X$ given by equation
$
u^2 L_{11}+2uv L_{12}+v^2L_{22}+2uQ_{1}+2vQ_{2}+C=0
$
is nonsingular along the line $l=\{x=y=z=0\}$ if and only if the quadratic form $Q_\Theta=L_{11}L_{22}-L_{12}^2$ is neither a
complete square nor zero. If $X$ contains a singular point outside $l$, then either $\det A_s$ is zero for every $s\in P^2$
or the curve $\det A_s=0$
has a singular point.
\endlemma

\proof
At the points of the line $l$ we have
$
df= u^2\, dL_{11} +2 uv \, dL_{12}+v^2\, dL_{22}
$.
If $df$ vanishes for some $(u,v)\ne(0,0)$, then by a projective $[u:v]$-coordinate change we can obtain vanishing at $[u\!:\!v]=[0\!:\!1]$,
that is equivalent to $L_{22}=0$,
which in its turn implies that $Q_\Theta$ is a complete square $-L_{12}^2$, if $L_{12}\ne 0$, or zero otherwise.
Conversely, assume that $Q_\Theta=L_{11}L_{22}-L_{12}^2$ is a complete square, say $-L^2$. If $L=L_{12}$, then either $L_{11}$ or $L_{22}$ vanishes,
which implies that the point of $l$ with $[u\!:\!v]=[1:0]$ or $[0:1]$ respectively is a singular point of $X$.
If $L\ne L_{12}$, then $L_{11}L_{22}=(L_{12}+L)(L_{12}-L)$
from where $L_{11}=\lambda(L_{12}\pm L)$ and $L_{22}=\lambda^{-1}(L_{12}\mp L)$ for some $\lambda\in\C\sm\{0\}$.
Then it is the point $[1:-\lambda]\in l$
that is singular
in $X$.

To check the second statement, let us select projective coordinates in a way
 that it is the point $u=0,v=0,x=0,y=0,z=1$ which is a singular point of $X$ and assume that $\det A_s$ is not
identically zero. Then, $C, Q_1, Q_2$, and $dC$ are zero at this point. This immediately implies that $\det A_s$ is zero at $s=(0,0,1)$, and
after taking the differential of the matrix determinant line by line we observe that this differential  is also zero at this point.
\endproof

\proposition\label{quintic+conic}
For any non-singular cubic threefold $X\subset P^4$ and any line $l\subset X$, the following properties hold:
\roster\item
the spectral quintic $\spec$ cannot have singularities other than nodes;
\item
the associated curve $\Theta$ is a conic which is not a double line;
\item
the conic $\Theta$ passes through every node of $\spec$;
\item
the local intersection index $\ind_s(\spec,\Theta)$ at a
common point $s\in\spec\cap\Theta$ is even for any $s$ not belonging to a common irreducible component of $\spec$ and $\Theta$.
\endroster
\endproposition

\proof
The fact that $\spec$ may have only nodes as singular points follows from Proposition \ref{nodality} applied to the blow-up of $X$ along $l$.
Lemma \ref{l-nonsing} implies (2).

To prove (3)
 note that according to Proposition \ref{nodality} if a point $s\in\spec$ is a node,
then
the residual conic $r_s$ is a double line,
which means that $A_s$ is of rank 1 and thus implies the vanishing of
$
\det\pmatrix
L_{11} & L_{12} \\
L_{12} & L_{22}
\endpmatrix
$ at $s$.

The property (4) is well-known
for generic lines on $X$ (see, for example, \cite{White}; in particular, for a generic line the quintic $\spec$ is non-singular (see \cite{Murre})
and then (4) can be derived from the Laplace formula
$$L_{11}\D=\D_{22}\D_{33}-\D_{23}^2,$$
where $\D$ is the determinant defining
the
quintic $\spec$ and $\D_{ij}$ are the $ij$-cofactors
in the matrix (\ref{spec-matrix}),
among which $\D_{33}$ defines $\Theta$).
After that the property (4) can be extended to any line $l\subset X$ by continuity.
\endproof

\subsection{Multiple lines and a characterization of the nodes of $\spec$ and $\Theta$}
We say that a plane $P\subset P^4$ is {\it tritangent to $X$} if $X$ cuts on $P$
a triple of lines: $P\cap X=l+ l'+ l''$. By definition, $\spec$ parameterizes the tritangent planes of the net $\{P_s\}_{s\in P^2}$.

If for a line $l$ there exists a plane $P$ containing it such that
$P\cap X=2l+l'$, we say that $l$ is a {\it multiple line} and $l'$ the {\it residual line}.
If $l=l'$ (that is $P\cap X=3l$), then $l$ is called a {\it triple line}, and if $l\ne l'$ then a {\it double line}.
Non-multiple lines $l\subset X$ (for which such a plane $P$ does not exist)
will be called {\it simple lines}.

\proposition\label{quintic+conic_via_l}
For any non-singular cubic $X\subset P^4$ and a line $l\subset X$, with the associated spectral curve $\spec$
and the theta-conic $\Theta$, the following holds:
\roster\item
$s\in P^2$ is a node of $\spec$ if and only if $P_s\cap X=l+2l'$ for some line $l'$; hence,
$\spec$ has $\ge1$ nodes if and only if $l$ is a residual line;
\item
$s\in P^2$ is a node of $\Theta$ if and only if $P_s\cap X=2l+l'$ for some line $l'$; hence,
$\Theta$ splits into
two lines if and only if $l$ is a multiple line;
\item
$s\in P^2$ is a common node of $\spec$ and $\Theta$ if and only if $P_s\cap X=3l$;
in this case the two lines forming $\Theta$ are tangent to the two
branches of $\spec$ at the node.
\endroster
\endproposition

\proof
Item (1) follows from Proposition \ref{nodality}.

For proving items (2) and (3) choose the coordinates so that $s=[x\!:\!y\!:\!z]=[1\!:\!0\!:\!0]$, then
in the plane $P_s$ with coordinates $[u\!:\!v\!:\!x]$ the equation of the residual conic is
$$r_s=u^2L_{11}^x+2uvL_{12}^x+v^2L_{22}^x+2uxQ_1^{xx}+2vxQ_2^{xx}+x^2C^{xxx}=0,$$
where $L_{ij}^x=L_{ij}(1,0,0)$, $Q_i^{xx}=Q_i(1,0,0)$, and $C^{xxx}=C(1,0,0)$
are respectively the $x$-coefficients of $L_{ij}$, the $x^2$-coefficient of $Q_i$, and the $x^3$-coefficient of $C$.

The line $l$ in $P_s$ has equation $x=0$, thus, it is double, namely, $P_s\cap X=2l+l'$,
if and only if $x$ divides $r_s$, that is iff $L_{11}^x=L_{12}^x=L_{22}^x=0$.
 This vanishing implies that the equation of $\Theta$ in the affine chart $\{x\ne0\}$,
 $$L_{11}(1,y,z)L_{22}(1,y,z)-L_{12}^2(1,y,z), $$
has no linear terms, which shows that $s=[1:0:0]$ is a node of $\Theta$.
Conversely, if $\Theta $ has a node at $s$, then, as is evident, $L_{11}$, $L_{12}$ and $L_{22}$ are linearly dependent
and each of these linear forms is vanishing at $s$; therefore, the coefficients $L_{11}^x=L_{12}^x=L_{22}^x$
vanish, which shows that $l$ is a double line. This proves item (2).

The first part of item (3) follows from items (1) and (2). For proving the second part of (3),
note that $P_s\cap X=3l$ if and only if $r_s$ is a multiple of $x^2$, that is iff
in addition to $L_{11}^x=L_{12}^x=L_{22}^x=0$ we have $Q_1^{xx}=Q_2^{xx}=0$, whereas $C^{xxx}\ne0$
(the latter is due to non-singularity of $X$). These conditions imply that the degree $\le2$
truncation of equation (\ref{spec-curve}) of $\spec$ in the above affine $(y,z)$-chart reduces to
$$
C^{xxx}(L_{11}L_{22}-L_{12}^2),
$$
where $L_{ij}$ depend only on $y$ and $z$.
Thus, $\spec$ and $\Theta$ both have nodes at $s$ and the branches of $\spec$ are tangent to those
of $\Theta$. \endproof

\subsection{The exceptional divisor $E_l$}\label{exceptional-divisor}
Simple lines $l\subset X$ can be also distinguished from multiple ones by their normal bundles
in $X$, or equivalently, by the exceptional divisor $E_l$ of blowing up $X_l\to X$.

\proposition$($\cite{Clemens-Griffiths}$)$
If a line $l$ on a non-singular cubic threefold $X$
is simple, then its normal bundle in $X$ is trivial (a sum of trivial line bundles  $O_l\oplus O_l$)
and $E_l\cong P^1\times P^1$.
Otherwise, the normal bundle is isomorphic to the sum of line bundles $O_l(1)\oplus O_l(-1)$ and $E_l$ is a geometrically ruled surface
$\S_2$.
\qed\endproposition

\proposition\label{E_L-covering}
In the case of a simple line $l$ the theta-conic $\Theta$ is non-singular and
the restriction
$$\CD E_l\cong P^1\times P^1@>\pi_{E_l}>> P^2 \endCD$$
of the canonical projection $\pi_l\: X_l\to P^2$ is
a double covering branched along $\Theta$.

In the case of a multiple line $l$, conic $\Theta$ splits in two lines and
such a restriction
$$\CD E_l\cong \Sigma_2@>\pi_{E_l}>> P^2\endCD$$
is composed of a contraction of the $(-2)$-curve of $\Sigma_2$
and subsequent double covering
over $P^2$ branched along $\Theta$.
The two lines of $\Theta$ are the images of two fibers of the ruling $\S_2\to P^1$ and the node of $\Theta$ is the image of the $(-2)$-curve.
\endproposition

\proof
The points of $E_l$ can be viewed as pairs $(P_s,t)$, where
$s\in P^2$
and $P_s$ is the plane passing through $s$, containing the line $l$ and tangent to $X$ at
the point $t\in l$.
Projection $\pi_{E_l}$ sending $(P_s,t)$ to $s=P_s\cap P^2$ is generically two-to-one, because $P_s$ is generically tangent at a pair of points
of $l$, namely, at the common points of $l$ with the residual conic $r_s$
(see Section \ref{conic-bundles-for-cubics}). It is one-to one over the non-singular points of $\Theta$, since for them the intersection
$r_s\cap l$ consists of one point. This gives item (1), and for proving item (2) it remains to
notice that for a nodal point $s_0$ of $\Theta$, $P_{s_0}\cap X=2l+l'$, and so, the preimage of $s_0$
consists of all the pairs $(P_{s_0},t)$, $t\in l$.
\endproof

Proposition \ref{E_L-covering}
together with
Lemma \ref{l-nonsing} imply
the following statement.

\corollary\label{double-lines}
For a non-singular cubic threefold $X$ and a line $l\subset X$ the following properties are equivalent:
\roster\item
The line $l$ is multiple.
\item
There is a plane passing through $l$ which is tangent to $X$ at each point of $l$.
\item
The linear forms $L_{11}, L_{12}, L_{22}$ are linearly dependent.
\item
The linear forms $L_{11}, L_{12}, L_{22}$ forms a rank 2 linear system.
\item
The theta-conic $\Theta$ corresponding to $l$ splits into two lines.
\qed\endroster
\endcorollary

\proposition\label{theta-parametrization}
The theta-conic $\Theta$ is formed by the points $s\in P^2$ for which
the binary quadric $l_s=r_s\vert_ l$
is degenerate, that is a square or identically zero.
If $\Theta$ is non-singular, then $l_s$ is non zero for any $s\in \Theta$,
and the map sending $s\in\Theta$
to $p\in l$ such that $r_s(p)=0$
gives an isomorphism
between $\Theta$ and $l$. The inverse map sends $p\in l$ to the point $s\in \Theta$ where the hyperplane $T_p$ tangent to $X$ at $p$ is tangent to
$\Theta$.
\endproposition

\proof
By definition, $l_s=L_{11}u^2+2L_{12}uv+L_{22}v^2$ and
degeneration of $l_s$ means vanishing on $l$ of the determinant that defines $\Theta$.
Non-singularity of $\Theta$ means linear independence of the linear forms $L_{11},L_{12},L_{22}$ (see Corollary \ref{double-lines}), and we can choose our coordinates so that
$x=L_{11}$, $y=L_{12}$ and $z=L_{22}$. Then the inverse map
$l \to \Theta$ to the above correspondence
is defined by $[u:v]\to[x:y:z]$, where $x=v^2$, $y=-uv$, $z=u^2$. It remains to notice that $x=v^2$, $y=-uv$, $z=u^2$ is the only point of intersection of $\Theta$ with the hyperplane $T_p$, $p=[u:v]$, since the latter is given by equation $u^2x+2uvy+v^2z=0$.
\endproof

\subsection{Fano surface and curves $\Gres$, $\GII$}
\label{unbalanced}
The lines on a non-singular cubic threefold, $X$, form a surface, $F(X)$, named after G.~Fano \cite{Fano}, who established its following properties
(see also \cite{Clemens-Griffiths}).

\theorem\label{FanoTh} For any non-singular cubic threefold $X$, the Fano surface $F(X)$ is irreducible, non-singular
and has the following numerical characteristics:
$$
h^{1,0}=h^{0,1}=h^{1,2}=h^{2,1}=5,\quad
h^{2,0}=h^{0,2}=10,\quad h^{1,1}=25.
$$
\endtheorem

As is known (see, for example, \cite{Murre}), the multiple lines $l\subset X$ form a curve on $F(X)$ that we denote $\GII$,
the residual lines form a {\it residual curve}, $\G_{res}\subset F(X)$,
and the triple lines form a finite subset $\GIII\subset \Gres\cap\GII\subset F(X)$.

\proposition\label{gamma_res} The spectral curve
 $S$ is non-singular if and only if $l\notin \G_{res}$.
\endproposition

\proof Straightforward from definition of $\G_{res}$ and Proposition \ref{quintic+conic_via_l}.
\endproof

\subsection{The spectral covering via theta-conics}\label{theta-conics}
Consider a non-singular cubic threefold $X$ and a line $l\subset X$ that is neither multiple nor residual, and so the spectral curve $S$ and the theta-conic $\Theta$
are non-singular. The set of lines on $X$ that intersect $l$ form a curve $\til S\subset F(X)$ and its projection $p_S\:\til S\to S$ induced by the conic bundle $\pi_l$
is clearly a double covering that will be called the {\it spectral covering} ({\it cf.}, Subsection \ref{conic-bundles}).

On the other hand, we consider the pull-back $S_{E_l}=\pi_{E_l}^{-1}(S)$ where $\pi_{E_l}\:E_l\to P^2$ is the double covering
branched along $\Theta$ ({\it cf}.,
Subsection \ref{exceptional-divisor}). This pull-back curve has
$\pi_{E_l}^{-1}(S\cap\Theta)$
as the singular locus. Namely, for a point $x\in S\cap\Theta$,
with the local intersection index $\ind_x(S,\Theta)$ being $2k$, the pull-back $S_{E_l}$ has at the point of $\pi_{E_l}^{-1}(x)$
 two local branches that intersect with the local intersection index $k$.

Note also that for any line $l'\in\til S$, its proper image in $X_l$ intersects the exceptional divisor $E_l\subset X_l$ at a point of $S_{E_l}$,
which yields a map $g_l\:\til S\to S_{E_l}$. Non-singularity of $\til S$ implies the following.

\proposition\label{Prym-cover-identification} The spectral covering $\til S\to S$ is the composition $\pi_{E_l}\circ g_l$ where
$g_l\:\til S\to S_{E_l}$ coincides with
the normalization map and $\pi_{E_l}\:E_l\to P^2$ is the double covering of $P^2$ branched along $\Theta$.
\qed\endproposition

\subsection{The case of real pairs $(X,l)$}
If a cubic threefold $X$ and line $l\subset X$ are real, then the associated spectral quintic $S$ and theta-conic $\Theta$ are also real.
In the most of constructions we assume that  $X$ and $S$ are non-singular, then
$\Theta_\R\ne\oo$ may be either (1) non-singular or (2) split into a pair of distinct lines,
which can be either (2a) both real, or (2b) form an imaginary complex-conjugate pair.
By a {\it region bounded by $\Theta_\R$} we mean the closure of a connected component of $\Rp2\sm\Theta_\R$; there are two such regions in the cases (1) and (2a)
and one region in the case (2b).
A non-singular quintic $S_\R$ contains always a one-sided connected component, $J\subset\Rp2$, and may have in addition, a few two-sided components, called {\it ovals}.
Since conic $\Theta$ has with $S$ even local intersection indices at the common points,
each connected component of $S_\R$ is contained precisely in one region bounded by $\Theta_\R$,
and the region containing $J$-component is denoted by $F_\Theta^+$ and
called the {\it exterior of $\Theta_\R$}, while the other region (non-empty in the cases (1) and (2a)) is denoted by  $F_\Theta^-$
and called the {\it interior of $\Theta_\R$}.

We say that an oval of $S$ is {\it visible} (respectively,  {\it invisible}) if it lies in the exterior (respectively, interior) of $\Theta_\R$.
Let $W_S^+=\pi_l(X_{l\R})\subset\Rp2$ and $W_S^-=\Cl(\Rp2\sm W_S^+)$.

\proposition\label{visible} If $X$ is a real non-singular cubic threefold and $l\in F_\R(X)\sm \G_{res}$, then $W_S^+$, $W_S^-$ are compact surfaces with common boundary
$W_S^+\cap W_S^-$  which is the union of invisible ovals.

\proof By Proposition \ref{gamma_res}, the spectral curve $S$ is non-singular.
According to Proposition \ref{nodality}, a point $s\in \Rp2\sm S_\R$ is an interior point of $W_S^+$ (respectively, $W_S^-$) if the real part of the residual conic $r_s$ is a topological circle (respectively, empty), while
for $s$ passing from one side of $S_\R$ to another the conic $r_s$ experiences a Morse transformation. The latter implies that
 a point $s\in S_\R$ is an interior point of $W_S^+$ (respectively, belongs $W_S^+\cap W_S^-$) if $r_s$ is a pair of real lines (respectively, a pair of imaginary lines). Thus, it remains to notice
that, for $s$ belonging to the $J$-component, the conic $r_s$ is a pair of real lines, since the $J$-component being one-sided, the topology of the real locus of the residual conics is the same on the both sides of
this component. \endproof


\section{Theta-charactersitics and the spectral correspondence}\label{Theta}

\subsection{Theta-conics and theta-charactersitics}\label{theta-conics-and-characteristics}
Recall that a theta-characteristic on a smooth projective curve $S$  over $\C$ is defined as the isomorphism class of
pairs $\theta=(L,\phi)$ formed by a line bundle $L$ and a quadratic map $\phi \: L\to K_S$ to the  canonical bundle $K_S$ that provides an isomorphism $L^{\otimes2}\cong K_S$
(in other words, that are solutions of equation $2D=K_S$ in the divisor class group).

Since the quadratic map $L \to K_S$ can be viewed as a double covering between the associated $U(1)$-bundles, theta-characteristics can be seen as Spin-structures in the general context of Spin-manifolds
(see \cite{Atiyah}, \cite{Mumford}).
There is also a well-known bijection
(see \cite{Johnson}) between the set of
Spin-structures on $S$ and the set of quadratic extensions of the intersection
 form on $H_1(S;\Z/2)$,
that is the set $\Quadr(S)$ of functions
$q\:H_1(S;\Z/2)\to\Z/2$
satisfying the relation $q(x+y)=q(x)+q(y)+xy$
for all $x,y\in H_1(S;\Z/2)$, where $xy\in\Z/2$ is the intersection index of $x$ and $y$.
 Possibility to identify these three objects, theta-characteristics, Spin-structures, and quadratic functions, will let us
easily switch from one language to another.

Curves endowed with a theta-characteristic are called {\it Spin-curves}, and the set of
all the theta-characteristics for a given $S$ is denoted $\Theta(S)$.
A theta-characteristic
$\theta\in\Theta(S)$
can be even or odd, depending on the parity of
$h^0( \theta)=\dim H^0(L)$, and
we denote by $\Theta_0(S)$ and $\Theta_1(S)$ the sets of even and odd theta-characteristics.
In terms of quadratic functions the parity coincides with the {\it Arf-invariant}.

The set  $\Theta(S)$ has a natural free transitive action of $H^1(S;\Z/2)$: in terms of
quadratic functions an element $h^*\in H^1(S;\Z/2)$ acts
on $q\in\Quadr(S)$ by addition $q\mapsto q+h^*$, if $h^*$ is interpreted as a linear function $h^*\:H_1(S;\Z/2)\to\Z/2$.
By the difference $\theta_1-\theta_2$ of $\theta_1,\theta_2\in\Theta(S)$ (or equivalently, the difference
$q_1-q_2$ of $q_1,q_2\in\Quadr(S)$) we mean the class $h^*\in H^1(S;\Z/2)$ such that $\theta_1=\theta_2+h^*$
(respectively, $q_1=q_2+h^*$). Such a {\it difference class} defines the {\it difference
double covering} $\til S\to S$. When $\theta_1,\theta_2$ are represented by divisor classes $D_1, D_2$ (with $2D_1=2D_2=K_S$)
this difference class $h^*$ is Poincar\'e dual to $h\in H_1(S,\Z/2)$ represented
(via the natural identification of the 2-torsion part of the Jacobian with $H_1(S;\Z/2)$)
by the image $[D_1-D_2]$ of the order two divisor class $D_1-D_2$
under the Abel-Jacobi map (extended to arbitrary divisors by linearity).

A real structure $\conj\: S\to S$ on a curve $S$
induces an involution $c_\Theta\:\Theta(S)\to\Theta(S)$ on the set of its theta-characteristics, and we call $\theta\in\Theta(S)$ {\it real}
if it is a fixed point of $c_\Theta$, and use notation $\Theta_\R(S)=\Fix(c_\Theta)$ for the set of all real $\theta$.
In terms of the quadratic function $q_\theta$, this is equivalent to $q_\theta(x)=q_\theta(c(x))$, for all $x\in H_1(S)$,
where $c$ denotes $\conj_* : H_1(S)\to H_1(S)$.

As is known (see, for example, \cite{Atiyah}),
for real curves $S$ with non-empty $S_\R$ a theta-characterisitic $\theta=(L,\phi)$ is real if and only if
there exists a fiberwise antilinear
involution of the line bundle $L \to L$ that covers $\conj : S\to S$.

Plane curves $S$ of an odd degree $d$ carry a special theta-characteristic, $\theta_0=\theta_0(S)$,
traced on $S$ by curves of degree $\frac{d-3}2$.
If $S$ is real, then $\theta_0$ is real too.
In the language of quadratic functions
this theta-characteristic was introduced
by Rokhlin \cite{Ro}, who also observed that $q_0$ takes value $1$ on every vanishing cycle of $S$
and its parity
is $\frac18(d^2-1)\mod2$. We name $\theta_0$ {\it Rokhlin's theta-characteristic}
and denote $\Theta_1(S)\setminus \{\theta_0\}$
by $\Theta_1^*(S)$.
If $S$ is real, then $\Theta_{\R 1}^*(S)$ states for $\Theta_\R(S) \cap \Theta_1^*(S) $.

For $h\in H_1(S;\Z/2)$, we let $q_h=q_0+h^*$, $\theta_h=\theta_0+h^*$ where $h^*$ is the Poincar\'e dual to $h$, and conclude that $h\mapsto q_h\mapsto \theta_h$ identifies
$H_1(S;\Z/2)$
with $\Quadr(S)$ and $\Theta(S)$. We say that $h^*$ is the {\it difference class associated with the theta-characteristic $\theta=\theta_h$.}

\subsection{Spectral theta-characteristic}\label{def-spectral-theta}
Given a non-singular plane quintic $S$,
a conic $\Theta$ is said to be {\it contact to}
$S$ if $\Theta $ is reduced and at every common points $s\in S\cap\Theta$ the
local intersection index $\ind_s(S,\Theta)$ is even.
To each conic $\Theta$ contact to $S$ we associate its
{\it contact divisor} on $S$ defined as
the half of the intersection divisor $S\cdot\Theta$.
Each contact divisor is a theta-characteristic, since $S\cdot\Theta$ represents the canonical divisor class of $S$.

 \proposition\label{contact-divisors}
 For any non-singular quintic $S$, the correspondence defined by taking the contact divisor of a contact conic  gives a
 bijection between  $\Theta^*_1(S)$
 and the set of contact conics to $S$. If $S$ is real, this correspondence gives also a bijection between
 $\Theta_{\R 1}^*(S)$ and the set of real contact conics.
 \endproposition

 \proof
Each odd theta-characteristic of $S$ is represented by a divisor $D$ such that $2D$ is cut on $S$ by a conic. If the latter conic is not reduced, the corresponding theta-characteristic is Rokhlin's one $q_0$, and then 
 $h^0(\theta_0)=3>1$. On the other hand, according to \cite{Beau-Det} Proposition 4.2,
a theta-characteristic $\theta\in  \Theta^*_1(S)$ does correspond to a contact conic, if $h^0(\theta)=1$.
Hence, to prove the first statement it is sufficient to show
that two distinct contact divisors can not be linear equivalent on $S$. Let us assume on the contrary
 that two distinct contact divisors belong to the same divisor class, that is provide two distinct sections of the same square root $\theta$
 of the canonical bundle of $S$. By definition of contact divisors, the contact conics are reduced. Hence, the linear combinations of the given
 sections provide, by squaring,  a local one-parameter family of equisingular reduced conics $\Theta_t$ having the same number of tangencies with $S$
 such that the order of tangencies is preserved in the family. Then,
we can apply the Gudkov-Shustin inequality to this family (its idea goes back to \cite{GS}, Theorem 2;
a bit more elaborated version that we make use of appeared in \cite{GLS}, Lemma II.2.18)
that gives a bound
$
 \Theta_s\circ \Theta_t\ge \sum ((\Theta_t\circ S)_{z}-1),
$
where the sum is taken over all the points $z\in \Theta_t\cap S$. This bound
leads to a contradiction, since $ \Theta_s\circ \Theta_t =4$ but  $\sum ((\Theta_t\circ S)_z-1)\ge 2\times 5 - 5=5$.

Since the real part of a real quintic is always non empty, for each
$\theta=(L,\phi)\in \Theta_{\R 1}^*(S)$
there exists an
antilinear involution
$L\to L$ which covers
$\conj : S\to S$, and therefore
the second statement follows from the first one.
\endproof

\proposition\label{difference-class-covering}
Let $S$ and $\Theta$ be the spectral quintic and the theta-conic associated with a line $l$ on a non-singular cubic threefold $X$.
Assume that $l\notin\Gres$.
Then the difference class $h^*\in H^1(\spec;\Z/2)$ associated with the theta-characterstic $\theta$ defined by $\Theta$ is non-zero.
The difference double covering $\til{\spec}\to\spec$ defined by the difference class $h^*$
is non-trivial and coincides with the spectral
covering associated with $(X,l)$.
\endproposition

\proof The first statement follows from $\theta\in\Theta^*_1(\spec)$ (see Proposition \ref{contact-divisors}).
Since, by definition, the difference class $h^*$ associated with the
theta-characterstic $\theta$ defined by $\Theta$ (that is $h^*=\theta -\theta_0$) is Poincar\'e dual to
$h=[D - H]$, where $D$ is the contact divisor
and $H$ stands
for a hyperplane ({\it i.e}., line) section,
the second statement is a straightforward consequence of
the first one and Proposition \ref {Prym-cover-identification}.
\endproof

The odd theta-characteristic $\theta\in\Theta^*_1(\spec)$ corresponding to $\Theta$
will be called the {\it spectral theta-characteristic}.

\subsection{Recovery of $(X,l)$ from $\spec$ and $\Theta$}
The following result belongs essentially to F.~P.~White \cite{White}.

\theorem\label{tang-conic}
Given a non-singular quintic $S\subset P^2$
and a conic $\Theta\subset P^2$ contact to $S$,
there exists a non-singular cubic threefold $X$ and a line $l\subset X$ whose associated spectral curve is $S$
and the theta-conic
is $\Theta$. Such a pair $(X,l)$ is unique up to a projective equivalence.

If $S$ and $\Theta$ are real, then
$X$ and $l$ can be chosen real, and they are defined uniquely up to a real projective equivalence.
\endtheorem

\proof Over $\C$, the existence statement is proved in \cite{White} by an
explicit construction of the
polynomials $L_{ij}$, $Q_{i}$, $i,j\in\{1,2\}$, and $C$ of degrees $1, 2, 3$ respectively,
such that
$\Theta$ and $S$ are  given by equations (\ref{theta-conic}) and (\ref{spec-curve})
(a modern proof is given in \cite{Beau-Det}). If
$S$ and $\Theta$ are real, then by choosing the auxiliary lines and auxiliary collections of points involved in this construction to be invariant under the complex connugation
one makes the above polynomials, and hence the above equations, real.
The absence of singular points on the cubics $X$ constructed in this way
follows from Lemma \ref{l-nonsing}.

The uniqueness  statement is
also well known, see, for example, Appendix C in \cite{Clemens-Griffiths}. Unfortunately, the proof is not explicit there, and we
could not find
any complete one elsewhere.
The proof that we give below
is different (although should be known to experts) and can
be easier adapted to the real setting.

Let $S$ and $\Theta$ be the spectral quintic and the spectral conic associated with a real line $l$ on a real non-singular cubic threefold $X\subset P^4$.
We assume, like in the statement, that $S$ is non-singular and show how
to reconstruct the pair $(X,l)$ from $(S,\Theta)$ canonically up to real projective equivalence.

Recall, first,
that for each $s\in S$ the projective plane $P_s$ generated by $s$ and $l$ intersects $X$ along $l+l_s'+l_s''$ where $l$, $l'_s$, and $l''_s$ are three distinct lines
unless  $\Theta$ is singular and $s$ is its node: in this exceptional case $l$ is equal to one of $l_s', l_s''$.
Consider the mapping $\phi : S\to  P^4$ that sends $s\in S$ to the intersection point $l'_s\cap l''_s\in P^4$.
It is well known to be
an embedding associated with the linear system $\vert L+D \vert$,
where $L$ stands for sections of $S$ by the lines in $P^2\supset S$ and $D$ is the contact divisor,
$D= \frac12 (S\cdot \Theta)$
{\it cf.,} \cite{Beau}, Remarque 6.27.
So, the pair $(S,\Theta)$ determines the curve $\phi(S)$ up to real projective equivalence.
Furthermore, since $l'_s\cap l''_s\in l$ if and only if $s\in S\cap \Theta$, we can also reconstruct the line $l$ from $\phi(S)$
as the unique line passing through $\phi(D)$.

Next, we consider the projection $\pi_{E_l}\:E_l\to P^2$ given in Proposition \ref{E_L-covering}.
In the case of non-singular $\Theta$, it
identifies the line $l$ with $\Theta$ so that the two points, $t'=l\cap l'_s$ and
$t''=l\cap l''_s$, coincide with the tangency points to $\Theta$ of the two tangent lines through the point $s$.
If the conic $\Theta$ has a node at $s_0$, then $E_l$ can be obtained by blowing up $\hat P^2\to P^2$ at $s_0$ and then taking
the double covering $E_l\to \hat P^2$ branched
along the proper image of $\Theta$.
The line $l$ is then identified with the $(-2)$-curve $C_2\subset E_l$ that covers the exceptional $(-1)$-curve $C_1\subset\hat P^2$,
and  the points $t',t''\in l$ become represented by the inverse image under the covering $C_2\to C_1$ of the point corresponding to
the line $s_0s$.

Thus, we have reconstructed
the lines $l'_s, l''_s\subset X$ for each generic $s\in S$.
Finally, we can conclude that this collection of lines
determines the whole $X$, because it determines all its generic
hyperplane sections passing through $l$ as it follows from
Schl\"afli's theorem \cite{Sch}:  given five skew lines in $P^3$ and a line intersecting them all, there exists a unique cubic surface containing these six lines.
\endproof

\subsection{Deformation class $[S]$ as a real deformation invariant of $(X,l)$}\label{spectral-cor}
Consider the variety $\Cal C$ formed by pairs $(X,l)$ where $X\subset P^4$ is a non-singular cubic 3-fold and $l$ is a line on it, as well as its two
hypersurfaces
$$\D_{\Cal C}^\res=\{(X,l)\in\Cal C\,|\,l\in\Gres(X)\},\quad \D_{\Cal C}^{I\!I}=\{(X,l)\in\Cal C\,|\,l\in\GII(X)\}$$
which are responsible for the associated spectral curve and theta-conic being singular.
Put $\Cal C^{*}=\Cal C\sm\D_{\Cal C}^\res$.

In addition, consider  the variety $\Cal S$ formed by pairs $(S,\Theta)$
where $S\subset P^2$ is a non-singular quintic and $\Theta$ is a conic contact to $S$.

The {\it spectral correspondence} associates to $(X,l)\in\Cal C^*$ its spectral curve and theta-conic,
that is a pair $(S,\Theta)\in\Cal S$.
The latter pair depends of course on the choice of a plane $P^2$
chosen as the target of the central projection from $l$.
So, it is natural to pass to the varieties of projective classes (moduli spaces), and
to
speak of the spectral correspondence as
a well-defined regular morphism between $\Cal C^*/PGL(5; \C)$ and
${\Cal S}/PGL(3; \C)$. This morphism is defined over the reals.
According to
Theorem \ref{tang-conic}, both morphisms, $\Phi \: \Cal C^*/PGL(5; \C)\to {\Cal S}/PGL(3; \C)$
and $\Phi_\R \:\Cal C_\R^*/PGL(5; \R)\to {\Cal S_\R}/PGL(3; \R)$, are bijective (in fact, the both are isomorphisms;
this stronger statement will be commented a bit more in
Subsection \ref{full_correspondence}, but not actually used in our paper).

The rest of this section is devoted to proving
the following result.

\proposition\label{real-spectral-correspondence}
If $(X_i,l_i)\in\Cal C^*_\R, i=0,1,$ are real deformation equivalent (can be connected by a path in $\Cal C_\R$)
and $(S_i,\Theta_i) \in \Phi_\R [(X_i,l_i)], i=0,1$, then
$S_i$ are also real deformation equivalent (can be connected by a path in $\Cal C ^{5,1}_\R\sm \Delta^{5,1}_\R$).
\endproposition

First, we introduce and analyze codimension in $\Cal C$ of the following subvarieties:
$$\D_{\Cal C}^{I\!I\!I}=\{(X,l)\in\Cal C\,|\,l\in\GIII(X)\} \quad \text{(where  $\GIII(X)$ is defined in Subsection \ref{unbalanced})}
$$
that is responsible for $S$ and $\Theta$ sharing a common nodal point,
and
$$\D_{\Cal C}^{\roman{split}}=\{(X,l)\in\Cal C\,|\, S \text{ and } \Theta \text{ have a common irreducible component}
\}.
$$

\lemma\label{deep-strata}
The codimension of the varieties $\D_{\Cal C}^{I\!I\!I}$ and $\D_{\Cal C}^{\roman{split}}$ in $\Cal C$ is $\ge 3$.
\endlemma

\proof
If $(X,l)$ belongs to
$\D_{\Cal C}^{I\!I\!I}$ or $\D_{\Cal C}^{\roman{split}}$ then either $\Theta$ splits into two lines or $\Theta$ is an irreducible component of $S$.
Indeed, if $\Theta$ belongs to $\D_{\Cal C}^{I\!I\!I}$, then $\Theta$
and $S$ have a common node (see Proposition \ref{quintic+conic_via_l}(3)), and in particular,
$\Theta$ is reducible.

Consider first the case of $\Theta$ splitting into a pair of lines.
Let us choose a coordinate system $[x\!:\!y\!:\!z\!:\!u\!:\!v]$ in $P^4$ so that the
line $l$ is given by equations $x=y=z=0$ and the conic $\Theta$ in $P^2$ with coordinates $[x\!:\!y\!:\!z]$ by equation $xy=0$.
Then, after a suitable linear coordinate change in $u$ and $v$,  the defining
minor (\ref{theta-conic}) of the fundamental matrix (\ref{spec-matrix}) of $X$
can be transformed into
$$
\pmatrix
L_{11} & L_{12} \\
L_{12} & L_{22}
\endpmatrix
=\pmatrix
x&0\\
0&y
\endpmatrix.
$$
Thereby, in the affine chart $(x,y)$ in $P^2$ the linear part of the polynomial (\ref{spec-curve}) defining $S$ is equal to
$-(xQ_2^{zz}+yQ_1^{zz})$, where $Q_i^{zz}$ is the coefficient of $Q_i$, $i=1,2$, at $z^2$.
So, $S$ is singular at $(x,y)=(0,0)$ if and only if $Q_1^{zz}=Q_2^{zz}=0$.
These two conditions along with the third condition of nodality of $\Theta$
(that is a linear dependence of the three linear forms $L_{ij}$)
are independent, and thus we conclude that $\D_{\Cal C}^{I\!I\!I}$ has codimension 3.

If $S$ and $\Theta$ have a common line-component, for instance in the above coordinates, $S$ contains
line $x=0$, then the summand $yQ_1^2$ in the determinant equation of $S$ must be divisible by $x$, which
means vanishing of the three coefficients $Q_1^{yy}$, $Q_1^{yz}$, and $Q_1^{zz}$ in the monomials that do not contain $x$.
Independence of these three coefficients implies that $\D_{\Cal C}^{\roman{split}}$ has codimension $\ge3$ (in fact, $4$ because of an extra condition of reducibility of $\Theta$)
in the case of a common line component.

At last, we consider the case of $\Theta$
being an irreducible component of $S$. Then,  by a coordinate change the theta-conic minor (\ref{theta-conic})
can be transformed into
$$
\pmatrix
L_{11} & L_{12} \\
L_{12} & L_{22}
\endpmatrix
=\pmatrix
x&z\\
z&y
\endpmatrix.
$$
Whereafter vanishing of $\det A_s$ on $\Theta=\{xy-z^2\}$ is equivalent to vanishing
of $2zQ_{13}Q_{23}-xQ_{23}^2-yQ_{13}^2$. To conclude the proof it is sufficient to notice that
already vanishing of it at the points $(1,0,0)$, $(0, 1,0)$, and $(1,1,1)$
imposes 3 independent conditions.
\endproof

Lemma \ref{deep-strata} implies in particular that a generic path in the real locus $\Cal C_\R$ does not intersect
$\D_{\Cal C\R}^{I\!I\!I}\cup \D_{\Cal C\R}^{\roman{split}}$. We let
$$
\Cal C^{**}_\R=\Cal C_\R\sm(\D_{\Cal C\R}^{I\!I\!I}\cup \D_{\Cal C\R}^{\roman{split}}).
$$

\lemma\label{regularity-S-corresppjdence}
Let $\{(X_t,l_t)\}_{t\in[0,1]}$ be a path in $\CC_\R^{**}$ that intersects
 $\D_{\CC\R}^\res$
 just once, at $0<t=t_0<1$, and not at a point of $\D_\CC^{II}$, and let
$(S'_i,\Theta'_i) \in \Phi_\R [(X_i,l_i)], i=0,1$. Then $S'_0$ and $S'_1$ are non-singular and real deformation equivalent.
\endlemma

\proof Let us choose continuously varying with $t$ projective coordinates $(x_t,y_t,z_t,$ $u_t,v_t)$ in $P^4$ so that $l_t$ is defined by equations $x_t=y_t=z_t=0$.
It associates with $\{(X_t,l_t)\}_{t\in[0,1]}$  a well defined path $(S_t,\Theta_t)$ in $\Cal S_\R$. Connectedness of $PGL(3;\R)$ allows us to adjust the path so that
$S_0=S'_0$ and $S_1=S'_1$.

By assumption, and according to Propositions \ref{quintic+conic} and  \ref{quintic+conic_via_l}, the quintic $S_{t_0}$ is nodal, the other quintics $S_t$ in the family are non-singular,
and the conic $\Theta_{t_0}$ is smooth at the nodes of $S_{t_0}$.
We need to justify that for each node of $S_{t_0}$ the topological type of smoothing for $t<t_0$ and for $t>t_0$
is the same.

In fact,
if a node of $S_{t_0}$ is cross-like, then smoothness of
$\Theta_{t_0}$ and
the properties (3)--(4) from Proposition
\ref{quintic+conic} imply that $\Theta_{t_0}$ has an odd local intersection index with
each of the real branches of $S_{t_0}$ at the node. This forbids the topological type of
smoothing shown on Figure \ref{nodal-smoothing}, since it generates real points of odd intersection index between $S_t$ and $\Theta_t$
near the node, for
$t\ne t_0$ close to $t_0$.

\midinsert
\hskip30mm
\epsfbox{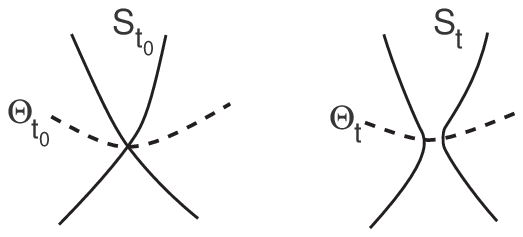}
\vskip-3mm
\figure{Forbiden smoothing}
\label{nodal-smoothing}\endfigure
\endinsert

If a node is solitary, then at such a node the local intersection index $\ind_x(\spec_{t_0},\Theta_{t_0})$ is equal to $2$. Thus, in a small neighborhood of the node, for
$t\ne t_0$ sufficiently close to $t_0$, the curves $\spec_t, \Theta_t$ have one and only one tangency point. Hence, this point is real, which
allows only the birth-of-an-oval type of smoothing both for $t<t_0$ and for $t>t_0$.
\endproof
One can say that
the path $\{S_t\}$ appearing in the proof of this Lemma is just a return
back after touching the discriminant
(see Figure 2).

\hbox{\hskip20mm\vtop{\hbox{\hskip0mm a path in $\CC_\R$ }
\hbox{$\xymatrix@1@=20pt{
{}&{}\ar@{-}[d]&\\
&{}\ar@{-}[d]\ar@{~>}[ul]\ar@{<~}[dr]^{(X_t,l_t)}&\\
&\ \ \D^\res_{\CC\R}&{\phantom{AA}}\\
}$}}
\hskip20mm
\vtop{\hbox{\hskip5mm a path in $\Cal C^{5,1}_\R$ }
\hbox{$\xymatrix@1@=20pt{
{}&{}\ar@{-}[d]&{}\ar@{<~}[dl]\\
&{}\ar@{-}[d]\ar@{<~}[dr]^{S_t}&\\
&\ \ \D^{5,1}_\R&{}\\
}$}}}
\figure{
Crossing of the walls in $\CC_\R$
and in $\Cal C^{5,1}_\R$.}\label{crossing-walls}\endfigure
\vskip3mm

\demo{Proof of Proposition \ref{real-spectral-correspondence}}
 Clearly, $\Cal C$ is a non-singular (quasi-projective) variety. Hence,
any path $(X_t,\ell_t)$, $t\in[0,1]$, in $\Cal C_\R$ with endpoints in $\CC^*_\R$ can be approximated by a smooth path in $\Cal C^{**}_\R$
that intersects
 $\D_{\Cal C\R}^\res$
 at a finite number of points, none of which belongs to $\D_\CC^{II}\cap \D_{\Cal C\R}^\res$. By Lemma \ref{regularity-S-corresppjdence}
 the deformation class of the associated spectral quintics remains the same after each crossing.
 \qed\enddemo


\section{Spectral matching correspondence}\label{SMC}

Throughout this section, for any given real algebraic curve $X$,
we denote by $c$ the involution induced in $H_1(X;\Z/2)$
by the complex conjugation, $\conj\:X\to X$.

\subsection{Smith discrepancy \cite{DK}}\label{Smith}
The Smith theory yields the {\it Smith inequality} $b_*(X_\R;\Z/2)\le b_*(X;\Z/2)$ for any real algebraic variety $X$ and its real locus $X_\R$,
where $b_*(\boxdot,\Z/2)=\dim_{\Z/2}H_*(\boxdot;\Z/2)$ is the {\it total $\Z/2$-Betti number}.

 The gap $b_*(X;\Z/2)-b_*(X_\R;\Z/2)$ is always even and we call its half $d_X=\frac12(b_*(X;\Z/2)-b_*(X_\R;\Z/2))$ the {\it Smith discrepancy of $X$}.
We say that $X$ is an {\it M-variety} (an M-curve, an M-surface, etc.) if $d_X=0$, and in the case $d_X>0$, an {\it $(M-d)$-variety}.

\proposition\label{quintic_discrepancy}
For any real non-singular algebraic curve $X$ with $X_\R\ne \emptyset$  the Smith discrepancy $d_X$ is equal to
the rank of the homomorphism
$\id+c$ in $H_1(X;\Z/2)$. For any real non-singular projective $n$-dimensional hypersurface $X$ of odd degree
the Smith discrepancy $d_X$ is equal to the rank of the homomorphism $\id+c$ in $H_n(X;\Z/2)$.
\qed\endproposition

\subsection{Klein type \cite{DK}}\label{Klein-type}
A real non-singular algebraic variety $X$ is said to be of
{\it Klein type I} if $X_\R$ is null-homologous in $H_n(X;\Z/2)$, $n=\dim_\C(X)$,
and of {\it Klein type II} otherwise. Note that a non-singular irreducible real algebraic curve $X$ with $X_\R\ne\oo$
is of type I if and only if $X\sm X_\R$ splits in two components.

One can reformulate this definition in terms of characteristic elements of the linear map
$H_n(X;\Z/2)\to\Z/2$,
$x\mapsto x\circ c(x)$,
that (due to Poincar\'e duality)
can be represented by a unique homology class $w_c=w_c(X)\in H^n(X;\Z/2)$ such that $x\circ c(x)=x\circ w_c$ for any $x\in H_n(X;\Z/2)$.
By famous Arnold's lemma, $w_c$ is given by the
$\Z/2$-fundamental class $[X_\R]\in H_n(X;\Z/2)$, which implies
the following result.

\proposition\label{arnold_typeI}
A real non-singular algebraic variety $X$
is of Klein type I if and only if $w_c(X)=0$.
\qed
\endproposition

\subsection{Topological and deformation classification of non-singular real plane quintic curves \cite{Kh-quintics},\cite{DK}}\label{quadrocubics-on-quadrics}
The real locus $A$ of a non-singular real plane curve $A$ of odd degree has always one, and only one, one-sided component, which we denote $J$.
The other, two-sided, components (if exist) are called {\it ovals}. Each oval bounds a disc in $\Rp2$ called its {\it interior}
and the complementary M\"obious band
called its {\it exterior}. A pair of ovals is said to be
{\it nested} if one of them contains the other in its interior, otherwise they are called {\it disjoint}.

The Smith inequality applied to a non-singular real plane quintic says that the number of its ovals is  $\le 6$.
We use the code $J\+k$ to refer to the case of $k$ disjoint ovals, and the code  $J\+1\la1\ra$ for the case of two nested ovals.
No other arrangement is possible.

Non-singular real plane quintics form 9 deformation classes.
One of the classes is formed by quintics of type  $J\+1\la1\ra$, quintics in the other 8 classes have no nested ovals.
Among them two ``twin-classes'' refer to quintics with 4 ovals, but in one class the quintics are of  Klein type I,
and in the other of type II. We use the code $J\+4_I$ to refer to the first of these classes and code $J\+4_{II}$
 for the second class of quintics.
The remaining 6 deformation classes of quintics are distinguished by the number $k$ of ovals and have the codes $J\+k$, $k=0,1,2,3,5,6$.
We frequently use $J$ as abbreviation for $J\+ 0$, and $J\+ 4$ for $J\+4_{II}$.

There is a simple criterion to distinguish quintics of types $J\+4_I$ and $J\+4_{II}$
in terms of convexity of the oval arrangement.
Namely, we say that four ovals of a real plane quintic are {\it not in convex position},
if, for three points $A,B,C\in P^2_\R$ chosen in the interior of three distinct ovals between the four, the triangle ABC whose edges do not intersect the J-component,
contains the fourth oval in its interior.
Otherwise, we say that the four ovals are in {\it convex position}.

\lemma\label{convexity-criterion}\cite{Fiedler}
A real plane nonsingular quintic having exactly four ovals is of type I, if and only if its four ovals are not in a convex position.
\qed
\endlemma

\subsection{Topological and deformation classification of non-singular real cubic
threefolds \cite{Kr-rigid}, \cite{Kr-top}, \cite{FK} }\label{real_cubics}
Non-singular real cubic threefolds form 9 deformation classes.
Two of these classes are ``twins'' in the sense that
they are formed by cubics $X$ with the same topological type of the real locus,
$X_\R = \Rp3\#_3(S^1\times S^2)$. These twins are distinguished by their Klein type, and we denote by $C^3_I$
that twin-class
which has cubics of Klein type I, and reserve notation $C^3_{II}$ for the other twin-class, with cubics of type II. We use also $C^3$ as abbreviation for $C^3_{II}$.

For the other seven deformation classes of cubics,
the topology of $X_\R$ determines $X$ up to a real deformation.
For five of these deformation classes, $X = \Rp3\#_k(S^1\times S^2)$, $k=0,1,2,4,5$ and we
notate them $C^k$. There is
one deformation class, where the real locus of cubics is disconnected,
$X_\R=\Rp3\+S^3$, and we denote it by $C^1_{I(2)}$.
The remaining deformation class is called {\it exotic} and is denoted by $C^1_I$,
in this case $X_\R$  is a Seifert manifold
$\Sigma(2,4,6)$.

The above classifications are summarized in the first two columns of Table 1 (see also Figure \ref{adjacency}).
\midinsert\topcaption{Table 1. Real cubic threefolds and their Fano surfaces}\endcaption
$\boxed{\matrix
 \text{deformation type}&\quad \text{real cubic locus } X_\R&\quad \text{real Fano locus } F_\R(X)\\
 \text{-----------------------}&\quad \text{-----------------------}&\quad \text{-----------------------------}\\
C^{k},0\le k\le5& \quad \Rp3\#k(S^1\times S^2)&N_5\+\binom{k+1}2T^2\\
 C^{3}_I& \quad \Rp3\#3(S^1\times S^2)&N_5\+6T^2\\
 C^{1}_I& \quad\Sigma(2,4,6)&\Rp2\+N_6\\
 C^1_{I(2)}& \quad\Rp3\+S^3&N_5\\
\endmatrix}$
\endinsert

\subsection{Topology of real Fano surfaces \cite{Kr-Fano}}\label{real_Fano}
For any non-singular real threefold $X$ of type $C^k$,
$0\le k\le 5$, the real locus of the Fano surface, $F_\R(X)=N_5\+\binom{k+1}2 T^2$, is
formed by disjoint union of $\binom{k+1}2$ copies of a torus (in particular, no tori for $k=0$) and
a non-orientable surface $N_5$ that is a connected sum of $5$ copies of $\Rp2$.
For $X$ of type $C^3_I$ the topological type $F_\R(X)=N_5\+6T^2$ is the same as for
for cubics of type $C^3$.
The case of type $C^1_I$ is the only one where the Fano surface has two non-orientable components, $F_\R(X)=N_6\+\Rp2$.
Finally, $F_\R(X)=N_5$ for cubics of type $C^1_{I(2)}$.
This classification is summarized in the third column of Table 1.

\subsection{Adjacency graphs \cite{Kr-rigid}, \cite{Kh-quintics} }
Our aim here is to relate two adjacency graphs (see their definition in ``Conventions'' in Section 1):
graphs $\G_{5,1}$ and $\G_{3,3}$ related to the plane quintics and cubic threefolds respectively.
These two graphs are shown on Figure \ref{adjacency}. They are combinatorially isomorphic
although their additional decoration with the Smith discrepancies and Klein type, as indicated on this Figure, makes a difference at one vertex.
 To explain this similarity was actually one of our principal motivations for this research
 and the matching correspondence theorem stated at the end of this section achieves this goal.
\midinsert
\hbox{\hskip7mm\hsize60mm\vtop{\hbox{\hskip8mm$\G_{5,1}$ \eightrm(plane quintics)}
\hbox{$\xymatrix@1@=12pt{
&\boxed{J\+6}\ar@{-}[d]&\\
&{J\+5}\ar@{-}[d]&\\
\boxed{J\+4_I}\ar@{-}[dr]&{J\+4}\ar@{-}[d]&\\
&{J\+3}\ar@{-}[d]&\\
\boxed{J\+1\la1\ra}\ar@{-}[dr]&{J\+2}\ar@{-}[d]&\\
&{J\+1}\ar@{-}[d]&\\
&{J}&\\
}$}}
\hskip10mm
\vtop{\hbox{$\G_{3,3}$ \eightrm(cubic threefolds)}
\hbox{
$\xymatrix@1@=9pt{
&\boxed{\ C^5}\ar@{-}[d]&\\
&{\ C^4}\ar@{-}[d]&\\
\boxed{C^3_I}\ar@{-}[dr]&{\ C^3_{}}\ar@{-}[d]&\\
&{\ C^2}\ar@{-}[d]&\\
\boxed{C^1_{I}}\ar@{-}[dr]&{\ C^1}\ar@{-}[d]&\boxed{C^1_{I(2)}}\ar@{-}[dl]\\
&{\ C^0}&\\
}
$}}\hskip10mm
\vtop{\hbox{\hskip-7mm $d$ \eightrm(discrepancy)}\vskip3pt\hbox{0}\vskip4mm
\hbox{1}\vskip4.4mm \hbox{2}\vskip4mm \hbox{3}\vskip4.5mm \hbox{4}\vskip4.5mm \hbox{5}\vskip4mm \hbox{6}
}}
\figure{Adjacency graphs {\tenrm \ (vertices of type I are boxed).}}\label{adjacency}
\endfigure
\endinsert

Note that these two examples demonstrate two fundamental properties of the graphs $\G_{d,k}$ that hold for many values of $d$ and $k$:
\roster\item
{\it Smith adjacency}: the endpoints of the edges have Smith discrepancies that differ by 1; in particular, there is no loop-edges;
\item
{\it Klein principle}: if one of the endpoints of an edge has Klein type I, then its Smith discrepancy is smaller than for the other endpoint.
\endroster

\subsection{Quadrocubics representing edges of $\G_{3,3}$ \cite{Kr-rigid}}\label{quadrocubics-edges}
Consider a cubic threefold $X\subset P^4$ with a singular point $s$, and choose coordinates $x,y,z,u,v$ in $P^4$ so that $s$
acquires the coordinates $(0,0,0,0,1)$. Then $X$ is defined by equation $f_3+vf_2=0$
where $f_2$, $f_3$ are homogeneous polynomials in $x,y,z,u$ of degree $2$ and $3$ respectively.

The intersection $A$
of the quadric $f_2=0$ with the cubic $f_3=0$ in the local projective $3$-space, $P^3_s$, centered at $s$
describes the locus of lines contained in $X$ and passing through $s$.
We call $A$ the {\it quadrocubic associated with $s$}.
If $X$ and $\ss$ are real, then $A$ is real also, so that
this correspondence provides a well defined map from the space of real projective classes of real pairs $(X,\ss)$ where $X$ is a real cubic threefold with a chosen real node $\ss$
to the space of real  projective classes of real pairs $(Q, A)$ where $Q$ is a real non-singular quadric surface in $P^3$ and $A$ is a real curve traced on $Q$ by a real  cubic
surface.

\proposition\label{cubics-quadrocubics-association} The map defined above is an isomorphism.
\endproposition

\lemma\label{nodes-to-quadrocubics}
Assume that real cubic threefolds $X_0$ and $X_1$ have a common node
$\ss_0$ and
that at this node both cubics have the same local quadric  and the same quadrocubic. Then there exists a
continuous family of
real projective transformations
$T_t\: P^4 \to P^4$, $t\in[0,1]$, such that $T_0=\id$,  $T_1(X_0)=X_1$,
and $T_t(\ss_0)=\ss_0$ for each $t\in[0,1]$.
\endlemma

\proof
Consider a real affine chart centered at $\ss_0$. The degree 3 real affine equations of $X_0$ and $X_1$
can be chosen then in the form
$f_2+f_3=0$ for $X_0$ and $f_2+(cf_3+f_1f_2)=0$ for $X_1$, where
$c\ne0$ is a real constant and $f_1,f_2, f_3$ are real homogeneous equations in 4 variables
of degrees 1, 2, 3 respectively.
Affine homotheties include $X_0$ into
a continuous
family of cubics $f_2+c f_3=0$,
with $c>0$. Alternation of the sign of $c$ is equivalent
to a simultaneous alternation of sign for
all the affine coordinates,
which can be
achieved by a continuous family of rotations with center in $\ss_0$.
So, it is left to eliminate the term $f_1f_3$ by replacing the affine chart  in a way
that  $1 + f_1=0$ becomes ``the infinity hyperplane'', which can be
done via a continuous family of real projective transformations.
\endproof

\demo{Proof of Proposition \ref {cubics-quadrocubics-association}}
Straightforward consequence of Lemma \ref{nodes-to-quadrocubics}.
\qed\enddemo

If $X$ is a one-nodal cubic, then the quadrocubic $A$
 is non-singular, it has genus $4$, and its embedding $A\subset P^3_s$ is canonical.
Conversely, for any complex non-singular not hyperelliptic genus $4$ projective curve $A$, its canonical embedding presents $A$ as a quadrocubic in $P_s^3$
(with not necessary degenerate $f_2$).
If $A$ is real, then the latter embedding
can be made real and gives an intersection of a real quadric $f_2=0$ with a real cubic $f_3=0$.

This observation implies that the edges of $\G_{3,3}$ are in one-to-one correspondence with
the deformation classes of real genus $4$ projective curves (note that in genus 4 hyperelliptic ones form a codimension $2$
subvariety, and therefore their removal does not change
the set of real deformation classes).

The deformation classification of real non-singular
genus $g$ projective curves is well-known: the deformation class of such a curve
 is detected by the number $0\le k\le g+1$
of components of its real locus and by its Klein type.
In the case of $g=4$, the Klein type is determined by
the number of components for $k\ne1,3$: for $k=5$ the type has to be $I$ (as it is always for M-curves),
and for even $k$ the type is $II$ (as it is always if the number of components has the same parity as the genus of
the
 curve). For $k=1,3$ both Klein types are realizable, and so,
 we have 8 deformation classes of real genus $4$ real curves in total;
they  will be denoted by {\bf k} and {\bf k}${}_I$ according to the number $k$ of real components,no
where subscript ``$I$'' will be used only for $k=1,3$ for the curves of Klein type I.

\theorem\label{quadrocubics_as_edges}\cite{Kr-rigid}
The edge in $\G_{3,3}$ represented by a quadrocubic of type {\bf k}, $1\le k\le 5$ connects the vertices $C^{k-1}$ and $C^k$. The edge of type {\bf 0} connects
$C^0$ with $C^1_{I(2)}$.
The edges of type {\bf k}${}_I$, $k=1,3$, connect $C^{k-1}$ with $C^k_I$
{\rm (}see Figure \ref {3fold-graph}{\rm )}.
\qed\endtheorem

\vskip-5mm
\midinsert
\figure{
Graph $\G_{3,3}$ with indication of types of edges (marked by
{\bf k} and {\bf k}${}_I$)  and vertices (type one classes are shown in black).}
\label{3fold-graph}\endfigure
\hskip30mm\epsfbox{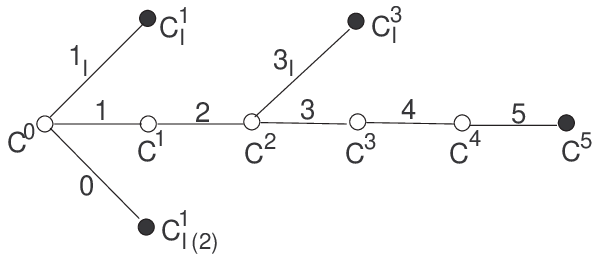}
\endinsert

Among the two endpoints of the edge of $\G_{3,3}$, the one with the greater value of the Smith discrepancy
will be called {\it ascending} and the one with the lesser value {\it descending}.
\corollary\label{Smith-Klein-preserving}
 If $A$ is a non-singular quadrocubic with $A_\R\ne\oo$ and $X$ a real cubic threefold
 representing the ascending endpoint of the edge of $\G_{3,3}$ given by $A$, then
the Smith discrepancy and the Klein type of $A$ are
the same as for $X$.
\qed\endcorollary

\subsection{Spectral matching correspondence}
By a {\it matching} we mean
a pair of deformation classes $([X],[S])$, where $[X]\in V_{3,3}$ is represented
by a non-singular cubic threefold $X$ and $[S]\in V_{5,1}$ by a non-singular plane quintic $S$.
 A matching $([X],[S])$ is called {\it spectral}, if $S$ is the spectral quintic with respect to a
 real line $l\subset X$ such that $(X,l)\in\CC_\R^*$,
in this case we say that the line $l$ {\it realizes the matching} $([X],[S])$.

A  matching $([X],[S])$ is
called {\it perfect}, if $X$ and $S$ have the same Klein type and same Smith discrepancy, $d_X=d_S$.
If $d_X=d_S-2$, the
matching $([X],[S])$ is called {\it skew}. We say that a skew matching
is {\it type-preserving}, if both $X$ and $S$
are of the same Klein type, and {\it mixed} otherwise.
A mixed skew matching is called {\it admissible}, if $X$
is of Klein type $I$, $X_\R$ is connected, and $S$
is of Klein type $II$.

\theorem\label{main}
$(1)$ The spectral matching $([X],[S])$ realized by a line $l\subset X$ is a real deformation invariant
of the pair $(X,l)$.

$(2)$ This invariant is complete in the sense that
all the pairs $(X,l)$ that give the same spectral matching belong to the same
real deformation class.

$(3)$ A matching is spectral if and only if it belongs to one of 3 types:
perfect, type-preserving skew, or admissible mixed skew.
\endtheorem

Comparison of the graphs $\G_{3,3}$ and $\G_{5,1}$ (see Figure \ref{adjacency}) implies immediately that for any
$[X]\in V_{3,3}$ there exists precisely one $[S]\in V_{5,1}$ such that $([X],[S])$ is a perfect matching.
In particular, there exist precisely $|V_{3,3}|=9$ such matchings.
Note also that there is no skew matching involving $[X]$ of type $C^0$ (because $d_S\le5$,
while $d_X+2=7$),
and that in the case of type $C^1_{I(2)}$
the only skew matching is not admissible.
In the cases of types $C^5$ and $C^3_I$ there exist two skew matchings: one is type-preserving and another is admissible mixed.
The remaining $5$ deformation classes $[X]$ admit only one skew matching, which is admissible mixed in the case of $C^1_I$ and
type-preserving in the other 4 cases.
So, by Theorem \ref{main} there exist totally $9+4+5=18$ spectral matchings and we obtain the following result.

\corollary\label{18classes}
There exist precisely 18
real deformation classes of pairs $(X,l)$, where $X$ is a real non-singular cubic threefold and $l\subset X$
is a real line. Namely, they are represented by:
\roster
\item 9 perfect matchings
$$(C^1_{I(2)}, J\+1\la1\ra ),  (C^1_{I}, J\+1\la1\ra ), (C^3_{I}, J\+4_I ), \,\,\text{and} \,\,(C^{k-1}, J\+k)\,\,\text{ with} \,\,1\le k\le 6; $$
\item 3 admissible mixed skew  matchings
$$(C^1_{I}, J), (C^3_{I}, J\+2 ), \,\,\text{and} \,\,(C^{5}, J\+4); $$
\item and 6 type-preserving skew matchings
$$(C^3_I, J\+1\la1\ra ),  (C^5,  J\+4_I ),  \,\,\text{and} \,\,(C^{k+1}, J\+k)\,\,\text{ with} \,\,0\le k\le 3. \quad \qed $$
\endroster

\endcorollary

\subsection{Proof of Theorem \ref{main}(1)}\label{scheme-main}
Proposition \ref{real-spectral-correspondence} implies that $[S]$ is an invariant of pairs $(X,l)\in\CC_\R^*$ with respect to deformations
in the (larger) space $\CC_\R$,
and so the pair $([X],[S])$ is also such an invariant.

\subsection{Strategy of proving (2)--(3) of Theorem \ref{main}}
Claim (3) of Theorem \ref{main}, which
describes the image of the {\it spectral matching correspondence}
that associates
a  pair $([X],[S])$ to a real deformation class $[X,l]$ of $(X,l)$,  is proved
in Subsection \ref{surjectivity of SMC}.

Claim (2) is proved in Subsection \ref{proof-main-easy-case} for cubics with $d_X\ge4$, while
 the more involved case $d_X\le3$
is derived in Section \ref{summary}
from Theorems \ref{N-correspondence} and \ref{T-orbits}.
In fact, the latter theorems
enumerate the deformation classes $[X,l]$ in terms of orbits of the monodromy action on  real Fano components
and provide, thus, an upper bound on the number of deformation classes. This bound happens to coincide with the lower bound given by
claim (3), which ends the proof of claim (2).

\subsection{Smith discrepancy in terms of real Fano components}
The remaining part of this section is devoted to formulation of two Theorems which develop our main Theorem \ref{main}.

As is known (see Subsection \ref{real_Fano}),  the real locus of the Fano surface of any
non-singular real cubic threefold $X$ has
precisely one connected component with
odd Euler characteristic.
Namely, this component is homeomorphic to $\Rp2$ in the case of an exotic cubic and to a connected sum $\#_5\Rp2$
in the other cases. Let us denote this component by $N(X)$ and call it the {\it odd Fano component}.

\theorem\label{N-correspondence}
Assume that $X$ is a non-singular real cubic threefold, and $l\in F_\R(X)\sm \Gres$.
Then:
\roster\item
If $l\in N(X)$, then $\spec$
has the same Smith discrepancy and the same Klein type as $X$ {\rm (}cf., Table 2 {\rm ).}
\item
If $l\notin N(X)$, then the Smith discrepancies $d_\spec$ of $\spec$ and $d_X$ of $X$ differ by 2:
$$d_\spec=d_X+2.$$
\endroster
In other words,
the perfect spectral matchings $(X,S)$ are realized by
$l\in N(X)$, while the skew matchings by $l\in F_\R(X)\sm N(X)$.
\endtheorem

$$
\resizebox{130mm}{12mm}{
\boxed{\matrix\text{Table 2. $X_\R$, $F(X)$ and $S_\R$ in the case of $l\in N(X)$,
\ here $H=S^1\times S^2$}\\
\text{-----------------------------------------------------------------------------------------------------------------------}\\
\matrix
{} & C^1_I&C^{1}_{I(2)}&C^0&C^1&C^2&C^3 \text{ or } C^3_I&C^4&C^5\\
X_\R &\S(2,4,6)& \Rp3\+S^3 &\Rp3& \Rp3\#H&\Rp3\#2H
&\Rp3\#3H&\Rp3\#4H&\Rp3\#5H\\
F_\R  &\Rp2\+N_6 &N_5 &N_5 &N_5\+T^2&N_5\+_3T^2&N_5\+_6T^2&N_5\+_{10}T^2&N_5\+_{15}T^2\\
S_\R &J\+1\la1\ra&J\+1\la1\ra &J\+1&J\+2   &J\+3    &J\+4    &J\+5     &J\+6\\
\endmatrix
\endmatrix}
}
$$

\subsection{The monodromy}\label{monodromy_groups}
Having chosen a connected component of the variety parametrizing real non-singular cubic threefolds and a particular cubic
$X$ as a base point in it, we get a well defined monodromy action homomorphism, $\mu_X$, from the
fundamental group of the  chosen connected component into the mapping class group $\Map(X_\R)$.
Besides, it yields a homomorphism, $\mu_F$, into
the mapping class group $\Map(F_\R(X))$, and, in particular, defines an action
of the fundamental group on the set of connected components of $F_\R(X)$. The image of the latter action is a subgroup $\mon$
of the full permutation group of the components of $F_\R(X)$,
and we call it  the {\it Fano real component monodromy group}.
We say that two components of $F_\R(X)$ are {\it monodromy equivalent} if they can be permuted by $\mon$.

Let us denote by $\TT(X)$ the set of toric components of $F_\R(X)$.
The next theorem describes the orbits of
the action of $\mon$ on the set of toric components, if their number
$|\TT(X)|$ is $\ge 2$, which is the case for $X$ of
deformation types $C^k$, $k\ge2$, and $C^3_I$.
By Theorem \ref{N-correspondence}, a choice of a line $l$
in $(F_\R(X)\sm N(X))\sm\Gres$
provides a mixed matching.
So, for cubics $X$ of Klein type I, namely
for those of deformation types $C^5$ and $C^3_I$
(among the ones with $|\TT(X)|\ge 2$), we may speak on
two types of toric components: a choice of $l$ on a {\it torus of type I} provides
a type-preserving skew matching (that is a spectral quintic $S$ of type I),
while choice on a {\it torus of type II} yields a mixed skew matching ($S$ is of type II).
These two types of tori form a partition of the set $\TT(X)$ into
two subsets,
$\TT_I(X)$ and $\TT_{II}(X)$.

After completing proving Theorem \ref{N-correspondence} in Section \ref{summary}, we finalize there
our proof of the following Theorem \ref{T-orbits}.

\theorem\label{T-orbits}
Let $X$ be a non-singular real cubic threefold. Then:
\roster\item
$\mon$ acts transitively on $\TT(X)$ unless $X$ has deformation type $C^5$ or $C^3_I$;
\item
$\mon$
has two orbits, $\TT_I(X)$ and $\TT_{II}(X)$, in $\TT(X)$ if $X$ is of deformation type $C^5$;
$\TT_I(X)$ is formed
by $6$ tori, and $\TT_{II}(X)$ by $9$ tori;
\item
$\mon$ has
two orbits, $\TT_I(X)$ and $\TT_{II}(X)$, if $X$ is of deformation type $C^3_I$; each orbit
consists of $3$ toric components.
\endroster
\endtheorem

Table 3 indicate explicitly the correspondence between the monodromy orbits of the tori in $F_\R(X)$
and the real deformation type of the corresponding spectral quintics (when the quintics are non-singular).
$$
\boxed{\matrix\text{Table 3. $S_\R$ for $l\in F_\R(X)\sm N(X)$}\\
\text{------------------------------------------------------------------------------------------------------}\\
\matrix
X & C^1_I&C^1&C^2&C^3_I&C^3_I&C^3&C^4&C^5&C^5\\
l &N_6 &\TT&\TT&\TT_I&\TT_{II}&
\TT&\TT&\TT_I&\TT_{II}\\
S_\R & J_{h_I}& J_{h_{II}}&J\+1& J\+1\la1\ra&J\+2_{h_I} &J\+2_{h_{II}}   &J\+3    &J\+4_I&J\+4_{II} \\
\endmatrix\endmatrix}
$$
In this table symbols $h_I$ and $h_{II}$ indicate the types of the difference class $h=\theta-\theta_{0}\in H^1(\spec;\Z/2)$: type I if $h$
is the characteristic class $w_c$ of involution $c$, and type II otherwise.


\section{Comparison of $X_\R$ and $\spec_\R$}\label{Prymian-sec}

\subsection{Homological correspondence }\label{HC}
As before, we consider a  real non-singular cubic threefold $X\subset P^4$,
a real line $l\subset X$, $l\notin\Gres$, the associated
spectral curve (a real  non-singular plane quintic) $\spec$ and the theta-conic $\Theta$.
Recall that for the spectral double covering
$p_S\:\spect\to\spec$ (see Subsection \ref{theta-conics})
 the curve $\spect$ is
formed by the lines intersecting with $l$, and that this covering is defined by
the {\it spectral difference class} $h^*$ that is the difference class associated with $\Theta$
(see Subsection \ref{def-spectral-theta}).

In what follows, we denote by $h\in H_1(\spec;\Z/2)$ the Poincar\'e dual to $h^*$ and  by $H_1^h(\spec;\Z/2)$ the orthogonal complement of $h$. We equip the quotient $\Z/2$-vector space
$H_1^h(\spec;\Z/2)/h$ with the induced by the intersection pairing non-degenerate bilinear form.

The involution induced  in $H_1(\spect)$
by the deck transformation of the theta-covering is denoted by
$\tau$
and its eigen-lattice $\Ker (1+\tau)$
by $H_1^-(\spect)$.
The involutions induced by the complex conjugation in $H_3(X;\Z/2)$, $H_1^h(\spec;\Z/2)/h$, $H_1(\spect)$, and $H_1^-(\spect)$ are denoted by $c_X$, $c_\spec^h$, $\cs$,
and $\cs^-$, respectively.

\theorem\label{homology_correspondence} There exists an isometry
$H_3(X;\Z/2)\to
H_1^h(\spec;\Z/2)/h$ which commutes with the involutions $c_X$ and
$c_\spec^h$.
\endtheorem

We postpone the proof of this theorem to the end of this subsection.

Let us start from considering the following piece of the Gysin exact sequence
$$H_2(\spec;\Z/2)@>\cap h^* >>H_1(\spec;\Z/2)@>p_S^*>>H_1(\spect;\Z/2)@>p_{S*}>>H_1(\spec;\Z/2)
@>\cap h^*>>H_0(\spec;\Z/2)$$
where the first and the last group are $\Z/2$, the first
map sends the
generator of $\Z/2$ to $h$, and the last map sends $a\in
H_1(\spect;\Z/2)$ to
$h^*(a)\in\Z/2$.
In addition, consider an {\it anti-symmetrization homomorphism}
$1-\trans\:H_1(\spect)\to H_1(\spect)$, $x\mapsto x-\trans x$.
 We denote by $\t_2$ the modulo $2$ reduction of $\t$ and
 put $B=\Im(1+\t_2)\subset H_1(\spect;\Z/2)$.

\lemma\label{h_versus_Prym}
\roster
\item $\Im (p_S^*\circ p_{S*})=B$;
\item $p_S^*$ induces an
isomorphism $H^h_1(\spec;\Z/2)/h\to B$;
\item  $H_1^-(\spect)=(1-\t)H_1(\spect)$;
\item $B$ is the image of
$H_1^-(\spect)$ under the modulo 2
reduction homomorphism
$H_1(\spect)\to H_1(\spect;\Z/2)$;
\item the latter
reduction identifies $H_1^-(\spect)\otimes \Z/2$
 with $B$;
\item all the above isomorphisms are equivariant with respect to
complex conjugation involutions.
\endroster
\endlemma
\proof Item (1) follows from the identity $1+\t_2=p_S^*\circ p_{S*}$. The same identity and the exactness of the Gysin sequence at its second and
forth terms implies item (2). Items (3) - (5) follow from $ (1-\t)H_1(\spect ) = B \mod 2$ and $\rank B=
\rank H_1^-(\spect) $. Item (6)  becomes now evident.
\endproof

Due to item (3) in above Lemma, $\frac 12 (a\cdot b)$ is an integer for all $a,b\in H_1^-(\spect)$. Namely, it is equal to $-x\cdot y+ x\cdot \t y$ where $x-\t x=a, y-\t y=b$.
Thus, taking the residue modulo 2 of $\frac 12(a\cdot b)$ we obtain a well defined symmetric bilinear form on $H_1^-(\spect)\otimes \Z/2$. We call it the {\it Prym pairing}.
On the other hand, $x\cdot (y+\t_2y)=(x+\t_2x)\cdot y$ for all $x,y\in B$ and gives a well defined symmetric bilinear form on $B$, which we call the {\it $B$-pairing}.

\proposition\label{isometry} The inclusion homomorphism
$H_1^-(\spect)\otimes \Z/2\to H_1(\spect;\Z/2)$ transforms the Prym pairing in $H_1^-(\spect)\otimes \Z/2$ into the $B$-pairing in $B$. The pull-back homomorphism $p_S^*$
transforms the intersection index pairing in $H^h_1(\spec;\Z/2)/h$ into $B$-pairing in $B$.
\endproposition

\proof The first part is straightforward from definitions and Lemma \ref{h_versus_Prym}(5). The second part follows from Lemma \ref{h_versus_Prym}(1-2) and
the identities $x\cdot (y+\t_2y)=x\cdot (p_S^*\circ p_{S*})y= (p_{S*}x) \cdot (p_{S*}y)$ for all $x,y\in H_1(\spect;\Z/2)$.
\endproof

Let us recall, finally, that given a 1-cycle in $\spect$ represented by a loop
$\gamma\:[0,1]\to\spect$,
one can associate with it a 3-cycle traced in $X$
by the lines $\ell(\gamma(t))$, $t\in[0,1]$. This yields a
so-called {\it Abel-Jacobi homomorphism} $\Psi\:H_1(\spect)\to H_3(X)$.

\proposition\label{Abel-Jacobi}
Abel-Jacobi homomorphism $\Psi$ has the following properties:
 \roster \item $\Psi(a)\cdot
\Psi(b)=a\cdot(\t_{\spect}(b)-b)=(\t_{\spect}(a)-a)\cdot b$;
\item $\Ker\Psi=H_1^+(\spect)$;
 \item $\Im\Psi=2H_3(X)$;
 \item
$\frac12\Psi$ yields an isometry between
$H_1^-(\spect)$
equipped with the pairing $(x,y)\mapsto -\frac12 x\cdot y$ and
$H_3(X)$
equipped with the Poincar\'e pairing. This isometry is equivariant with respect to $\cs^-$ and $c_X$.
\endroster
\endproposition

\proof For items (1) and (2), see, for example, \cite{Tyu} Lecture 2, \textsection 2. Items (3) and  (4) follow from
(1) and (2). Indeed, since  each element
in  $H_1^-(\spect)$ is of the form $x-\tau x, x\in H_1(\spect)$,
(2) implies that $\Psi(H_1^-(\spect))\subset 2H_3(X)$;  finally (1), due to the latter inclusion, implies both (3) and the isometry property.
Equivarience is tautological.
\endproof

\demo{Proof of Theorem \ref{homology_correspondence}}
It is left to combine the isometry from Proposition \ref{Abel-Jacobi}(4)
with the isometries from Proposition \ref{isometry}.
\qed\enddemo

\subsection{Vanishing classes and bridges}
Given a continuous family $\spec_t$, $t\in[0,1]$, of real plane curves with non-singular $\spec_t$ for $t\in(0,1]$ and $\spec_0$ having only one node as a singularity,
we say that $\spec_1$ is a {\it real perturbation} of $\spec_0$ and $\spec_0$ is a {\it one-nodal real degeneration} of $\spec_1$.
In this case, the so-called {\it vanishing cycles} can be represented by simple closed curves $\g_t\subset S_t, t>0,$
invariant with respect to the real structure,
and we call them
{\it real vanishing cycles}. If the node is a solitary point, then one kind of its perturbations leads to a birth of an oval
and this oval does represent
a real vanishing cycle. Another kind of perturbations leads to ``disappearing'' of the solitary point: in this case
$\g_t$ have no real points, and the action of $\conj$ on them is modeled by an antipodal map on $S^1$.
The homology class $v_t\in H_1(\spec_t;\Z/2), t>0,$ represented by such a $\g_t $ is called a {\it vanishing oval class} for the first kind of perturbations
and a {\it vanishing invisible class} for the second kind. Clearly, in the both cases $c=\conj_*$ acts identically on $v$.

A closed curve on $\g\subset\spec_t$
is called  a {\it bridge} if it  is $\conj$-invariant and has precisely two
real points, or in other words $\g=\b\cup \conj(\b)$ where $\b$ is a path between points of $S_{t\R}$
that passes in $S_t\sm S_{t\R}$. Their homology classes $v\in H_1(S_t;\Z/2)$, which are also $c$-invariant, are called {\it bridge classes}.
A special kind of bridges
appear after a perturbation of a cross-node on $S_0$.
These classes are called {\it vanishing bridge-classes}.

\subsection{Subgroups $I$, $K$ and the characteristic element $w_c$}
Let $K$ and $I$ denote the kernel and the image of the linear endomorphism $\phi_c=1+c: V\to V$, $V=H_1(S;\Z/2)$.
Note that $\phi_c^2=0$, or equivalently $I\subset K$.

\proposition\label{K-generators} For any real non-singular algebraic curve $S$, the subgroup
$K\subset H_1(S;\Z/2)$ is
spanned by the oval-classes and the bridge-classes, provided $S_\R\ne\oo$.
\endproposition

\demo{Proof} It is equivalent to proving that
$K$ is the image of the Smith homomorphism
$$ H_1(S/\conj,S_\R;\Z/2)\oplus H_1(S_\R;\Z/2)\to H_1(S;\Z/2).$$
The latter is a straightforward well-known consequence of the Smith exact sequence. \qed
\enddemo

We denote by $w_c\in V$ the characteristic class of the
quadratic form $x\mapsto x\cdot cx$, and, for any class $h\in V$,
denote by $h^*\in V^*=H^1(S;\Z/2)$
its Poincar\'e dual
({\it cf.}, Subsections \ref{Klein-type} and \ref{theta-conics-and-characteristics}).
We put $V_0=\{x\in V\,|\, x \cdot cx=0\}$, $I_0=\phi_c(V_0)$,
$I_1=I\sm I_0$, and $d=\rank I$. We call the latter the {\it discrepancy of $c$} (cf., Proposition \ref{quintic_discrepancy}).

\proposition\label{symmetrization-characteristic}
\roster\item
$K=\{h\in V\,|\,h\cdot x=0\,\,
\text{for every}\, \,\, x\in I\}$, and so $K$ and $I$ are mutual annihilators.
\item
$w_c=[S_\R]\in I$
and $K\subset V_0$.
\item
If $w_c\ne0$, then $I_0$ is an index 2 subgroup of $I$.
\item
$x\cdot cx=y\cdot cy$ if $\phi_c(x)=\phi_c(y)$.
\item
$w_c$ belongs to $I_0$ if the discrepancy $d=\rank I$ is even and to $I_1$ if $d$ is odd.
\endroster
\endproposition

\demo{Proof} Since $(h+ch)\cdot y= h\cdot (y+cy)$, the validity of $h\cdot x=0$ for all $x=y+cy\in V$ is equivalent to $x=cx$, that is $x\in K$, which proves claim (1).

The first part of claim (2) is a particular case of Arnold lemma (cf., Subsection \ref{Klein-type}).
Since $x=cx$ implies that $x\cdot w_c=x\cdot cx=x^2=0$, we see that $K\subset V_0$ and $w_c$ annihilates $K$. Due to claim (1), the latter implies in its turn,  that $w_c\in I$.

Since $K\subset V_0$ and thus $I_0\cong V_0/K$, we obtain claim (3).

Claim (4) follows from $x\cdot cx= x\cdot (x+cx)=x\cdot (y+cy)= (x+cx)\cdot y=(y+cy)\cdot y=y\cdot cy$.

For proving claim (5),
we introduce a {\it residual pairing in $I$} between $a=x+cx$ and $b=y+cy$ by the formula
$\la a,b\ra_I=x\cdot b$. It is well defined by (1), symmetric due to the identity $x\cdot (y+ cy) =(x+cx)\cdot y$,
and clearly non-degenerate.
Its characteristic element is $w_c$ because $w_c\in I$ by (2) and
$\la a,a\ra_I=x\cdot a=x\cdot cx=w_c\cdot x=\la w_c,a\ra_I.$
It remains to notice that the
square of any characterstic element of a non-degenerate $\Z/2$-valued inner product equals to the rank
of the underlying $\Z/2$-vector space,  to represent $w_c\in I$ as  $w_c=z+cz$, and to conclude that
$
\rank(I)=\la w_c,w_c\ra_I=z\cdot w_c=z\cdot cz\mod 2.$
\qed \enddemo

\subsection{Comparison of discrepancies $d$ and $d'$}\label{Smith-comparison}
For $h\in K$, we consider $V'=V^h/h$, equip it  with the induced, nondegenerate, pairing and induced involution $c'\:V'\to V'$,
and put $K'=\Ker(1+c')$, $I=\Im (1+c')$, $d'=\rank I'$.

\proposition\label{change-of-defect}
\roster\item
If $h\in I_0\sm \{0\}$, then $d=d'+2$;
\item
if $h\in I_1\sm \{0\}$, then $d=d'+1$;
\item
if $h\in K\sm I$, then $d=d'$.
\endroster
\endproposition

\demo{Proof}
Consider the restriction $\phi_c|_{V^h}\:V^h\to V^h$.
If $h\in I\sm \{0\}$, then $K\subset V^h$, and
$\phi_c(V^h)\subset I=\phi_c(V)$ is a codimension one subspace since $V^h\subset V$ is.
In the case (1), we get $h\in\phi_c(V^h)$, so that
$$\dim I'=\dim\phi_c(V^h)/h=\dim\phi_c(V^h)-1=\dim I-2.$$
 In the case (2), $h\notin\phi_c(V^h)$ and the projection $V^h\to V^h/h$ is
 injective
on $\phi_c(V^h)$, so that
$$\dim I'=\dim\phi_{c'}(V^h/h)
=\dim\phi_c(V^h)=\dim I-1.$$
Next,  in the case (3),
Proposition \ref{symmetrization-characteristic}(1)
implies that $K\cap V^h\ne K$, and thus,
$V^h/(K\cap V^h)=V/K$ or, in other words,
$\phi_c(V^h)=\phi_c(V)=I$.
Furthermore, the quotient map $V^h\to V^h/h$ is injective on
$I\cap V^h$, since $h\notin I$, and therefore
$$\dim I'=\dim\phi_c(V^h)=\dim I.
\qed$$
\enddemo
\subsection{The case of real plane quintics}
Here, and for the rest of this Section, $S$ is a real non-singular quintic and $q_0$ its Rokhlin quadratic function.
We let also  $K_i=\{x\in K\,|\, q_0(x)=i \}$, $i=0,1$.

\proposition\label{real-quadratic-functions}
\roster\item
 A quadratic function $q\:H_1(X;\Z/2) \to\Z/2$ is real
 if and only if its difference class $h^*=q-q_0 \in H^1(\spec;\Z/2)$ is real, that is Poincar\'e-dual to some $h\in K$.
\item
All real quadratic functions $q\:H_1(X;\Z/2)\to\Z/2$ have the same restriction to $I$;
this restriction is linear and
$q(x)=i$ for any $x\in I_i$, $i=0,1$. In particular, $I_i=K_i\cap I$, $i=0,1$.
\item
The Arf-invariant $\Arf(q_0+h^*)\in\Z/2$ is equal to $q_0(h)+1$. In particular,
$$\Arf(q_0+h^*)=\cases 1 &\text { if }  h\in K_0, \\
0 &\text{ if } h\in K_1.
\endcases$$
\endroster
\endproposition

\proof
Proof of Item (1) is straightforward from definitions.

The identity $q\vert_I=(q_0+h^*)|_I=q_0|_I$ holds for any $h\in K$
because
$h^*(x)=h\cdot x=0$ for any
$x\in I$, as it follows from
Proposition \ref{symmetrization-characteristic}(1).
The other part of item (2) follows from
$
q(x+cx)=q(x)+q(cx)+x\cdot cx=x\cdot cx\mod2.
$

Item (3) follows from two well known facts: a general formula
$$\Arf(q+h^*)=\Arf(q)+q(h), \quad q\in\Quadr(S), h\in H_1(S;\Z/2)$$
that we apply to $q=q_0$,  and
the fact that $\Arf(q_0)=1$ in the case of plane quintics.
\endproof

\corollary\label{discrepancy-relation} If $S$ is the spectral curve of a non-singular real cubic threefold $X$ with respect to a real line
$l\in F_\R(X)\sm\Gres$, then,
for the Smith discrepancies
$d_X$ and $d_\spec$ of $X$ and $\spec$, respectively,
we have two options,
depending on the spectral difference class $h^*\in H^1(S;\Z/2)$ dual to $h\in H_1(S;\Z/2)$:
\roster\item
$d_\spec=d_X+2$, if $h\in I$, that is if $h=x+cx$ for some $x\in H_1(S;\Z/2)\sm K$;
\item
$d_\spec=d_X$ otherwise, that is for any $h\in K\sm I$.
\endroster
\endcorollary

\proof Recall, that according to Proposition \ref{quintic_discrepancy} the Smith discrepancies $d_X$ and $d_S$ are equal to the rank of the homomorphism $\id +c_*$
in $H_3(X;\Z/2)$ and $H_1(S; \Z/2)$, respectively. Note also that $h$ is always non zero.
By Proposition \ref{real-quadratic-functions}(1), the quadratic function $q$ associated with the spectral theta-characteristics of $(X,l)$
is of the form
$q=q_0+h^*$, where $h\in K$. It remains to apply Proposition \ref{change-of-defect}
and to notice that its second case, $h\in I_1$, is impossible, since otherwise, by Proposition
\ref{real-quadratic-functions}(3), we get $\Arf(q)=0$, while the spectral theta-characteristics are odd, so the Arf-invariant should be equal to 1.
\endproof

\subsection{Comparison of Klein types}\label{Klein-types}
\lemma\label{Klein-type-relation}
In notation of Corollary \ref{discrepancy-relation},
\roster\item
if $d_X=d_S$, then the Klein types of $X$ and $S$ are the same;
\item
if $d_X=d_S-2$ and $S$ is of Klein type I, then $X$ is of  Klein type I too.
\endroster
\endlemma

\proof According to Corollary  \ref{arnold_typeI},
being of Klein type I for a non-singular real quintic $S$
(respectively, non-singular real cubic threefold $X$) is equivalent to
vanishing of the characteristic element
$w_c(S)\in H_1(S;\Z/2)$ (respectively, $w_c(X)\in H_3(X;\Z/2)$).
On the other hand, Theorem
\ref{homology_correspondence} says that
$H_3(X;\Z/2)$ and $H_1^h(\spec;\Z/2)/h$ are equivariantly isometric,
which implies  that $w_c(X)$ is the image of $w_c(S)$ under this isometry.
So, vanishing of $w_c(S)$ implies vanishing of $w_c(X)$. In the other direction,
if $w_c(X)=0$ then $w_c(S)\in\{0,h\}$, and since $w_c(S)\in I=\Im(1+c_S)$ (see Proposition \ref{symmetrization-characteristic}
(2)), the Klein types
of $S$ and $X$ must be the same if $d_X=d_S$ (since it is equivalent to $h\in K\sm I$, see Proposition \ref{discrepancy-relation}).
\endproof

\proposition\label{skew-matchings-and-wc} Let $l\subset X$ provide
a spectral skew matching $([X], [S])$,
and let $q=q_0+h^*$ with $h\in H_1(S;\Z/2)$ be the quadratic function
representing the theta-characteristic associated with the pair $(X,l)$.
Then this matching is  admissible mixed if $h=w_c$ and
type-preserving if $h\ne w_c$.
\endproposition

\demo{Proof}
According to Proposition  \ref{arnold_typeI}, $S$ is of Klein type I if and only if
$w_c=0$. Similarly, using in addition Theorem \ref{homology_correspondence}, we conclude that
$X$ is of the Klein type I if and only if  the image of $w_c\in V=H_1(S;\Z/2)$ in the quotient
$V^h/h\cong H_3(X;\Z/2)$
vanishes. The latter happens if either $w_c=0$ or $w_c= h$.
Recall also that $h\ne 0$, see Proposition \ref{difference-class-covering}.
Thus, we have 3 cases:
 \roster\item
 $w_c=0$, then both $S$ and $X$ have Klein type I and thus, the matching is type-preserving;
\item
$w_c\ne0$ and $h\ne w_c$, then both $S$ and $X$ have Klein type II and thus, the matching is also type-preserving;
\item
$w_c\ne0$ and $h=w_c$, then $S$ has type II but $X$ has type I and the matching is admissible mixed.
\qed\endroster
\enddemo

\corollary\label{bottom-skew-matchings}
The only spectral skew matchings $([X],[S])$ with $[S]=J$
are with $[X]=C^1$ and $[X]=C^1_I$. Namely,
\roster\item
$[X]=C^1_I$ if the spectral difference class $h^*$ is dual to
$h=w_c$,
\item
$[X]=C^1$ if $h\ne w_c$.
\endroster
In particular, there is no  spectral matching $([X],[S])$ with $[S]=J$ and $[X]=C^1_{I(2)}$.
\endcorollary

\proof
For $[S]=J$
we have $d_S=6$, while by
Lemma \ref{Klein-type-relation}
$d_X$ must be either $6$ or $4$. But real cubic threefolds $X$ with $d_X=6$ do not exist,
and $X$ with $d_X=4$ must be of
type $C^1$, $C^1_I$, or $C^1_{I(2)}$ (see Figure \ref{adjacency}).

To exclude the latter, let us consider
the conic bundle
$\pi_l\:X_l\to P^2$, see Subsection \ref{conic-bundles-for-cubics}.
For $X$ of type $C^1_{I(2)}$ we have $X_\R=S^3\+\Rp3$,
and, thus, the real locus $(X_l )_\R$ is also disconnected.
On the other hand, the restriction $\pi_l|_{(X_l )_\R}$ projects $(X_l )_\R$ to
$\Rp2$ with connected fibers, they are real conics, and a connected critical locus, that is
the only real component of $S$, which is impossible for a disconnected $(X_l )_\R$.

So, it remains to consider the skew matchings $(C^1,J)$ and $(C^1_I,J)$.
By Proposition \ref{symmetrization-characteristic}(5), $w_c\in I_0$, which together
with Proposition \ref{real-quadratic-functions} implies that $\Arf(q_h)=1$
if $h=w_c$. Hence, choosing $h=w_c$ we may apply Theorem \ref{tang-conic}
to $(S,q_h)$ with $[S]=J$, which
provides us with $X$ of type I and $d_X=4$, which must be of type $C^1_I$, since we
have already excluded the type $C^1_{I(2)}$.
Similarly, any choice of $h\in I_0\sm\{0,w_c\}$ gives $X$ of type II and $d_X=4$, which yields $[X]=C^1$,
and it is left to notice that $I_0\ne\{0,w_c\}$. In fact $I_0$ is a subgroup of $I$
of index $\le2$ (as a kernel of $q_0$ that is linear on $I$) and $\rank I=6$.
\endproof

\subsection{The values of $h$ and $q$ on the ovals and bridges of $\spec$}\label{OB}
We say that a connected component $C\subset\spec_\R$ has
{\it even contact with the theta-conic $\Theta$}
if the number of tangency points of  $\Theta$ with $C$ is even
and {\it odd contact}
otherwise. Tangency points 
are counted here with multiplicities: if the local intersection index is $2k$, such a tangency point is counted with multiplicity $k$.

\proposition\label{values-of-h-on-ovals}
Consider a real non-singular cubic threefold $X$ and a real line $l\in F_\R(X)\sm\Gres$.
Let $q= q_0 +h^*$ and $h\in H_1(\spec;\Z/2)$ be the corresponding quadratic function and spectral difference class, respectively.
Then for each connected component $C\subset\spec_\R$ the following properties are equivalent:
\roster\item
$C$ has even contact with $\Theta$;
\item
the spectral covering $p_S \:\spect\to\spec$ is trivial over $C$;
\item
$h\cdot[C]=0$;
\item
$q([C])=\cases 1 &\text{ if $C$ is an oval,}\\
0 &\text{ if $C$ is the one-sided component.}\endcases$
\endroster
\endproposition

\proof
Equivalence of (1) and (2) follows from
Proposition \ref{Prym-cover-identification}, and of (2) and (3) from Proposition \ref{difference-class-covering}.
As it follows immediately from original Rokhlin's definition of $q_0$ (via the indices of normal vector fields over membranes), $q_0$ takes value $1$ on oval components.
If $S$ is of type I, this implies that $q_0$ is equal to $1$ on the one-sided component. It takes the same value on the one-sided component for any real non-singular plane quintic, since
(once more according to Rokhlin's definition) the value of $q_0$ on the one-sided component does not change under any variation
of $S$ through real plane quintics without singular points on the one-sided component,
Now, equivalence of (3) and (4) follows from that $q=q_0+h^*$.
\endproof

\proposition\label{values-of-h-on-bridges}
In notation of Proposition \ref{values-of-h-on-ovals},
for any $v=[\g]\in H_1(\spec;\Z/2)$, $cv=v$ represented by some bridge $\g\subset\spec$ the following properties are equivalent:
\roster\item
$\g$ connects real points on the components not separated by $\Theta$;
\item
the theta-covering $p_S\:\spect\to\spec$ is trivial over $\g$;
\item
$h^*(v)=h\cdot[\g]=0$.
\endroster
\endproposition

\proof
Analogous to Proposition \ref{values-of-h-on-ovals}.
\endproof

\corollary\label{Skew-matching-criterion}
Consider a real non-singular cubic threefold $X$ and a real line $l\in F_\R(X)\sm\Gres$.
Then the matching $([X],[S])$ is skew, that is $d_S=d_X+2$, if and only if the following two conditions are satisfied:
\roster\item
every oval of $\spec_\R$ has even contact
with $\Theta$;
\item
the interior of each oval of $\spec_\R$
does not contain $\Theta_\R$ and all the ovals lie in the same connected component of $P^2_\R\sm \Theta_\R$ as $J\sm \Theta_\R$.
\endroster
\endcorollary

\proof
By Corollary \ref{discrepancy-relation}(1), the matching is skew if and only if $h\in I$,
which according to Proposition \ref{symmetrization-characteristic}
is equivalent to
$h^*(K)=0$.
By Proposition \ref{K-generators}, the latter is equivalent to the vanishing of $h^*$ on the oval-classes and bridge-classes.

By Propositions \ref{values-of-h-on-ovals} and \ref{values-of-h-on-bridges},
$h^*$ vanishes on an oval-class
$[C]$
if and only if $C$ has even contact with $\Theta$,
and
it vanishes on a bridge-class $[\g]$
if and only if the the two real points of the bridge
are not separated by $\Theta$. Since any oval can be connected by a sequence of bridges with the J-component, these two equivalences
imply the statement.
\endproof

\subsection{Proof of Theorem \ref{main}(3)}\label{surjectivity of SMC}
The ``only if'' part follows immediately from Corollary \ref{discrepancy-relation} and Lemma \ref{Klein-type-relation}.
The proof of ``if'' part is split into Lemmas \ref{perfect-are-spectral} and \ref{skew-matchings-are-spectral} below.

\lemma\label{class-to-cubic}
Let $S\subset P^2$ be a non-singular real plane quintic and $h\in H_1(S;\Z/2)$ be an element of $K_0\sm\{0\}$.
\roster\item
There exist a non-singular real cubic threefold $X$ and  a real line $l\in F_\R(X)\sm\Gres$
such that
$S$ and $h$ are the spectral curve and  the spectral difference class
associated to $(X,l)$.
\item
If $h\in K_0\sm I$ then for any such pair $(X,l)$ representing $(S,h)$ the matching $([X],[S])$ is perfect.
\item
If $h\in I_0^*=(K_0\cap I)\sm\{0\}$ then for any such pair $(X,l)$ representing $(S,h)$ the matching $([X],[S])$ is skew.
\endroster
\endlemma

\proof
(1) By Proposition \ref{real-quadratic-functions}(3), the theta-characteristic represented by the quadratic function
$q=q_0+h^*$
belongs to
$\Theta^*_{\R1}(S)$. Hence, according to Proposition \ref{contact-divisors},
this theta-characteristic
 can be represented by a contact conic,
and it remains to apply Theorem \ref{tang-conic}.

(2) is an immediate consequence of Corollary \ref{discrepancy-relation}, while (3) is an immediate consequence of Corollary \ref{discrepancy-relation}
and Lemma \ref{Klein-type-relation}(1).
\endproof

\lemma\label{K-I-nonempty}
Let $S$ be a non-singular real plane quintic.
\roster\item If $S_\R$ has at least one oval, then
the set $K_0\sm I$ is non-empty.
\item If $S_\R$ has at most four ovals, then
the set $I_0^*$ is non-empty.
\endroster
\endlemma

\proof
(1) If $\spec$ has $k\ge1$ ovals, then $\rank I=d_S=6-k$ by Proposition \ref{quintic_discrepancy}
and then $\rank K=12-\rank I=6+k$
in accordance with \ref{symmetrization-characteristic}(1).
Thus, $K\sm I$ is non-empty.

Assume now first that $w_c\ne0$. Then, by Proposition \ref{symmetrization-characteristic}(3) we have
$I_1\ne\oo$.
Choosing $y\in I_1$ and $x\in K\sm I$  we deduce from $x+y\in K\sm I$ and
$q_0(x+y)=q_0(x)+q_0(y)=q_0(x)+1
$
that $q_0$ is not constant on $K\sm I$ and, hence,
$K_0\sm I\ne\oo$.

If $w_c=0$,
then $k\ge 2$, so that $\rank K/I\ge 4$.
Due to Proposition \ref{symmetrization-characteristic}(1),
it implies that
there exist $x,y\in K\sm I$ such that $x+y\in K\sm I$ and $x\cdot y=0$.
Then, for such a pair $x,y$ either $q_0(x)$, or $q_0(y)$, or
$q_0(x+y)=q_0(x)+q_0(y)$ is $0\in\Z/2$, and thus $\{x,y,x+y\}\cap (K_0\sm I)\ne\oo$.

(2) Since $k\le4$ implies $\rank  I=6-k\ge2$, the kernel $I_0$ of $q_0$ (whose restriction to $I$ is linear) has rank $\ge1$,
and so $I_0^*=I_0\sm\{0\}\ne\oo$.
\endproof

\lemma\label{3-cases-algebraically}
Let  $S$ be a real non-singular plane quintic of type $J\+1\la1\ra$, and let $J$, $O_1$ and $O_2$ denote its one-sided component, the external and the internal oval respectively.
Then:
\roster\item
$[J]\in K_0\sm I$ and $[O_1],[O_2]\in K_1\sm I$;
\item
there exists $h_1\in K_0\sm I$ such that $h_1^*(J)=1$.
\endroster
\endlemma

\demo{Proof}
Rokhlin's function $q_0$ takes value $1$ on all oval-classes. Hence, we have $q_0(O_1)=q_0(O_2)=1$, and the identity $[J]+[O_1]+[O_2]=0\in H_1(S;\Z/2)$
($S$ is of Klein type I) implies $q_0(J)=q_0(O_1)+q_0(O_2)+O_1\circ O_2=0$
so, $[O_i]\in K_1$, $i=1,2$, and $[J]\in K_0$. To complete proving (1) note that for all
bridges $b_i$, $i=1,2$, between $J$ and $O_i$
we have $b_i\circ J=b_i\circ O_i=1$, which implies that $[J],[O_i]\notin I$, $i=1,2$, since
bridge classes belong to $K$, while $I$ is orthogonal to $K$ (see Proposition \ref{symmetrization-characteristic}(1)).

Existence of a nodal degeneration of $S$ merging components $J$ and $O_1$ implies that as a
bridge $b_1$ we can take the vanishing class of such a degeneration. Then, $q_0(b_1)=1$. With such a choice, $h_1=[b_1]+[O_2]$ satisfies the required properties of (2):
$q_0(h_1)=q_0(b_1)+q_0(O_2)+b_1\circ O_2=0$
and $
h_1\circ J=1$, which implies both $h_1\notin I $ and $h_1^*(J)=1$.
\qed\enddemo

\lemma\label{3-cases}
There exist precisely 3 deformation classes of non-singular real cubic threefolds $X$ that match the class $[S]$
of non-singular real quintics $S$ of type $J\+1\la1\ra$. They are as follows:
\roster
\item
$[X]=C^3_I$, the matching is skew and type-preserving, the theta conic $\Theta$ has even contact with the ovals,
the interior of each oval of $\spec_\R$
does not contain $\Theta_\R$ and the both ovals are visible.
\item
$[X]=
C^1_{I(2)}$, the matching is perfect, the interior of each oval of $\spec_\R$
does not contain $\Theta_\R$ and the both ovals are invisible.
\item
$[X]=C^1_I$, the matching is perfect, the displacement of $\Theta_\R$ with respect to the ovals is different from those described in (1) and (2).
\endroster
\endlemma

\proof
(1) By Proposition \ref{quintic_discrepancy},
$\rank I= d_S=4$.
Hence, $\rank I_0\ge3$ which implies $I_0^*\ne\oo$.
By Lemma \ref{class-to-cubic}, any  $h\in I_0^*$
leads to a pair $(X,l)$ that realizes a skew matching such that $h$ is its spectral difference class.
Then, $[X]=C^3_I$ since $C^3_I$ is the only deformation class that may give a skew spectral matching:
it is type-preserving, and admissible mixed ones do not exist at all. The rest of claim (1) follows from Corollary \ref{Skew-matching-criterion}.

The remaining matchings, as in (2) and (3), can be obtained by the following choices of $h$: $h=[J]$ for (2), and $h=h_1$ from Lemma \ref{3-cases-algebraically}
for (3).

Indeed, let us consider the real part
$\pi_{l\R}\:X_{l\R}\to P^2_\R$ of the conic bundle $\pi_l$ defined by $l$.
Recall that a residual conic $r_s$ with $s\in P^2_\R$ has a non-empty real locus if and only if $s\in \pi_{l\R}(X_{l\R})$.
Since the real locus of any real conic is connected, connectedness of $X_\R$ is equivalent to connectedness of $\pi_{l\R}(X_\R)$.
By Proposition \ref{visible}
the boundary of $\pi_{l\R}(X_\R)$
is formed by the invisible ovals of $S_\R$. Thus, in the case $[S]=J\+1\la1\ra$
the only possibility for $W$ to be disconnected is if the both ovals are invisible, which, according to Proposition \ref{values-of-h-on-bridges},
happens for $h=h_0$, and does not happen for $h=h_1$. The rest of claims (2) and (3) follows from Corollary \ref{Skew-matching-criterion}.
\endproof

\lemma\label{perfect-are-spectral}
Any perfect matching is spectral.
\endlemma

\proof
Given a perfect matching $([X],[S])$, the quintic $S$ should have at least one oval, because quintics without ovals have Smith discrepancy $6$,
but cubic threefolds $X$ with $d_X=6$ do not exist. Hence, according to Lemma \ref{K-I-nonempty}(1) there exists $h\in K_0\sm I$, which leads to a
spectral perfect matching
$([X'],[S])$ by Proposition \ref{class-to-cubic}(2). The only case in which several perfect matching with the same $[S]$ exist
is that of
$[S]=J\+1\la1\ra $. These perfect matchings are spectral due to Proposition \ref{3-cases}.
\endproof

\lemma\label{skew-matchings-are-spectral}
Any type-preserving or admissible mixed skew matching is spectral.
\endlemma

\proof
For any skew matching $([X],[S])$ we have $d_S=d_X+2\ge2$, so $S$ should have $6-d_S\le4$ ovals.
Then by Lemma \ref {K-I-nonempty}(2) there exists $h\in I_0^*$,  which in its turn produces by Proposition \ref{class-to-cubic}(3) a
spectral skew matching $([X'],[S])$. The only cases of several skew matchings with the same $[S]$ are those of
$[S]=J$ and $[S]=J\+2$, see Figure \ref{adjacency}.

The case $[S]=J$ is treated in Corollary \ref{bottom-skew-matchings}.
Assume that $[S]= J\+2$. Then,  there exist two skew matchings: one is type-preserving and another is admissible mixed.
By Proposition \ref{symmetrization-characteristic}(5), $w_c\in I_0$. In addition, $w_c\ne0$ since
$S$ is of Klein type II. Thus, we have two options for the choice of $h$:
$h=w_c$ and $h\in I_0\sm\{0,w_c\}$. The second option does exist, since $\rank I_0 \ge \rank I-1=d_S-1=3$.
Finally, according to Proposition \ref{skew-matchings-and-wc}
the first choice gives admissible mixed and the second one gives type-preserving skew matching.
So, both matchings are spectral.
\endproof

\subsection{Proof of Theorem \ref{main} in the case $d(X)\ge4$}\label{proof-main-easy-case}
Part (1) and part (3)
are already proved in Subsections \ref{scheme-main} and \ref{surjectivity of SMC}, respectively.
Due to part (3) the number of deformation classes of $(X,l)$ with $d_X\ge 4$ is at least 6 (see Figure \ref{adjacency}). On the other hand,
due to part (1), this number is bounded from above by the sum of the numbers obtained by counting the number of
connected components in a real locus $F_\R$ of a Fano surface for each of deformation classes with Smith discrepancy $\ge 4$.
The latter is one for each of the deformation types $C^0$ and $C^1_{I(2)}$, and it is equal to two for types $C^1$ and $C^1_I$ (see
Table 1). Thus our upper and lower bounds fit, and the result follows.


\section{The conormal projection}\label{conormal-section}

\subsection{Binary codes of lines in a cubic surface}
Consider a non-singular cubic surface $Y\subset P^3$ and a line $l\subset Y$.
The pencil of planes $\{P_s\}_{s\in P^1}$ in $P^3$ passing through $l$ defines a conic bundle over $P^1$, $\pi_l\: Y\to P^1$, whose fibers are
the residual conics $r_s\subset P_s\cap Y$.
The set of points $s\in P^1$ for which the residual conic is singular consists of 5 distinct points, $\{s_1, \dots, s_5\}\subset P^1$, and for each of them the residual conic splits into a pari of distinct lines, $r_{s_i}=l^0_i\cup l^1_i$ (see Figure \ref{conic-splitting}).
This set will be called the {\it spectrum of} $(Y,l)$.

\midinsert
\hskip40mm\epsfbox{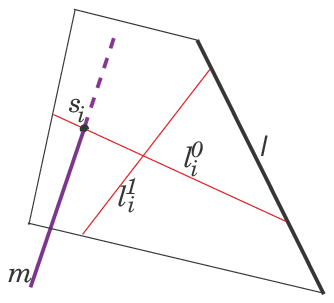}
\figure{One of the five tritangent planes passing through the  line $l$.}\label{conic-splitting}\endfigure
\endinsert

In other words, the spectrum specifies the five {\it tritangent planes} $P_s$ that contain $l$.
Note that a line $m\subset Y$ disjoint from $l$ intersects one of the two lines, $l_i^0$ or $l_i^1$, for $i=1,\dots,5$,
since $l+l_i^0+l_i^1$ form a hyperplane section divisor. Moreover, $m$ cannot intersect both lines, since a triple line intersection point
on a non-singular cubic surface  may appear only for a triple of coplanar lines.
We denote by $b_i=b_i(m)\in\{0,1\}$ the corresponding upper index, so that $l_i^{b_i}\cap m\ne\oo$ for $i=1,\dots,5$,
and thus associate to each line $m\subset Y$ disjoint from $l$ its {\it binary code}
$b_1\dots b_5$.
These codes depend on our choice of the lines $l_1^0,\dots,l_5^0$, or in other words,
on the order of points $s_1,\dots, s_5$, and the the order of lines in each pair $l_i^0, l_i^1$,
and  we call such a choice an {\it $l$-transversal marking} of $(Y,l)$.

\proposition
\label{binary-codes}
Among the 27 lines on a non-singular cubic surface $Y$ exactly 16 do not intersect a chosen line $l\subset Y$.
For any $l$-transversal marking of $(Y,l)$, the binary codes $b_1\dots b_5$
 of these 16 lines are all distinct, and the sum of their bits $b_1+\dots+b_5$ taken modulo $2$ is the same for all these 16 lines.
\endproposition

\proof The first statement is contained in the discussion that precedes the definition of binary codes. The Gram determinant of the  homology classes $[l],[ l_1^0], \dots , [l_5^0] \in H_1(Y)$ is equal to 4.
Hence, a  combination $a[l] +\sum a_i[l^0_i]$ of these 6 elements
with integer coefficients $a, a_0, \dots, a_5$ is divisible by 2 if and only if $a+\sum a_i=0\mod 2$. This implies the second statement.

\subsection{Real spectrum and truncated codes}\label{real-spectrum}
Assume now that the non-singular cubic surface $Y$ and  the line $l\subset Y$ are real.
Then its spectrum $\{s_1,\dots,s_5\}\subset P^1$ is
invariant under complex conjugation, namely, contains $c\le2$
pairs of conjugate imaginary points and $r=5-2c$ real ones.
A real spectral point $s_i$ is said to be of {\it real crossing} type if the two lines
$l_i^j$, $j=0,1$, are real, and of {\it imaginary crossing} type if they
are conjugate imaginary. Let us denote by $r_{re}$ and $r_{im}$ the number
of the real spectral points of corresponding type, $r_{re}+r_{im}=r$.

The binary codes that enumerate real lines $m\subset Y$, $m\cap l=\oo$, will be called {\it real codes}.
For a pair of conjugate imaginary  points of the spectrum, $s_j=\bar s_i$,
we may choose an $l$-transversal marking so that $l_j^0=\bar l_i^0$, then for any real line $m\subset Y$
$m\cap l_j^0\ne\oo$ if and only if $m\cap l_i^0\ne\oo$, and thus,  $b_j(m)=b_i(m)$.
So, only one of these two bits is {\it informative} and by dropping the other one from the real binary code we do not
loose any data. After dropping one bit for each pair of conjugate imaginary spectral points, we may additionally drop one of the bits
representing a real spectral point, since four bits determine the fifth one due to parity check
(see Proposition \ref{binary-codes}).
So, for real codes there are only $(r-1)+c=4-c$ informative bits and we obtain the following.

\proposition
\label{real-binary-codes}
Given a non-singular real cubic $Y$ and a real line $l\subset Y$, the set of real lines $m\subset Y$ disjoint from $l$
is empty if $r_{im}>0$ and consists of $2^{4-c}$ lines distinguished by their truncated codes if $r_{im}=0$.
\qed\endproposition

\subsection{Binary codes under nodal degeneration}\label{codes-degeneration}
Let $Y$ be a one-nodal cubic surface. Then, $Y$ contains precisely 21 lines: 6 of them,
called {\it double lines},  pass through the node, while 15 others, called {\it simple lines}, do not.
The term ``double line'' is motivated by the fact
that under any perturbation $\{Y_t\}_{t\in[0,1]}$ of $Y$ each of them splits into 2 lines on $Y_t$, $t>0$,
whereas a simple line $l\subset Y$ is varied univalently (see \cite{Seg}).

A pencil of planes $\{P_s\}$ passing through
a fixed simple line $l\subset Y$ contains precisely $4$ tritangent planes: one of them, $P_{s_0}$, passes through the node,
while 3 others, $P_{s_i}$, $i=1,2,3$, do not. Under a perturbation, the plane $P_{s_0}$, which is
called {\it double tritangent plane}, splits into a pair of tritangent planes, $P_{s_{01}}(t)$ and $P_{s_{02}}(t)$, whereas
$P_{s_i}$, $i=1,2,3$, which are called
{\it simple tritangent planes}, are varied univalently (see \cite{Seg}).

Like in the non-singular case,
we consider an $l$-transversal marking of such a pair $(Y,l)$ as a choice
of residual lines $l_i^0$ and $l_i^1$ in the planes $P_{s_i}$, $0\le i\le 3$, and
define the corresponding 4-bit codes $b_0b_1b_2b_3$ for
the 8 simple lines disjoint from $l$.
For the 4 double lines disjoint from $l$ (that is double lines different from $l_0^0$ and $l_0^1$) we use
the 3-bit codes $b_1b_2b_3$ (the bit $b_0$ is not well-defined for  these 4 lines, since they
all meet both $l_0^0$ and $l_0^1$ at the node).

Given a perturbation $\{Y_t\}_{t\in[0,1]}$ of $Y$ and an  $l$-transversal marking of $Y$, there exists a unique
perturbation $l(t)$ of $l$ and a unique $l(t)$-
transversal marking of $(Y_t, l(t)), t>0$, such that
$l_{i}^0(t)$ and $l_{i}^1(t)$ with $i=1, 2$ converge
to $l_0^0$ and $l_0^1$, respectively, while $l_{i+1}^j(t)$ with $i\ge 2$ converge to $l^j_{i-1}$.
We call such a family of markings
{\it coherent}.

\proposition\label{complex-merging} Let $\{Y_t\}_{t\in[0,1]}$ be a perturbation of a one-nodal cubic surface
$Y=Y_0$ and $l(t)\subset Y_t$ a perturbation
of a simple line $l\subset Y$. Then among the 16 lines disjoint from $l(t)$ on $ Y_t$ precisely 8 merge
pairwise forming 4 double lines in $Y$, namely, with respect to a coherent family of transversal $l$-markings, the lines encoded $01b_1b_2b_3$ and $10b_1b_2b_3$ merge together
and form a double line encoded $b_1b_2b_3$. The 8 other lines
encoded $bbb_1b_2b_3$ converge to the simple lines encoded $bb_1b_2b_3$.
\endproposition

\proof  It follows from
interpretation of the binary codes as intersection indices in $H_2(Y_t;\Z/2)$
and a simple observation that, for any pair $m'(t),m''(t)$ of lines on $Y_t$ merging
to a double line on $Y$, the difference between the classes realized by $m'(t),m''(t)$ in $H_2(Y_t;\Z/2)$
is equal (in accordance with the Picard-Lefschetz formula) to the class realized by the vanishing sphere, and the intersection of the latter class with each of the former ones is equal to 1.
\endproof

\subsection{Binary codes under real nodal degenerations}
Here, in the setting of the previous subsection,
we  assume in addition that $Y$, $l$, and $\{Y_t\}_{t\in[0,1]}$ are real.
This implies that the lines $l(t)\subset Y_t$,
the spectrum $\{s_0,s_1,s_2,s_3\}$ of $(Y,l)$, and the spectrum
$\{s_{01}(t),s_{02}(t),s_1(t),s_2(t),s_3(t)\}$ of $(Y_t,l(t))$ are also real.
Since the only node of $Y$ must be real, we have $s_0\in P^1_\R$ and
for the lines $l_0^0$ and $l_0^1$ there are 2 options: they may be either both real
or conjugate imaginary.
The points $s_{01}(t)$ and $s_{02}(t)$ can be also
either both real or conjugate imaginary.

\proposition\label{real-merging} Let  $(Y_t,l(t))$ be a real perturbation of a real one-nodal cubic surface $Y$ and let $l\subset Y$ be a real simple line.
The following criteria for merging of real lines $m(t)\subset Y_t$, $l(t)\cap m(t)\ne\oo$ in $Y_t$ as $t\to 0$ hold for any family of coherent $l(t)$-trransversal markings.
\roster\item
If the points $s_{01}(t),s_{02}(t)$
and the  lines $l_0^0$, $l_0^1$ are real, then among $2^{4-c}$ real lines $m(t)$ a half is univalent (not merging),
namely, the ones with equal bits $b_{01}=b_{02}$. The other half of real lines merge pairwise, namely, the one encoded
$01b_1b_2b_3$ merges with the one encoded $10b_1b_2b_3$.
\item
If the points $s_{01}(t)$, $s_{02}(t)$ are imaginary but the lines $l_0^0$, $l_0^1$ are real,
then all the $2^{4-c}$ real lines $m(t)$ have equal bits $b_{01}=b_{02}$ in their encodings and in particular
cannot merge.
\item
If the points $s_{01}(t)$, $s_{02}(t)$, as well as the lines $l_0^0$, $l_0^1$, are imaginary,
then all the $2^{4-c}$ real lines $m(t)$ have distinct bits $b_{01}\ne b_{02}$ and all merge pairwise, namely,
the one encoded $01b_1b_2b_3$ merges with the one encoded $10b_1b_2b_3$.
\endroster
\endproposition

\proof
Part (1) is a direct consequence of Propositions \ref{complex-merging} and \ref{real-binary-codes}. In the rest, it is sufficient to notice in addition that
$l_{01}^0(t)$  is conjugate to $l_{02}^0(t)$ in part (2)
(so that their limit $l_0^0$ is real), and conjugate to $l_{02}^1(t)$ in part (3) (so that their limit is imaginary).
\endproof

\subsection{Hyperplane sections of cubic threefolds}
Assume that $X\subset P^4$ is a non-singular cubic threefold and
$l\notin\Gres\cup\G_{II}$, so that the associated spectral quintic $S$ and theta-conic $\Theta$ are non-singular.
There is an obvious correspondence between the set of lines $m\subset P^2$
on a plane $P^2\subset P^4\sm l$
(that is used to numerate the projective planes containing $l$)
and
the set of
hyperplanes $H_m\subset P^4$ containing $l$, namely, $m=H_m\cap P^2$.
Consider the hyperplane section $Y_m=X\cap H_m$ and note that
when it is non-singular, it can be identified with
its proper inverse image $Y_{m,l}$ under the blow-up $X_l\to X$ along $l$.
 In such a case, the spectrum of $(Y_m,l)$ is the intersection $m\cap S$.

\proposition\label{nodal-cubic-fibers}
 For a non-singular cubic $X\subset P^4$ and a line  $l\notin\Gres\cup\G_{II}$, the following holds:
\roster\item
$Y_{m,l}$ is non-singular if and only if $m$
 is transverse to $S$.
\item
If $m$
is simply tangent to $S$ at a point $s\in S\sm\Theta$, that is $\ind_s(S,m)=2$,
then $Y_{m,l}$ has no other singularities than a node at some point $x\in\pi_l^{-1}(s)$.
\item
$Y_m=X\cap H_m$ is non-singular  if and only if $m$ is transverse both to $S$ and $\Theta$;
in particular, if $Y_m$ is non-singular, then so is $Y_{m,l}$.
\endroster
\endproposition

\proof
Claims (1) and (2) follow from Proposition \ref{nodality}(3).
Since $Y_{m,l}$ is the proper image of $Y_m$ under blowing up $X$ along $l$, the singularities on $Y_m$ which
are absent on $Y_{m,l}$ may appear only along $l$. On the other hand,
a singularity of $Y_m$ at a point $p\in l$ happens if and only if $H_m=T_p$. Due to Lemma \ref{theta-parametrization},
$H_m=T_p$ for $p\in l$ if and only if $m$ is tangent to $\Theta$. Wherefrom Claim (3).
\endproof

\subsection{Conormal projection}
Consider a non-singular cubic threefold $X$ and a line $\ell\subset X$ such that $\ell\notin\G_{II}$.
Let $F_l\to F$ denote the blow-up of the Fano surface $F$ of $X$ performed at
the point $l\in F$, and $E_l\subset F_l$ its
exceptional curve. Let $\spect=\spect_l\subset F$ denote, like in Section 2, the spectral covering curve formed by the
lines incident to $l$.
 The hyperplanes  in $P^4$ that contain $\ell$ form a {\it conormal projective plane} dual to $P^2$, and we denote it by $\hP$.
There is a map $\l\:F_\ell\to\hP$
that we define first\
in the complement of $\spect\cup E_\ell$.
If $m\in F\sm\spect$, and $m\ne l$,
then $m\cap l=\oo$ and we
let $\l(m)$ be the $3$-space spanned by $\ell$ and $m$.

\proposition\label{projection-well-defined}
\roster\item The above map $\l$ can be extended to
$\spect$ and $E_\ell$ and provide a well-defined regular morphism
$\l\:F_\ell\to\hP$ with finite fibers. It maps $E_\ell$ bijectively onto the conic $\hat \Theta$ dual to $\Theta$, and $\spect$ onto $\hat \Theta$
with degree 5.
\item The map $\l$ has degree $16$.\endroster
\endproposition
\proof Let $[x:y:z:u:v]$ be coordinates in $P^4$ such that the line $l$ is given by equations $x=y=z=0$.
Each line $l'$ close to $l$ intersects hyperplanes $v=0$ and $u=0$ at the points $[p_1:p_2:p_3:1:0]$ and $[q_1:q_2:q_3:0:1]$, respectively,
and we consider the 6-tuples $(p,q)$, $p=(p_1,p_2,p_3)$, $q=(q_1,q_2,q_3)$ as local coordinates
in the Grassmannian of lines $G(2,5)$.

As is known (see \cite{Clemens-Griffiths}),
the  tangent plane $T_lF$ to the Fano surface $F\subset G(2,5)$ is
defined in coordinates $(p,q)$ by
equations
$$
L_{11}(p)= 0, \, L_{22}(q)=0, \, L_{11}(q)+2L_{12}(p)= 0, \, 2L_{12}(q)+ L_{22}(p)=0 \eqtag\label{tangent-to-fano}
$$
where $L_{ij}$ are entries of
the fundamental matrix (\ref{spec-matrix}). The linear forms $L_{11}, L_{12}$, $L_{22}$ are linearly independent, as it follows
from our assumption that the line $l$ is not multiple, see Corollary \ref{double-lines}. Therefore,  the equations (\ref{tangent-to-fano}) mean that
the line in $P^2$ connecting $[p_1:p_2:p_3]$ and $[q_1:q_2:q_3]$, which is the projection of $l'$ to $P^2$ from $l$,
is tangent to the theta-conic $\Theta=\{L_{11}L_{22}-L_{12}^2=0\}$.
This proves the existence of a continuous extension of $\l$ to $E_\ell$ and implies that its sends bijectively
$E_\ell$ onto the conic dual to $\Theta$.

To extend $\lambda$ to $\spect$,  we let
 $\l(m)$ be equal to the $3$-space tangent to $X$ at the point $m\cap l$, for each $m\in \spect$.
Continuity follows from analysis of the same first order variation (\ref{tangent-to-fano}), applied this time to $m\in\spect$.
The extension of $\l$ to the whole $F_l$ is a regular map, as any continuous extension of a rational map.
Our description of $\l$ on $\spect$ implies that $\l(\spect)$ is  the conic $\hat\Theta$ dual to $\Theta$.
The restriction $\spect\to \hat\Theta$ has degree $5$ because there exist $6$ lines passing through a generic point of $X$
from which number we subtract one representing $l$.

To show finiteness of fibers of $\l$, assume that, conversely, some cubic surface
$Y=X\cap H$ traced by a hyperplane $H\subset P^4$, $l\subset H$ contains
infinitely many lines. It is well-known that such cubic $Y$
should be either a cone over irreducible cubic curve,
or contain a curve, $C$, formed by singular points of $Y$
(it can be a conic if $Y$ is reducible, or otherwise a line), see, for example, \cite{Ab}.
 In the first case, $l\subset Y$ should be one of the generators of the cone, which contradicts to that $l$ is a simple line in $X$ (which means that plane sections of $X$ cannot contain
$l$ as a double line).
 In the second case, one can easily show that $X$ must contain some singular point. Namely, we can choose coordinates $[x_0:\dots:x_4]$ in $P^4$ so that $H=\{x_4=0\}$, $X=\{f(x_0,\dots,x_4)=0\}$.
 Then $\frac{\partial f}{\partial x_i}=0$ along $C$ for $0\le i\le3$, and singular points appear on $X$ at the intersection of $C$ with the quadric $\{\frac{\partial f}{\partial x_4}=0\}$, if the latter derivative is not
 identically zero. If it is zero, then $X$ is a cone and the singular point appears at the vertex.

The degree of $\lambda$ is 16, since for any non-singular cubic surface $Y$, in particular, for $Y=X\cap H$ with a generic $H$,
this is  the number of lines in $Y$ disjoint from a given one.
\endproof

\proposition\label{folds-bifolds}
The set of critical values of $\l$ is the curve $\hspec\subset\hP$ dual to the spectral curve $S\subset P^2$.
Over each non-singular point of
$\hspec$ the projection $\l$ has $8$
unramified sheets and $4$ copies of folds
(double coverings branched along $\hspec$).
\endproposition

\proof
The first statement is a straightforward consequence of Propositions \ref{projection-well-defined} and
\ref{nodal-cubic-fibers} implying that the degree of $\l$ is $16$ over $\hP\sm\hspec$.

As it holds for any surjective proper finite holomorphic map $\lambda : X\to M$ with non-singular  covering space $X$ and non-singular target
$M$, the part $R=\lambda^{-1}(\hspec\sm \Sing \hspec)$ of the ramification locus is a non-singular curve and the ramification index of $\lambda$ is constant along each
connected component of $R$. Hence, it is sufficient to check the second statement on a dense subset of $\hspec$.

 Thus, pick a line $m\subset P^2$ simply tangent to $S$ at a point $s\notin \Theta$ and assume $m$ intersects $S$ at 3 more distinct points.
 In some pencil of lines containing $m$ take a small circle
 around $m$ and represent this circle as a path $m_\tau$, $\tau\in[0,1]$, $m_0=m_1$. Denote the points of $m_\tau\cap S$ by
 $s_\tau^i$, $\tau\in[0,1]$, $i=1,\dots,5$, where $i=1,2$ are chosen for the two points near $s$.
 Since the tangency is simple, the points $s^1_\tau$, $s^2_\tau$ alternate after a full twist, that is $s^1_1=s^2_0$ and
 $s^2_1=s^1_0$.

 Each of the surfaces $\pi_l^{-1}(m_\tau)$ is a non-singular cubic surface, and, considering the radial path from $\tau$
 to the center of the circle, we get, for each $\tau$, a perturbation $Y^\tau_t$ of the cubic surface $Y_0=\pi_l^{-1}(m)$ which is one-nodal
 by Proposition \ref{nodal-cubic-fibers}. Each of these surfaces
 contains the line $l$. Start from $\tau=0$ and choose a coherent family of $l$-transversal markings compatible with an ordering $l^0_0, l^1_0$ of the 2 lines lying on $Y_0$ over $s$, so
 that the lines $l^0_i(t), i=1,2$ that lie over $s^i_0$ for $t=1$ converge to $l^0_0$ as $t$ tends to $0$, while $l^1_i(t), i=1,2$ that lie over $s^i_0$ converge to $l^1_0$. This compatibility property implies,
 due to the unicity of coherent markings (see Proposition \ref{complex-merging}),
 that under varying $\tau$ from $0$ to $1$ and making a continuous choice of $l$-transversal markings the line $l^0_1(t)$ (respectively, $l^0_2(t)$) sends to $l^1_2(t)$ (respectively, $l^1_1(t)$).
 Therefore, the monodromy action transforms a line with initial binary code $b_1b_2b_3b_4b_5$ into the line which has the code $b_2b_1b_3b_4b_5$ with respect to the same, initial, marking.
 Thus, there are 8 invariant lines, they have binary codes
 $bbb_3b_4b_5$, and 4 transpositions of lines, $01b_3b_4b_5 \mapsto 10b_3b_4b_5$.
\endproof

Let $\l_\R\:F_{\ell\R}\to\hP_\R$ be the restriction of $\l$.
We say that a real generic line $m\subset P^2$ has
{\it real intersection type $(a,b)$}
if it has $a$ real intersection points with $\spec$ outside
$\Theta$ and $b$ real
intersection points inside (genericity means that
$m$ is not tangent to $\spec$ and $\Theta$ and has no common
points with $\spec\cap\Theta$).
Clearly, $a\ge1$ is odd, $b\ge0$ is even, and $a+b\le5$.
In fact,
$a=r_{re}$, $b=r_{im}$ are the characteristics introduced in Subsection \ref{real-spectrum}
applied to $(Y_m, l)$.

\proposition\label{folds-near-ovals} Assume that an oval $C\subset S_\R$ is convex, and
$m$ is a real line obtained by a small shift of a tangent to $C$
in the outward direction, so that the dual point $\hat m\in\hP_\R$ lies in the interior
of the dual oval $\hat C$. Let $(a,b)$ be the real intersection type of $m$ with $S$.
\roster
\item
If $b>0$, then $\l(F_{l\R})$  have no point in a neighbourhood of $\hat C$.
\item If $b=0$ and $C$ is a visible oval, then the number of unramified sheets over $\hat C$
is $2^{\frac{a+3}2}$, the number of folds is $2^{\frac{a+1}2}$, and the folds are directed towards the exterior of $\hat C$.
\item If $b=0$ and $C$ is invisible, then there are no unramified sheets, the number of folds is $2^{\frac{a+1}2}$, and the folds are directed towards the interior of $\hat C$.
\endroster
\endproposition

\proof
Proposition \ref{real-binary-codes} gives  $|\l_\R^{-1}(m)|=2^{\frac{3+a}2}$ if $b=0$ and $|\l_\R^{-1}(m)|=0$ if $b>0$.
If a point is chosen from the other side of $\hat C$, then the dual line
has the  real intersection type
$(a,b)+(2,0)$ if oval $C$ is visible
and $(a,b)+(0,2)$ otherwise. Since by Proposition \ref{folds-bifolds}, the critical locus over $\hat C$
is formed
only by folds,
the local degrees of $\l_\R$ prescribe
the number of such folds and their normal direction.
\endproof

\subsection{Special example: spectral quintic with a nest}
In what follows we investigate a pair $(S,\Theta)$ such that:
\roster\item
$S$ is a real non-singular plane quintic that has a nest of two ovals:
an oval $O_2$
in the interior of an oval $O_1$;
\item
$\Theta$ is a real non-singular conic tangent to $S$ at 5 distinct points;
\item
$O_2$ lies in the interior of $\Theta_\R$, which in its turn lies in the interior of $O_1$;
\item
at every point $s$ on $O_1$ or $O_2$,
the tangent line
to $S$ has a simple tangency
at this point
(in particular, both ovals are convex) and intersects the $J$-component of $S$ at only one point.
\endroster

As it follows from Proposition \ref{3-cases}, if the spectral pair $(S,\Theta)$ corresponds to a pair $(X,l)$ and satisfies
properties (1)--(3), then $X$ is of type $C^{1}_I$ (exotic cubic).

An example of a pair $(S,\Theta)$ that satisfies properties (1)--(4) is given by a nonsingular real conic $\Theta$ with an oval $\Theta_\R\ne\oo$
and a quintic $S$ obtained by a small real perturbation
$$
\Theta LQ\pm\varepsilon L_1\dots L_5=0, \  \ 0<\varepsilon<<1
$$
of the product $\Theta LQ $, where $Q$ is
a nonsingular  real conic whose oval $Q_\R$
contains $\Theta_\R$ in the interior, $L$ is a real line whose real locus $L_\R$ is disjoint from $Q_\R$ and $\Theta_\R$, and
$L_1\dots L_5$ are five real tangents to $\Theta$.

\proposition\label{covering-over-R}
For the conormal projection $\l_\R\:F_{l\R}\to\hP_\R$ associated to a pair $(X,l)$
whose spectral quintic and theta-conic satisfy conditions (1)--(4)
the following properties hold:
\roster
\item
$\l_\R(F_{l\R})$ is the disc bounded by $\hat O_2$; over the interior of the
annulus bounded
by $\hat O_1$ and $\hat O_2$
the map $\l_\R$ is a 8-sheeted non-ramified covering;
\item
over every point of $\hat O_2$ the number of folds is equal to 4 and all the folds are directed inward;
\item
over every point of $\hat O_1$ the number of unramified sheets is 4 and the number of folds is 2; the both folds are directed outward.
\endroster
\endproposition

\midinsert
\epsfbox{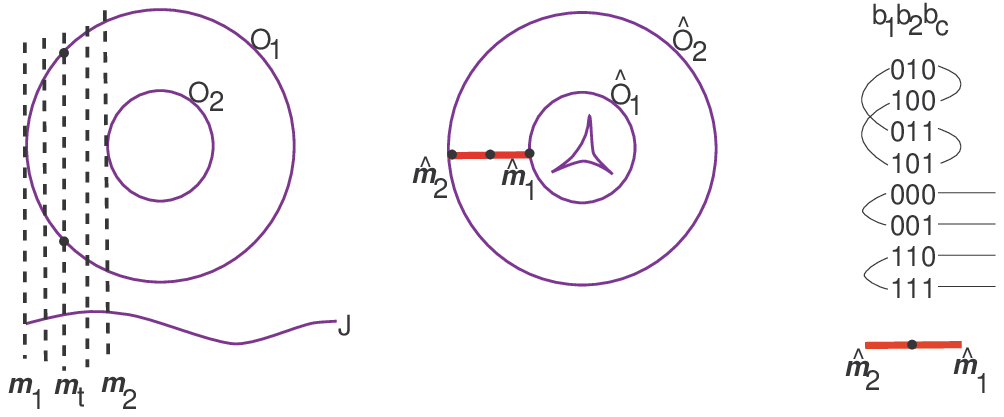}
\figure{Conormal projection for an exotic cubic threefold}\label{conormal-projection}\endfigure
\endinsert

\proof Straightforward consequence of Proposition \ref{folds-near-ovals}.
\endproof

Our aim now is to show that in the example given in Proposition \ref{covering-over-R}
the blown up Fano surface $F_\R(X)$ contains a component nicely projected to the annulus $\hat R$
bounded
by $\hat O_1$ and $\hat O_2$.

\proposition\label{Fano-component-zero-chi}
The inverse image $\l_\R^{-1}(\hat R)\subset F_{l\R}$ contains a connected component of $F_{l\R}$
whose Euler characteristic is $0$.
\endproposition

To prove this proposition we find a component of $F_{l\R}$ which is fibered over a circle.
Namely, we consider the composition $\psi\:\l_\R^{-1}(\hat R)\to S^1$, $\psi=\mu\circ\l_\R$,
of $\l_\R$ restricted over $\hat R$ with the fibration $\mu\:R\to S^1$ whose fibers
are line segments on the annulus $\hat R$ that are cut by the lines passing through a points inside $\hat O_1$.
Let us pick such a line segment $[\hat m_1,\hat m_2]$, connecting points $\hat m_i\in \hat O_i$.

\lemma\label{circle}
The inverse image $\l_\R^{-1}([\hat m_1,\hat m_2])$ is homeomorphic to a disjoint union of a circle with two closed intervals.
\endlemma

\proof
By Propositions \ref{projection-well-defined} and \ref{folds-bifolds}, $\l_\R^{-1}([\hat m_1,\hat m_2])$
consists of several copies of the line segment $[\hat m_1,\hat m_2]$ some of which are identified at the boundary points.
We will enumerate them using Proposition \ref{real-binary-codes} and apply Proposition \ref{real-merging} to precise the boundary identification.

By definition, the line segment $[\hat m_1,\hat m_2]\subset P^2_\R$
corresponds to a segment of a pencil $\{m_t\}_{t\in[1,2]}$ of lines between $m_1$ tangent to $O_1$ and $m_2$ tangent to $O_2$
(see Figure \ref{conormal-projection}).
An internal point $\hat m\in[\hat m_1,\hat m_2]$ represents a line $m=m_{t_0}$, $1<t_0<2$, crossing $O_1$ at a pair of points, $s_1,s_2$,
and not crossing $O_2$. In addition, $m$ intersects once the $J$-component of $S_\R$ at some point $s_J$
and contains a pair of conjugate imaginary points $s_c, \bar s_c\in S$.
These five points form the spectrum of the cubic $Y_{m,l}=\pi_l^{-1}(m)$.

Consider an $l$-transversal marking $l_1^0$, $l_2^0$, $l_J^0$, $l_c^0$, $l_{\bar c}^0$ corresponding to the spectral points $s_1$, $s_2$, $s_J$, $s_c$, $\bar s_c$
respectively. The points of $\l_\R^{-1}(\hat m)$ can be interpreted as real lines in $Y_{m,l}$ disjoint from $l$. Like in Subection \ref{real-spectrum},
we enumerate them via truncated codes containing only the three informative bits $b_1,b_2,b_c$ that correspond to $s_1$, $s_2$ and $s_c$.
Varying $\hat m$ along the interval $[\hat m_1,\hat m_2]$, we obtain 8 copies of $[\hat m_1,\hat m_2]$ that we denote $[1,2]_{b_1b_2b_c}$, with 8 combinations of the bits $b_1,b_2,b_c\in\{0,1\}$.

To apply Proposition \ref{real-merging}, we choose an $l$-transverse family of markings
coherently
with respect to the degeneration of $Y_{m_t,l}=\pi_l^{-1}(m_t)$ as $t\to m_1$.
Then by part (1) of that proposition, for each $b_c\in\{0,1\}$,
the intervals $[1,2]_{10b_c}$ and $[1,2]_{01b_c}$
are glued over the point $1$,  while the other 4 intervals denoted
$[1,2]_{00b_c}$ and $[1,2]_{11b_c}$ remain unglued over $1$.

Part (3) of the same proposition
implies pairwise gluing of all 8 copies $[1,2]_{b_1b_2b_c}$ over the point $2$, namely, the interval indexed with
$b_1b_20$ is glued over this point
to the interval indexed with $b_1b_21$.
Here, it is the spectral points $s_c$ and $\bar s_{c}$ that play the role of  $s_{01}$ and $s_{02}$,
but the bit corresponding to $\bar s_c$ is not informative, and is dropped.

Such a gluing of the 8 intervals $[1,2]_{b_1b_2b_c}$
produces a circle from intervals indexed with $010$, $011$, $100$, and $101$, and two line segments from the remaining 4 intervals (see Figure \ref{conormal-projection}).
\endproof

\demo{Proof of Proposition \ref{Fano-component-zero-chi} }
The monodromy of the fibration $\psi\:\l_\R^{-1}(\hat R)\to S^1$ sends the circle component (see Lemma \ref{circle}) to itself,
which yields a connected component of $F_{l\R}$ with zero Euler characteristic.
\qed\enddemo

\subsection{Proof of Theorem \ref{N-correspondence} for exotic cubics $X$}\label{N-correspondence-for exotic}
As was already seen in the previous subsection the cubic threefold $X$ treated there is of type $C_I^1$.
Hence, its Fano surface $F_\R(X)$ is a disjoint union of two components, $N(X)=\Rp2$ and $N_6$,
and $F_{l\R}$ is obtained by blowing up
$F_\R$ at its point $l$.
So, existence of a component with zero Euler characteristic in $F_{l\R}$ established in Proposition \ref{Fano-component-zero-chi}
implies that $l\in N(X)$. Together with
Theorem \ref{main} (proved for all cubics with Smith deficiency $\ge 4$ in Subsection \ref{proof-main-easy-case}), this completes the proof of Theorem \ref{N-correspondence} for exotic cubics. \qed


\section{Quadrocubics and spectral curves as their central projections}\label{node-to-quadrocubic}

\subsection{Quadrocubic associated to a singular point}\label{quadrocubic}
Given a cubic threefold
$X\subset P^4$
with a node  $\ss\in X$, the lines $l$ passing through the node
represent points $[l]$ of the associated quadrocubic
$A\subset P^3_\ss$ (see Subsection \ref{quadrocubics-edges}).
As usual, we choose
coordinates $x,y,z,u,v$ in $P^4$ so that $s$
acquires the coordinates $(0,0,0,0,1)$, present the equation of $X$ in the form $f_3+vf_2=0$, and denote
by $Q_\ss\subset P^3_s$  the quadric defined by equation $f_2=0$. Recall that $A$ is the intersection of $Q_\ss$
with the cubic surface defined by $f_3=0$.

Real nodes of real $n$-dimensional hypersurfaces are classified by the signature $(p,q), p+q=n+1$, of the Hessian matrix, or equally by the signature
of the degree two part of the real polynomial representing the hypersurface in a real local system of coordinates centered at the node.
For example, a real node $x\in A_\R$ of a real plane curve $A$ may have signature $(1,1)$, and then called a {\it cross-node},
or signature $(2,0)$ or $(0,2)$, and then called a {\it solitary node}.

In the case of a real nodal cubic threefold $X=\{f_2+f_3=0\}$
the real quadric $Q_\ss=\{f_2=0\}$ is a
hyperboloid
for signature $(2,2)$, an ellipsoid for $(3,1)$, and
a quadric without real points for $(4,0)$.
In the case of signature $(4,0)$, the  real locus of $A$ is also empty.

\proposition\label{sing-quadricubic}
Assume that $\ss$ is a node of a cubic threefold $X$. Then,
a point $\ss'\in X$, $\ss'\ne\ss$, is singular on $X$ if and only if
the line $\ss\ss'$ represents
a singular point $[\ss\ss']$ of
$A\subset P_\ss^3$, and $\ss'$ is a node of $X$ if and only if $[\ss\ss']$ is a node of $A$.
In particular, $A$ is non-singular if and only if $X$ has
no other singular point than $\ss$.

If $X$, $\ss$, $\ss'$ are real and $\ss'$ has signature $(p,q)$, then the corresponding to it node of $A$ has signature
$(p-1,q-1)$. Furthermore, if $X$ undergoes a perturbation $X_a, a\in [0,1]$ which keeps $X_a$ nodal at $\ss$ and has a local model $x_0^2+\dots x_p^2 - (y_0^2+\dots +y_q^2)= a$ in Morse coordinates at $\ss'$,
then the associated with $\ss$ quadrocubic $A$ undergoes
a transformation which has, in appropriate Morse coordinates on $Q_\ss$ at $[\ss\ss']$,
a local model $u_1^2+\dots + u_p^2 - (v_1^2+\dots +v_q^2)= g(a)$, $g(0)=0$, $g(a)>0$ for $a>0$.
\qed\endproposition

\demo{Proof {\rm ({\it cf.}, \cite{Wall})}} Choose
a real system of affine coordinates $x,y,z,t$ in a way that $\ss=(0,0,0,0)$ and $\ss'=(0,0,0,1)$,
and write an equation of $X$ in a form $f_2+f_3=0$. Then, $df_2+df_3=0$ and $f_2+f_3=0$ at $\ss'$, and using the Euler relations we get $f_2=0, f_3=0$, which shows that $[\ss\ss']\in P^3_\ss$ is
a singular point of $A$. Reciprocally, if $df_2$ and $df_3$ are linear dependent at a point $(0,0,0,1)$ and $f_2=f_3=0$ at this point, while $df_2$ does not, then $X$ is singular at the point $(0,0,0,\lambda)$ with
$df_2+\lambda df_3=0$.

Assume that $\ss'$ is a node of $X$. Then $\ss'$ is a Morse point of $f_2+f_3$, since $\ss'$ is a node. In addition, $f_2=0$ is non-singular at $\ss'$, since $\ss$ is a node. Hence
the restriction of $f_3$ on $f_2=0$ is also Morse at $\ss'$, which means that $[\ss\ss']$ is a node on $A$.

To prove the converse statement and to compare the signatures, we choose  $ f_2(\frac{x}t, \frac{y}t, \frac{z}t, 1)$ as a first coordinate $h_1$ in a local system centered at
$\ss'$ and add $1-t$ and two generic linear forms $h_2,h_3$ in $x,y,z$ as supplementary coordinates.
Then the quadratic part in Taylor expansion of the equation $f_2(x,y,z,1)+tf_3(x,y,z,1)=0$ writes in these coordinates as
$h_1(1-t+L(h_1,h_2,h_3)) + Q(h_2,h_3)$ where $L$ is a linear form and $Q$ is a quadratic form. This combined quadratic form is a direct sum of $Q$ and a
non-degenerate form of rank 2. Thus, $\ss'$ is a node on $X$ if $\ss\ss'$ is a node on $A$. In the real case,
the signature of $h_1(1-t+L(h_1,h_2,h_3)) + Q(h_2,h_3)$, which
is the signature of the node of $X$, is the sum of the signature of $Q$, which is the signature of the node on $A$,
while the signature of the complementary rank 2 summand is  $(1,1)$.

To conclude, note that according to the already proved part of the proposition making a cubic non-singular in a neighborhood of $\ss'$
is equivalent to making the quadrocubic non-singular in a neighborhood of $[\ss\ss']$. Hence, it remains only to compare the directions of the perturbations, that is to identify the
signs in appropriate local models, and it is sufficient to do it just for one arbitrary chosen perturbation of $X$. We choose the one defined by $f_2+f_3=at^3$ and repeat literally the above formulas
with Taylor expansion, which gives us $h_1(1-t+L(h_1,h_2,h_3)) + Q(h_2,h_3)=a$ as a local model for the quadrocubics, and the result follows.
 \qed\enddemo

\proposition\label{quadrocubic-projection}
Let  $S$ and $\Theta$ be the spectral curve and theta-conic of a nodal cubic threefold $X$,
and let $l\subset X$ be a line passing through a node $\ss\in X$ and not containing other nodes of $X$.
Then:
\roster\item
the central projection from $[l]$ induces a regular birational map $A\to \spec$, this map is biregular over all points
of $S$ except two singular points of $S$ which are nodes or cusps;
\item
$\Theta =2L$, where $L$ is the line passing through
these two singular points of $S$; this line is traced on $P^2_l$ by
the plane $H\subset P^3_\ss$
tangent to
$Q_\ss$ at
$[l]\in Q_\ss$.
\endroster
\endproposition

\proof
If $t\in P^2_l$ is represented by a line $l'\in A$, $l'\ne l$, then the residual conic $q_t$ is reducible, since it contains $l'$.
This shows that the image of $A$ under the central projection from $[l]$ is contained in $\spec $.
Conversely, for any given $t\in\spec$, at least one of the lines $l'\subset q_t$ passes through the singular point $\ss$, and therefore $[l']\in A$.

The central projection induces a birational map $Q_\ss\to P^2_l$ which blows up the point $[l]\in Q_\ss$ and contracts the
two generators  of $Q_\ss$ that pass through $[l]$
into two points on $P^2_l$. Since, according to Proposition \ref{sing-quadricubic},
$[l]$ is non-singular point of $A$, this implies the claim (1). This shows also that the line $L$ is traced by the plane tangent to $Q_\ss$ at $[l]\in Q_\ss$.

To check that $Q=2L$, we choose the coordinates so that $l$ has equation $x=y=z=0$ and $s=[0:0:0:0:1]$.
Then, using notation from Subsection \ref{conic-bundles-for-cubics}, we get $L_{22}=0$ and $f_2(x,y,z,u)=2(Q_2(x,y,z)+uL_{12}(x,y,z))$.
Thus shows that the tangent plane at $[l]=[0:0:0:1]$ to the quadric $f_2=0$  is defined by $L_{12}=0$.
On the other hand,  according to (\ref{theta-conic}), $\Theta$ is defined by $L_{12}^2$.
\endproof

\corollary\label{quadrocubic-and-spec}
Let $X$ be a one-nodal real cubic threefold, $l\subset X$
a real line passing through the node, and $Q, A, \spec$
the associated quadric, quadrocubic, and spectral curve.
\roster\item
If the node has signature $(3,1)$, then $S$ has a conjugate pair of imaginary nodes.
\item
If the node has signature  $(2,2)$, then $S$ has two real nodal or cuspidal points.
Such a point is a solitary node, if the corresponding line-generator of $Q$ passing through $[l]$
intersects $A$ at imaginary points, a cross-like node if the
line generator intersects $A$ at 3 distinct real points, and a cusp if
at least two of the three points coincide.
\qed\endroster
\endcorollary

\proposition\label{binodal-line-projection}
For a nodal cubic threefold $X$ and a line $l=\ss_1\ss_2$ passing through a pair of nodes $\ss_1,\ss_2\in X$
the associated spectral quintic splits as $S=S'\cup m$ where $S'$ is a quartic and $m$ is a line.
The theta-conic is the double line $\Theta=2m$.
\endproposition

\proof
Pick the coordinates in $P^4$
so that $\ss_1=[0\!:\!0\!:\!0\!:\!1\!:\!0]$ and $\ss_2=[0\!:\!0\!:\!0\!:\!0\!:\!1]$, then
in the fundamental matrix $A_s$ (see (\ref{spec-matrix}) ) of $(X,l)$
we have $L_{11}=L_{22}=0$, and thus,
$\det A_s=2L_{12}Q_1Q_2-L_{12}^2C$. So, $S$
splits into a line, $m=\{L_{12}=0\}$ and a quartic $S'=\{2Q_1Q_2-L_{12}C=0\}$,
while $\Theta=\{L_{12}^2=0\}$.
\endproof

\proposition\label{spectral-perturbation} Each real two-nodal plane quintic is a spectral quintic of a real one-nodal cubic hypersurface. Any Morse-type perturbation of the nodes
of such a quintic is realized by some
real smoothing of the corresponding cubic hypersurface.
\endproposition

\proof Given a real two-nodal plane quintic, we get a non-singular real quadro-cubic and a real point on it by blowing up the nodes and contracting after that the proper image of the line joining the nodes.
The quadrocubic with a marked point thus obtained defines a real one-nodal cubic hypersurface and a real line on it, and the initial quintic is the spectral quintic of this pair.

Let $(X,l)$ be a pair such that $X$ is a real one-nodal cubic hypersurface, $l$ is a real line on $X$ going through the node, and $S$ is the spectral, real two-nodal, quintic of this pair.
Then, under an appropriate choice of projective coordinates,  $l$ is given by equations $x=y=z=0$, $X$ by equation $xu^2+2yu+2Q_{13}u +2Q_{23}v +C=0$, and
$S$ by
$$ \det\pmatrix
x&y&Q_{13}\\
y&0&Q_{23}\\
Q_{13}&Q_{23}&C
\endpmatrix=0.
$$
The two double points of $S$ are the intersection points of the conic $Q_{23}=0$ with the line $y=0$ (cf., the proof of Proposition \ref{quadrocubic-projection} above).
They are nodal under (binary form discriminant) conditions $pC(p,0,1)-Q_{13}^2(p,0,1)\ne 0, qC(q,0,1)-Q_{13}^2(q,0,1)\ne 0$, where $(p,0,1)$ and $(q,0,1)$
are coordinates of the nodes. Finally, one can get all the Morse types of modifications by perturbations of the form
$$ \det\pmatrix
x&y&Q_{13}\\
y&t(ax+by+cz)&Q_{23}\\
Q_{13}&Q_{23}&C
\endpmatrix=0
$$
where $a,b,c$ are real constants and $t$ a small real parameter. Indeed, with such a choice the perturbation term is $t(ax+by+cz)(xC-Q_{13}^2)$, and, since the last factor is non vanishing
at the points $(p,0,1)$ and $(q,0,1)$, one can achieve, by an appropriate choice of $a$ and $c$, any given signs of the perturbation term at these points, if they are real, and achieve
non-vanishing of the pertrubation term at them, if they are imaginary (complex conjugate).
\endproof

\subsection{Isomorphisms between the quadratic cones of nodes}
Assume that $X$ is a nodal cubic threefold with nodes $\ss_0,\dots,\ss_k$, $k\ge1$.
With each node $\ss_i$ in its local projective space $P^3_i$ (formed by the lines through $\ss_i$)
we associate as usual the {\it local quadric} $Q_i\subset P_i^3$ and the quadrocubic $A_i\subset Q_i$.
These quadrics can be related naturally by means of birational isomorphisms
$\phi_{ij}\:Q_j \dashrightarrow P^2$, where $P^2$ is any plane not passing through the nodes of $X$
and $\phi_{ij}$ is induced by the central projection $P^4 \dashrightarrow P^2$ with center $l=s_is_j$.

\proposition\label{quadratic-cone-isomorphism} For each pair $0\le i,j\le k$, the following holds:
\roster
\item The spectral quintic $S$ associated with the line $\ss_i\ss_j$ splits as $S=S'\cup m$ where $S'$ is a quartic and $m$ is a line.
\item $\phi_{ij}(A_j)=\phi_{ji}(A_i)=S'$.
\item $\phi^{-1}_{ij} : P^2\to Q_i$ is the blow up of two points of $S'\cup m$ followed by the blow down of the proper image of $m$,
while $\phi^{-1}_{ji} : P^2\to Q_j$ is the blow up of two other points of $S'\cup m$ followed by the blow down of the proper image of $m$.
\endroster
\endproposition

\proof Item (1) recalls the first part of Lemma \ref{binodal-line-projection}. To get item (2) it is sufficient to follow the definitions and to notice that for each line-direction $s_jx\in A_j$
the section of $X$ by the plane generated by $s_i$ and this direction splits into 3 lines: $s_jx$, $s_is_j$, and a third line, which goes through $s_i$, since such a section should
be singular at $s_i$. Item (3) is a straightforward consequence of the classical decomposition of the central projections $Q_i\dashrightarrow P^2$ into one blow up and two blow downs.
\endproof

\subsection{6-nodal Segre cubics}
We say that $k$ points in $P^n$ are {\it in linearly general position} if no
$m+1\le n+1 $ of them generate a subspace of dimension $<m$.
If $k\ge n+1$, it is clearly
enough to require that no hyperplane contains more than $n$ points.

In what follows the nodes of a multinodal cubic threefolds $X\subset P^4$ will be always
assumed to be in linearly general position.
Such 6-nodal cubics,
studied by C.~Segre, are of a particular interest for us (for a modern overview see \cite{Dolgachev-Segre}).
These cubics $X$, which we call {\it 6-nodal Segre cubics}, can be equivalently characterized as follows.

\proposition\label{simplicity-conditions}
If a cubic $X\subset P^4$ has $6$ nodes, $\ss_0,\dots,\ss_5$ and no other singular points, the following conditions are equivalent:
\roster\item
$X$ is a Segre cubic.
\item
For each $0\le i\le 5$,
the quadrocubic $A_i$ associated with $\ss_i$ is 5-nodal,
no two of its nodes lie on the same generator of the local quadric $Q\subset P^3_{\ss_0}$
and no four are coplanar in $P^3_{\ss_i}$.
\item
Each  $A_i$ splits into a pair of irreducible transversely intersecting components
of bidegrees $(2,1)$ and $(1,2)$.
\item
The spectral curve associated to any line $l_{ij}=\ss_i\ss_j\subset X$ splits into
a pair of non-singular transverse conics
and a line transverse to the conics and not passing through  their 4 common points.
\endroster
\endproposition

\proof (1) $\Leftrightarrow$ (2) is immediate, since due to Proposition \ref{sing-quadricubic},  for each $i$, the set of nodes of $X$ distinct from $s_i$ is in bijection with the set of nodes of $A_i$,
and this bijection consists in representing the nodes of $A_i$ as line-directions $s_is_j$, $j\ne i$.

(2) $\Rightarrow$ (3)
Since the quadrocubics are 5-nodal, they are reducible. They have no neither line or conic components: indeed, in the case of line-component
the 3 points of intersection with the complementary components give 3 of the nodes on a line, and in the case of conic-component
the 4 points of intersection with the complementary components give 4 in a plane. Therefore, our quadrocubics split in two irreducible components
of bidegrees $(2,1)$ and $(1,2)$. The intersection is transversal, since the quadrocubics are nodal.

(3) $\Rightarrow$ (2) Straightforward.

(3) $\Leftrightarrow$ (4) The spectral curve $S_{ij}$ splits into $S'\cup m$, where $S'$ can be seen as the image of $A_i$ under central projection $Q_i \dashrightarrow P^2$ with center $[l_{ij}]\in Q_i$
and $m$ as the trace of the tangent plane to $Q_i$ at $[l_{ij}]$; conversely, $A_i$ can be seen as the proper image of $S_{ij}$ under blow up at 2 points of intersection of $S$ with $m$
followed by the blow down of the proper image of $m$ (see Proposition \ref{quadratic-cone-isomorphism}). Therefore, $A_i$ splits into irreducible components of bidegrees $(2,1)$ and $(1,2)$ if and only
if $S'$ splits into conics.
Transversality properties of $S$ are equivalent to nodality of $A_i$.
\endproof

\corollary\label{6-nodal-construction} Let $f_2$ and $f_3$ be homogeneous polynomials of degrees $2$ and $3$ in four unknowns.
If the quadric $f_2=0$ is non-singular and the quadro-cubic $f_2=0, f_3=0$ splits into two transversal irreducible components of bidegrees
$(2,1)$ and $(1,2)$, then $f_2+f_3=0$ is an affine equation of a 6-nodal Segre cubic $X\subset P^4$.
\qed\endcorollary

\subsection{5-point configurations on an ellipsoid}\label{5-point-config}
As is well known, the result of blowing up a real non-singular quadric at a real collection of 5 points in a linearly general position is a real del Pezzo surface of degree 3,
which embeds via the anti-canonical linear system as a real non-singular cubic surface in $P^3$. It is
B.~Segre \cite{Seg} who distributed real lines on a real cubic surface in elliptic and hyperbolics ones, and
determined which of the five exceptional curves $E\subset Y$ are elliptic and which ones are hyperbolic.
We summarize these observations of Segre as follows.

Given a set of 5 linearly generic points $\Cal P=\{p_1,\dots, p_5\}$ on an ellipsoid $Q\subset P^3$, the following bipartitions of $\Cal P$ into a pair and a triple of points
appear naturally.

\roster
\item Bipartition by the valency of $p_i$ on the convex hull of $\Cal P$: 3 vertices have valency 4 and 2 vertices have valency 3.
\item For 3 points of 5 (respectively, 2 points of 5) the central projection of $\Cal P\sm\{p_i\}$ from a point $p_i$ to any affine plane that is complementary to the tangent plane to $Q$ at $p_i$
sends the set $\Cal P\sm\{[p_i]\}$ into a convex quadrilateral (respectively, a triangle with one point inside it).
\item Among 5 exceptional real lines $E_i\subset Y$ of the blow up $Y\to Q$ at $\Cal P$, there are 3 elliptic (respectively, 2 hyperbolic) ones.
\endroster

\proposition\label{2+3partition} All the four above partitions coincide.
\endproposition

\proof The equivalence between (1) and (2) is more or less explicit in \cite{Seg}, \textsection 40), and it is straightforward.
The equivalence between (1) and (3) is implicit in \cite{Seg}, but it follows directly from the definition of ellipticity through the Segre involution.
\endproof

We call such a bipartition of a 5-point configuration on an ellipsoid the {\it Segre 3+2 bipartition}.
Applying the central projection to reduce 5-point configurations on an ellipsoid to 4-point configurations
on an affine plane, one can easily obtain the following consequence (also contained in \cite{Seg}, \textsection 40).

\proposition\label{5-point-monodromy}
The space of pairs $(Q,\Cal P)$ formed by real ellipsoids $Q$ and collections of 5 real linearly generic points $\Cal P\subset Q_\R$
is connected. Given such a pair $(Q,\Cal P)$, each even permutation $\Cal P \to\Cal P$
 that preserves Segre 3+2 bipartition
can be realized by monodromy along a continuous loop in the space of pairs $(Q,\Cal P)$.
\qed\endproposition

The similar question for less than 5 points is much simpler.

\proposition\label{k-point-monodromy}
The space of pairs $(Q,\Cal P)$ formed by real ellipsoids $Q$ and linearly generic k-point configurations $\Cal P\subset Q_\R$, $k\le4$,
is connected. Given such a pair
$(Q,\Cal P)$, each even permutation $\Cal P \to\Cal P$ if $k=4$,  and each permutation $\Cal P \to\Cal P$ if $k<4$,
can be realized by monodromy along a continuous loop in the space of pairs $(Q,\Cal P)$.
\endproposition
\proof It follows from connectedness of $SO(3)$.
\endproof

\subsection{Another classification of quadrocubics}
\label{real-quadrocubics-classification}
The object of this subsection is a {\it coarse deformation classification} of pairs $(Q,A)$ where $Q$ is a real non-singular quadric and $A\subset Q$ is a real non-singular curve of bi-degree $(3,3)$, that is,
the classification of such pairs up to deformation and projective equivalence.

Here,
by an {\it oval} of $A$ we mean a component of $A_\R$ which is null-homologous in $Q_\R$; clearly, non-oval components may
appear
only if $Q_\R$ is a hyperboloid. We say that a pair of ovals of $A_\R$ is {\it separated}
if they bound disjoint discs on $Q_\R$,  and that a set of 3 ovals $O_1$, $O_2$ and $O_3$ forms
a {\it 3-nest} on $Q_\R$ if a disc $D_1$ bounded by $O_1$ contains a disc $D_2$ bounded by $O_2$,
and $D_2$ contains $O_3$.
A non-oval component of $A_\R$ on a hyperboloid is said to be of type $(p,q)$, $0\le q\le p$, if, under an appropriate orientation of this component and line-generators of $Q_\R$, it realizes the
class $(p,q)\in\Z^2\cong H_1(Q_\R)$
with respect to the basis in $H_1(Q_\R)$ given by the line-generators.

The following theorem is a straightforward consequence of deformation classification of real bi-degree $(3,3)$ curves on real quadrics, as established
in \cite{Degt-Zvon}, \cite{Zvon}.

\theorem\label{DZ}
\roster
\item Pairs $(Q,A)$ with $Q_\R=\emptyset$ are all coarse deformation equivalent to each other.
\item
Pairs $(Q,A)$ where $Q_\R$ is an ellipsoid form 8 coarse deformation classes: one of type {\bf 1}${}_1$, six of type  {\bf k}, $0\le k\le5$, with $A_\R$ consisting of  $k$ disjoint ovals, and one of type {\bf 3}${}_I$ with
 $A_\R$ consisting of 3 ovals forming a nest.
\item
Pairs $(Q,A)$ where $Q_\R$ is a hyperboloid form 8 coarse deformation classes: five of type {\bf k}, $1\le k\le5$ with $A_\R$ consisting of $k-1$ disjoint ovals
and component of type $(1,1)$, one of type {\bf 1}${}_I$  with $A_\R$ consisting of a sole $(3,1)$-component, and two of type {\bf 3}${}_I$
with $A_\R$ consisting either of three homologous to each other
type $(1,1)$ components, or of two separated ovals and a type $(1,1)$ component.
\qed\endroster
\endtheorem


\subsection{Ascending and descending perturbations of a node on $X_\R$}
Each real one-nodal cubic threefold $X_0$ admits two types of real perturbations $X_t$, $t\in[0,1]$:
as it is shown on Figure \ref{adjacency}, for one type of perturbations the Smith discrepancy of $X_t$, $t>0$, is equal to the Smith discrepancy
of the quadrocubic associated with $X_0$, for the other type it is less by 1 (see, for example, \cite{Kr-Fano}). We call {\it ascending} the perturbations of the first kind,
and {\it descending} for the other. In terms of Figure \ref{adjacency}, an ascending perturbation represents an ascending oriented edge of $\G_{3,3}$ (direction
of growing $d$) and a descending perturbation represent a descending oriented edge.

Perturbations of real one-nodal plane quintics have similar properties, and we apply to them a similar terminology: an ascending perturbation represents an ascending oriented
edge of $\G_{5,1}$, and a descending perturbation a descending edge.

These definitions are consistent with the spectral correspondence:  if we pick a real line $l\subset X_0$ not through the node and a family $l_t\subset X_t, t\ge 0$, of real lines accompanying a perturbation as above, then
the spectral curve $S_0$ experiences an ascending (respectively, descending) perturbation, if the perturbation of $X_0$ is ascending (respectively, descending). Such a relation follows, for example, from
the congruence $d(S_t)=d(X_t)\mod2$.

\proposition\label{S-smoothings}
If $X_0$ is a real one-nodal cubic threefold whose node has signature $(3,1)$ and $l_0\subset X_0$ is a real line passing through the node,
then each real perturbation $(X_t,l_t)$, $t\in[0,1]$, of $(X_0,l_0)$ and the associated with it perturbation of $S_0$ are ascending,
and the matchings $([X_t],[S_t]), t >0$, are skew.
\endproposition

\proof
According to Corollary \ref{quadrocubic-and-spec}(1), the nodes of $S_0$ are imaginary, and so the real locus $A_\R$ of the quadrocubic $A$
associated with $X_0$ is projected homeomorphically on $S_{0\R}$. Hence, $S_{0\R}$ has the same number $r$ of connected components as $A_\R$,
that is $r= 5-d_A$. Since the nodes of $S_0$ are imaginary, $S_{t\R}$ are homeomorphic
to $S_{0\R}$, and hence have the discrepancy $7-r=d_A+2$.

On the other hand, by Theorem \ref{quadrocubics_as_edges} and Corollary \ref{Smith-Klein-preserving}, ascending perturbations of $X_0$
have the discrepancy $d_X=d_A$, while the discrepancy for descending perturbations is $d_A+1$ (existence of a real line $l_0$ guarantees
that $A_\R\ne\oo$ and so, Corollary \ref{Smith-Klein-preserving} is applicable).
So, since $d(S_t)$ and $d(X_t)$ must have the same parity
(see Corollary \ref{discrepancy-relation}),
the considered perturbation is ascending, and the matching $([X_t],[S_t])$ is skew due to $d_{S_t}=d_{X_t}+2$.
\endproof


\subsection{Monodromy permutation of ovals of quadrocubics}
As is well known (see, for example, \cite{Zvon}), any real quadrocubic with 5 ovals is of Klein type I, and if it lies on an ellipsoid $Q$ then, with respect to the complex semiorientation, 3 ovals
have the same sign, while 2 other ovals have the opposite sign. We call this decomposition the {\it complex orientation bipartition}.

Denote by $V^{3,1}_k$, $0\le k\le5$, the space of real one nodal cubic threefolds with a node of signature $(3,1)$ whose quadrocubic, $A_\R$,
consists of $k$ real components bounding disjoint discs on the ellipsoid $Q_\R$ (the latter assumption excludes only quadrocubics of Klein's type I with 3 ovals
on this ellipsoid).

\proposition\label{monodromy-oval-transitivity}
\roster\item
The space $V^{3,1}_5$ is connected.
 The complex orientation bipartition of the ovals of the quadrocubics $A_\R$ is preserved along any continuous paths in $V^{3,1}_5$, and
any permutation preserving this bipartition can be realized by a continuous loop.
\item
For $k\le4$, the monodromy in $V^{3,1}_k$
acts transitively on the set of ovals of $A_\R$, and on the set of unordered pairs of  ovals.
\endroster
\endproposition

\proof
Part (1): According to Theorem \ref{DZ}, the pairs $(Q_\R, A_\R)$ where $Q_\R$ is an ellipsoid and $A_\R$ is a real quadrocubic with 5 ovals are all deformation equivalent. Combining this result
with Proposition \ref{cubics-quadrocubics-association} (and connectedness of the group of real projective transformations of $P^4_\R$)
we obtain the connectedness claim. The invariance of this bipartition is trivial.

To prove the claim on permutations, let us consider 5 points on $Q_\R$ in linearly general position. Then, there exists a unique curve $B_1$ of bidegree $(2,1)$ and a unique curve $B_2$ of bidegree $(1,2)$ passing through these points.
These curves are complex conjugate to each other and intersect each other transversally. Their union is a curve of bidegree $(3,3)$, and we denote by $f_C=0$ an equation of $C=B_1\cup B_2$ on $Q$, and an equation of $Q$ by $f_2=0$.
Consider then a small perturbation $C'$ of $B_1\cup B_2$ producing a real curve with 5 ovals, say a perturbation given by an equation $f'_{C}=0, f'_{C}= f_C+ \epsilon g_3$ where $g_3$ is any homogeneous polynomial of degree 3 without any
zero on $Q_\R$.

The bipartition of the 5 ovals according to complex semiorientations coincides with the bipartition transpoted to the ovals from the Segre bipartition of the points $p_1,\dots,p_5$. To justify it, it is sufficient:
to consider another real quadrocubic with 5 ovals, namely, a 5-ovals quadrocubic obtained by a perturbation of the union of 3 circles cut by 3 sides of a regular tetrahedron on a sphere concentric to the circumscribed
sphere, but of a bit bigger radius;  to notice that the limit of a complex orientation of this quadrocubic should give an orientation of the circles;  to observe that thus each plane intersecting the 3 ovals of the same sign separate the 2 other ovals;
and apply finally the first definition of Segre bipartition, see Subsection \ref{5-point-config}.

Now,  given any even permutation of 5 ovals of $C'$ preserving their bipartition, we apply Proposition \ref{5-point-monodromy} and realize this permutation by a family of quadrocubics given, as above,
by equations $f'_{C_t}=0, f'_{C_t}=f_{C_t}+\epsilon g_3$.
The passage from a loop of quadrocubics to a loop of cubic threefolds is ensured, for example, by Lemma \ref{nodes-to-quadrocubics}. To get the odd permutations preserving the bipartition it is sufficient to combine the above loops with a loop
in the space of real one-nodal cubic threefolds that moves the node in a way reversing the local orientations of $P^4_\R$.

Part (2) of the Proposition is proved similarly, modulo replacement of the reference to Proposition \ref{5-point-monodromy} by the reference to  Proposition \ref{k-point-monodromy}.
\endproof

\subsection{Atoric nodal cubics and monodromy permutations of their nodes}
We say that a real nodal
cubic threefold $X\subset P^4$ is {\it atoric} if:
\roster\item
its nodes are real and  in linearly general position;
\item all the nodes have signature $(3,1)$;
\item the real loci of quadrocubics associated with the nodes consist only of solitary points.
\endroster
 The name is motivated by Proposition \ref{multi-nodal-perturbation}, where we show that the real Fano surface $F_\R(X)$
 of such a cubic $X$
 has no torus components.

It turns out that if the conditions (2) - (3) hold for one node of $X$, then they are satisfied for the others.

\lemma\label{nodes_coherence}
Assume that $X$ is a real $k$-nodal cubic satisfying the above condition (1). Then if one of the nodes, $\ss\in X$
has signature $(3,1)$ and its quadrocubic $A_\R$ consists of $k-1$ solitary nodes, then
$X$ is atoric (that is, the conditions (2) and (3) hold for all the other nodes as well).
\endlemma

\proof
According to Proposition \ref{quadratic-cone-isomorphism},  the quadric $Q_{\ss'}$ associated with the node $\ss'\ne \ss$ is obtained from
the ellipsoidal quadric $Q_\ss$ associated with $\ss$ by the following sequence of birational transformations:
the blowing up of $Q_\ss$ at the point $p'$ represented by the line $\ss\ss'$, the blowing down of the pair of conjugate imaginary exceptional curves
arisen from linear generators, the blowing up of the pair of conjugate imaginary points on the image $T$ of the section of $Q_\ss$ by the plane tangent to
$Q_\ss$ at $p'$, and, finally, the blow down of the  real exceptional curve arised from $T$.
The result is an ellipsoid with a quadrocubic consisting of $k-1$ solitary nodes on it, since the real locus $(Q_{\ss})_\R$ experienced only one blowup and one blowdown.
\endproof

\proposition\label{monodromy-nodes-with-a fixed-one}
Assume that $X\subset P^4$ is an atoric real $k$-nodal cubic, $3\le k\le 5$, and $\ss\in X$
is a fixed node.
Then the real equisingular deformations of $X$ preserving the node $\ss$ fixed act transitively on the other nodes.
\endproposition

\proof Straightforward consequence of Proposition \ref{monodromy-oval-transitivity} and Lemma \ref{nodes-to-quadrocubics} .
\endproof

An immediate consequence of Proposition \ref{monodromy-nodes-with-a fixed-one} is the following.

\corollary\label{monodromy-of-nodes}
For an atoric real  $k$-nodal cubic $X\subset P^4$ with $3\le k\le 5$,
the monodromy action of real equisingular deformations is transitive on the set of nodes and the set of unordered pairs of nodes.
\qed\endcorollary

\subsection{Perturbation of atoric nodal cubics}\label{perturb-atoric} The main strata of the real locus of the discriminant hypersurface $\Delta^{3,3}\subset \Cal C^{3,3}$
that correspond to real one-nodal cubic threefolds with signature $(p,q)$ different from $(2,2)$ are coorientable. In the case of signature $(3,1)$,
real perturbations modeled in local Morse coordinates by $x_1^2+x_2^2+x_3^2-x_4^2=a$ with $a>0$, which we call {\it one-sheeted perturbations}, lead to one side,
and those with $a<0$, which we call  {\it two-sheeted}, to the opposite one.

Now, consider a real k-nodal cubic threefold $X_0$ whose nodes $\ss_1,\dots,\ss_k$ have all signature $(3,1)$.
If its real perturbation $\{X_t\}_{t\in[0,1]}$ is one-sheeted at $i$ nodes and two-sheeted at the other $k-i$ nodes, we call it
a {\it perturbation of real type $(i,k-i)$}.
Two perturbations of $X_0$ that are of the same type at each of the $k$ nodes are called {\it coherent}.

To extend the notion of ascending perturbation from $k=1$ to any $k$, we define a real perturbation of $X_0$ {\it ascending}, if it produces the minimal possible
value of Smith deficiency.
We say that two real perturbations $X_t$ and $X'_t$ of $X_0$ are {\it deformation equivalent}, if they can be connected by a
continuous family of perturbations.

\proposition\label{ascending-coherent}
Let $X_0$ be an atoric real nodal cubic threefold with $k\le 6$ nodes. Then:
\roster
\item
Any two coherent real perturbations of $X_0$ are locally deformation equivalent.
\item
Each ascending perturbation of $X_0$ is one-sheeted at every node.
In particular, the ascending perturbations of $X_0$ are all coherent and, thus, deformation equivalent.
\endroster
\endproposition

\proof For any $k\le 6$, each of the strata of nodal cubic threefolds with $k$ nodes is an intersection of $k$ transversal  branches of the discriminant.
This implies immediately the first claim.

For $k=1$, the second claim follows from the first claim and Proposition \ref{S-smoothings}.
If $k>1$, we perform the perturbation in two steps. First, we keep one node and perturb the other nodes. As it follows from Proposition \ref{sing-quadricubic},
to provide a quadrocubic with the maximal possible number of ovals, at the latter nodes such a perturbation should be one-sheeted. Note also that such a perturbation
always exists. Thus, the case $k>1$ is reduced to the case $k=1$.
\endproof


\section{Fano surfaces of nodal cubics}

\subsection{Fano surface of a nodal cubic threefold  \cite{Clemens-Griffiths}, \cite{Hadan}}\label{Fano_nodal}
Assume that $X$ is a nodal cubic threefold containing no planes, $\ss\in X$ is one of the nodes, and $A$ is the quadrocubic associated with
this node.
Then there appear a natural birational
map $\s\:\Sym^2(A)\dashrightarrow F(X)$ from the symmetric square of $A$ to the Fano surface.
This map sends a generic pair $\{[l_1],[l_2]\}\in\Sym^2(A)$ to the line $l\in F(X)$ such that $l_1+l_2+l$ is a plane section  of $X$.

There are also two natural embeddings $\phi_i\:A\to\Sym^2(A)$, $i=1,2$, defined as follows.
For every $a\in A$,
each of the two generators $g_i(a)$, $i=1,2$, of the quadric $f_2=0$ passing through $a$
intersects $A$ at two more points in addition to $a$,
and the embeddings $\phi_i$ send $a\in A$ to the pair of such additional intersection points of $A$ with $g_i(a)$.

\proposition\label{Fano_for_nodal_cubic}
If $X$ is a one-nodal cubic threefold, then $F(X)$ is a singular surface, its singular locus $\Sing F(X)$ is a double curve with a simple double crossing,
$\s\:\Sym^2(A)\to F(X)$ is the normalisation map, and, for each $i=1, 2$, the composition $\s \circ \phi_i$ maps $A$ on $\Sing F(X)$ isomorphically.
\qed\endproposition

Now, consider a multi-nodal cubic threefold $X$ with nodes $\ss_i$, $i=0,\dots,k$. Then, according to Proposition \ref{sing-quadricubic},
the quadrocubic $A$ associated to $\ss_0$ is nodal and its nodes
$p_i$, $i=1,\dots, k$ are represented by the lines $\ss_0\ss_i$.
Note that the embeddings  $\phi_i\:A\to\Sym^2(A), i=1,2,$
lift to embeddings $\hat\phi_i\:\hat A\to\Sym^2(\hat A)$, where
$\hat A$ is the normalization of $A$. Note also that each nodal point $p_i$ of $A$ lifts to a pair of points $p'_i, p''_i$ in $\hat A$.

\proposition\label{sym-square-as-Fano}
For any multi-nodal cubic threefold $X$, the
normalization map $\hat F(X)\to F(X)$ of the Fano surface $F(X)$ decomposes into the
composition of the map $\hat F(X)\to \Sym^2(\hat A)$
obtained by
blowing up the points $(p_i',p_i'')$, $i=1,\dots k$, the projection $\Sym^2(\hat A) \to \Sym^2(A)$ induced by the normalization of $A$,
and the mapping $\s\:\Sym^2(A)\to F(X)$.
\qed\endproposition

\subsection{Fano surface of a 6-nodal Segre cubic threefold \cite{Dolgachev-Segre}, \cite{Hasset-T}}\label{6nodal-complex-cubic-section}
We need to analyze in more details the normalization map from Proposition \ref{sym-square-as-Fano}
when $X$ is a
6-nodal Segre cubic threefold $X$,
in which case the quadrocubic  $A$
associated to a node $\ss_0$
splits into two non-singular rational irreducible components $B_1$, $B_2$
(see Proposition \ref{simplicity-conditions}) that intersect each other
at 5 points $p_i$, $i=1,\dots,5$ ,
given by the lines $\ss_0\ss_i$.
Thus,
 $$\Sym^2(\hat A)=\Sym^2(B_1)\+\Sym^2(B_2)\+(B_1\times B_2),$$
and Proposition \ref{sym-square-as-Fano} implies that  $$\hat F(X)=P_1^2\+P_2^2\+R,$$
where $P^2_i=\Sym^2(B_i)$, $i=1,2$, are projective planes and $R$ is a del Pezzo surface of degree $3$
obtained from $B_1\times B_2$ by blowing up
the 5 points
$(p_i^1, p_i^2)$ representing the intersection points $p_i\in B_1\cap B_2$, $i=1,\dots,5$.
Using an anticanonical embedding, $R$ can be viewed as a non-singular
cubic surface in $P^3$
(well-defined up to projective equivalence). In what follows we call $R$ (respectively, its image $R'$ in $F(X)$) the {\it central component} of
F(X) (respectively, $\hat F(X)$).

Each of the planes $P^2_j=\Sym^2(B_j)$, $j=1,2$, contains $5$ lines that are the images of $B_j\times p_i^j\subset B_j^2$
under the projection $B_j^2\to P^2_j$.
According to Proposition \ref{sym-square-as-Fano}, to obtain $F(X)$
these lines should be identified
with the corresponding five disjoint lines $L_{j,i}\subset R$ that are proper inverse images of
$B_1\times p_i^2$  (respectively, $p_i^1\times B_2$) if $ j=1$ (respectively, $j=2$).

One can easily check also that the curve $\phi_1(B_2)$
(respectively, $\phi_2(B_1)$) is also a line in $P^2_1$ (respectively, $P^2_2$)
and  that in $F(X)$ the six described lines in $P^2_1$ (respectively, $P^2_2$)
are  identified with the six lines
$L_{1,i}$ (respectively, $L_{2,i}$), $i=0,\dots,5$,  in $R$,  where
$L_{1,0}$ (respectively, $L_{2,0}$) is the proper inverse image in $R$
of $\phi_2(B_2)$ (respectively, $\phi_1(B_1)$).

Note that the two sextuples of lines $L_{1,i}$ and $L_{2,i}$ form a Schl\"{a}fli double six on $R$.
Note also that identification of the six disjoint lines $L_{1,i}\subset R$ with pairwise intersecting six lines on $P_1^2$
and that of lines $L_{2,i}\subset R$ with lines on $P_2^2$
results  on $R$ in
gluing of 15 pairs of points $L_{1,i}\cap L_{2,j}$ and $L_{1,j}\cap L_{2,i}$,
$0\le i<j\le5$, which
gives 15 singular (not normal) points
$q_{ij}$
on the result $R'$ of this factorization.
The points $q_{ij}$ are {\it quadruple} in the sense that
they have four local branches: two on $R'$ and one on each of the planes $P_1^2$, $P_2^2$
(where
$q_{ij}$ are identified with the 15 pairwise intersection points of the 6 lines on each of $P^2_j$).

The above observations belong to C.~Segre \cite{Corrado}
({\it cf.,} \cite{Dolgachev-Segre}) and can be summarized as follows.

\proposition\label{quadrocubics-on-quadrics}
The Fano surface $F(X)$ of a 6-nodal Segre cubic threefold $X$ splits into 3 irreducible components,
$F(X)=P_1^2\cup P_2^2\cup R'$, two of which are  projective planes $P_i^2=\Sym^2(B_i)$.
The third component $R'$ is
the quotient of the cubic surface $R$
obtained from $B_1\times B_2$
by blowing it up at 5 points,
$(p_i^1, p_i^2), i=1,\dots, 5$, representing 5 intersection points $p_i\in B_1\cap B_2$.
The factorization $R\to R'$
consists in gluing pairwise 15 pairs of intersection points  $L_{1,i}\cap L_{2,j}$ and $L_{1,j}\cap L_{2,i}$, $0\le i<j\le5$ and is induced by
identification of $L_{1,0},\,L_{1,1},\,L_{1,2},\,L_{1,3},\,L_{1,4},\,L_{1,5}$ with
6 generic lines in $P_1^2$ and $L_{2,0},\,L_{2,1},\,L_{2,2},\,L_{2,3},\,L_{2,4},\,L_{2,5}$ with 6  generic lines in $P_2^2$.
\qed\endproposition

In the above construction it is often useful to identify
the product $B_1\times B_2$ with the local quadric $Q_{0}$ associated to the chosen node $\ss_0$.

\proposition\label{cubic-surfaces-from quadrics}
Let $Q_0$ and $A=B_1\cup B_2$ be the quadric and
the quadrocubic associated to one of the nodes $\ss_0$ of a 6-nodal Segre cubic threefold $X$.
Then, $B_1\times B_2$ can be canonically identified with
$Q_0$, which yields a canonical isomorphism  between
the central component $R$ of $\hat F(X)$  and $Q_0$ blown up at the 5 nodes of $A$.
\endproposition

\proof
By Proposition \ref{simplicity-conditions}, one of the components of $A$, say $B_1$, has bidegree $(2,1)$, while the other component, $B_2$, has bidegree $(1,2)$.
Projections of $B_1$ to the first and $B_2$ to the second factor (line-generator) of $Q_0$
are isomorphisms, which identifies
$Q_0$ with $B_1\times B_2$ and induces
the required isomorphism.
\endproof

Note that although the Fano surface $F(X)$ is independent of the choice of a node $\ss_0$, the above construction makes use of such a choice, and
this choice yields a distinguished pair of lines,
$\binom{L_{1,0}}{L_{2.0}}$
(represented by the curves $B_1$ and $B_2$), in the double six
$\binom{ L_{1,0}\,L_{1,1}\,L_{1,2}\,L_{1,3}\,L_{1,4}\,L_{1,5}}{L_{2,0}\,L_{2,1}\,L_{2,2}\,L_{2,3}\,L_{2,4}\,L_{2,5}}$.
It follows from Proposition \ref {quadratic-cone-isomorphism}
that,
for any $i=1,\dots,5$, the pair $\binom{L_{1,i}}{L_{2.i}}$ represents the irreducible components of the
quadrocubic $A_i\subset Q_{i}$ associated with $\ss_i$. Thus, the next statement is a consequence of the two previous propositions.

\corollary\label{nodes-as-lines}
In notation and under assumptions of Proposition \ref{quadrocubics-on-quadrics}:
\roster
\item
The six nodes $\ss_i$ of $X$ are in a canonical one-to-one correspondence with the six pairs of lines,
$\binom{L_{1,i}}{L_{2,i}}$,
of the Schl\"{a}fli double six $\binom{ L_{1,0}\,L_{1,1}\,L_{1,2}\,L_{1,3}\,L_{1,4}\,L_{1,5}}{L_{2,0}\,L_{2,1}\,L_{2,2}\,L_{2,3}\,L_{2,4}\,L_{2,5}}$.
\item
Each of the local quadrics $Q_i$ can be canonically identified with the result of contraction of the 5 lines
$L_{i,j}$, $0\le j\le5$, $j\ne i$, on the
cubic $R$ that are incident with the lines $L_{1,i}$, $L_{1,j}$, $L_{2,i}$, $L_{2,j}$ and not incident with the other
lines of the above double six.
\item
Each of the points $[\ss_i\ss_j]\in F(X)$ represented by the line
$\ss_i\ss_j\subset X$ connecting the nodes $s_i, s_j$ is a quadruple point of $F(X)$ and
belongs to the image of
$L_{1,i},L_{2,i},
L_{1,j},L_{2,j}\subset R$ in $R'$. \qed
\endroster
\endcorollary

\subsection{Real Fano surfaces of one-nodal cubics \cite{Kr-Fano}}\label{real-Fano-subsection}
Assume that $X$ is a real one-nodal cubic threefold, and that each of $A$, $F(X)$, and $\hat F(X)$ is equipped with the induced real structure.
Then, the real locus of $\Sym^2(A)=\hat F(X)$ consists of  $\binom{r}2+1$ connected components.
One component, which we denote $N(X)$, is non-orientable; it is obtained from the quotient $A/\conj$ of $A$ by the complex conjugation
by filling its holes with the $r$ M\"obius bands $\Sym^2(C_i)$, $i=1,\dots, r$, where $C_1,\dots C_r$ denotes the
 connected components of $A_\R$. The other $\binom{r}2$ components are tori $T_{ij}=C_i\times C_j$, $1\le i<j\le r$.

\proposition\label{Collapsible and T-mergeable tori}
Let $X$ be
a real one-nodal cubic threefold.
If the node if $X$ is of signature $(p,q)$ with $p,q\ge1$,
then the signature of the node of
the normal real slice of $F(X)$ at a point of $A_\R$
(such a slice is nodal according to Proposition \ref{Fano_for_nodal_cubic})
is $(p-1,q-1)$.
If the node
of $X$
has signature $(4,0)$ then $A_\R=\oo$.
\qed\endproposition

Recall that for real surfaces with a simple double crossing
a normal real
slice of signature $(2,0)$ (respectively, $(1,1)$)
can be characterized by a local model $x^2+y^2=0$ (respectively, $x^2-y^2=0$)
in the 3-space with local coordinates $(x,y,z)$.
Let us call a component $C$ of $A_\R$ {\it solitary circle}
if the normal real slices along this component have
signature $(2,0)$ and {\it cross-like intersection circle} in the case of signature $(1,1)$.

It is well known also that
a perturbation of $X$ leads to a smoothing of the singular locus $A$ of $F(X)$ modeled by a usual smoothing of
normal slices (cf., Krasnov's \cite{Kr-Fano} where it was crucial for the proof of the results we reproduced in Subsection \ref{real_Fano}).
 In particular, if the node of $X$ has signature $(3,1)$, then
the real singular locus $A_\R\subset F(X)$ consists of $0\le r\le5$ solitary circles, and after ascending
real perturbation $X_t$, $t\in[0,1]$, of $X_0=X$ such circles $C_i\subset A_\R$, $1\le i\le r$, are perturbed simultaneously into $r$ toric components
$T_{0i}$ of
$F_\R(X_t)$, $t>0$, while in the case of descending perturbation, these circles simultaneously vanish.

Those $r$ toric component $T_{0i}$ in the case of ascending perturbation $X_t$
will be called {\it collapsing tori} of $F_\R(X_t)$, $t>0$, relative to this perturbation (that can be viewed also as a nodal degeneration
as $t\to0$).

\proposition
Assume that a real non-singular cubic threefold $X_1$ is obtained by
 an ascending perturbation $X_t$, $t\in[0,1]$, of
 a one-nodal cubic $X_0$ whose node has signature $(3,1)$.
 Then if $X_1$ has Klein type II,
then the monodromy action is
transitive on the
set of
collapsing tori of $F_\R(X_1)$ as well as on the set
of the
non collapsing ones.
In particular, there are at most two orbits
in the set of toric components of $F_\R(X_1)$.
\endproposition

\proof
The collapsible tori arise as result of perturbation of solitary circles represented by the components $C_1,\dots,C_r$ of $A_\R$ treated as the double curve of $F_\R(X)$
and non collapsing ones by the products $C_i\times C_j$, $1\le i<j\le r$.
Since the Klein type of $X_1$ is the same as that of $A$ (see Theorem \ref{quadrocubics_as_edges}),  we can apply
Proposition \ref{monodromy-oval-transitivity} and conclude, using Lemma \ref{nodes-to-quadrocubics}  together with Proposition \ref{ascending-coherent} and continuity argument, that the monodromy permutes the components $C_i$ as well as their pairs,
which implies the required claim.
\endproof

A somewhat different perturbation scheme appears if $X$ has a node of signature $(2,2)$.
Then Propositions \ref{sym-square-as-Fano} and \ref{Collapsible and T-mergeable tori} imply that
the real components of the associated quadrocubic $A$ look as simple normal
self-intersection curves of $F_\R(X)$.
More precisely, each singular connected component $\Sigma $ of $F_\R(X)$ contains
one component of quadrocubic, $C\subset A_\R$, as such self-intersection curve.
In particular, some singular components $\Sigma $ may be homeomorphic to the direct product
of a figure-eight (wedge of two circles) with a circle.
 The same arguments (due to Krasnov \cite{Kr-Fano}) as in the previous case imply that an ascending perturbation of $X$ leads to a perturbation of such component $\Sigma$
into a pair of tori in $F_\R(X_t)$, $t>0$,
(the figure-eight factor is smoothes into a pair of circles). The tori that appear in $F_\R(X_t)$, $t>0$, in such a way
will be call {\it T-mergeable} with respect to the nodal degeneration $X_t$, $t\to 0$.
If perturbation $X_t$ is descending, then
the figure-eight factor of $\Sigma $ is smoothed into one circle, so that $\Sigma $ is perturbed into one torus component on $F_\R(X_t)$.

\proposition\label{collapsible-mergeable}
Consider an ascending perturbation $X_t$, $t\in[0,1]$, of a one-nodal real cubic $X_0$ that yields a non-singular real cubic threefold $X=X_1$.
Assume that $l\subset X$ is a real line representing a point $[l]\in F_\R(X)$ on a torus which is
collapsible, if the node of $X_0$ has signature $(3,1)$, and  T-mergeable in the case of signature $(2,2)$.
Then $l$ realises a skew spectral matching, or in other words, for the spectral curve $S$ of $(X,l)$ we have
$$d(S)=d(X)+2.$$
\endproposition

\proof
If $A_\R$ has $r\ge1$ connected components, then $d_A=d_X=5-r$, and we need to show that $d_S=7-r$.
Consider a line $l_1\subset X_1$ representing point $[l_1]$ of
a torus component $T_1\subset F_\R(X_1)$ and their continuous variations, $l_t\subset X_t$,
$T_t\subset F_\R(X_t)$, $[l_t]\in T_t$ along with the spectral curves $S_t\subset P^2$ associated to $(X_t,l_t)$.

If the torus $T_1$ is contracting, then it degenerates as $t\to 0$ to one of the circle components
of  $A_\R\subset F_\R(X_0)$, which implies, in particular, that
the limit line $l_0$ contains the node $\ss\in X_0$.
By our assumption,  $\ss$ has signature $(3,1)$, which means that the local quadric $Q_\ss$ is an ellipsoid and, thus,
the line generators of $Q_\ss$ are imaginary. This implies that the two nodes of $S$ are imaginary and, hence, the components
 $C_1,\dots,C_r$ of $A_\R$ are projected into $r$ smooth real components of $S_{0\R}$.
 Their smoothness implies that $S_{t\R}$ with $t>0$ has also $r$ components and so,  $d_{S}=d_{S_1}=7-r$.

If the torus $T_1$ is T-mergeable, then its limit $T_0$ is merging along some component $C_i$ of $A_\R\subset F_\R(X_0)$ with another torus, and
the path $l_t$ can be chosen so that $[l_0]\in C_i$.
 Since the node $\ss$ in this case is of signature $(2,2)$, the quadric $Q_\ss$ is a hyperboloid.
 Therefore, the two nodes of $S$ are real.

 None of these nodes is solitary or a self-intersection of
 the image of one of the real compensuredonents of $A_\R$. Indeed, otherwise,
the  torus $T_0$ is merging not with a torus but
with the non-orientable component of $\hat F_\R(X_0)$,  that is $$N(X_0)=A/\conj\cup_{j=1}^r\Sym^2(C_j)\subset\Sym^2(A)=\hat F(X_0).$$
This is because a solitary point of $S_0$ can appear only as the projection
of an imaginary pair, $\{z,\bar z\}$, and then
$[l_0]\in A/\conj$, while if a node is formed by
projection of some component
$C_j\subset A_\R$, then
$[l_0]\in \Sym^2(C_j)$.

Thus, the two nodes of $S_\R$ are intersection points of the projections $C_j'$ and $C_k'$ of two different components $C_j$ and $C_k$ of $A_\R$.
Therefore, the quintics $S_{t\R}, t>0$, which are obtained by smoothing the nodes of $S_\R$, have
either $r-1$ or $r$ connected components. The first option is excluded by Corollary \ref{discrepancy-relation}, and we conclude that
$d_{S_t}=7-r$, for $t>0$.
 \endproof

\subsection{Cubic threefolds of type $C^3_I$}\label{permutation-type-C3I}\label{3to6tori}
Here, we consider a real one-nodal cubic threefold $X_0$ with a node of signature $(2,2)$ such that the associated quadrocubic is of type {\bf 3}${}_I$
(see Subsections \ref{quadrocubics-edges} and \ref{real-quadrocubics-classification}). Such a cubic exists; indeed, one may start from constructing $A$ as an intersection of a hyperboloid $f_2=0$ with
three planes $L_1=0, L_2=0, L_3=0$ from a pencil of planes whose real base-line $L_1=L_2=L_3=0$ does not intersect $Q_\R$ (such a curve $A$ is of type I, see, for example, \cite{Zvon}),
and define $X_0$ by affine equation $f_2+L_1L_2L_3=0$.

\proposition\label{Fano-for one-nodal-M-2}
\roster
\item
$F_\R(X_0)$ consists of 4 connected components, one of which is a non-orientable smooth surface and each of the other components is the product of a figure-eight with $S^1$.
\item
There exists an equisingular deformation $\{X_t\}_{t\in[0,1]}$ of $X_0$ preserving the node of $X_0$ at a fixed point of $P^3$,
such that $X_1=X_0$ and whose monodromy gives a cyclic permutation of the three non-smooth components of $X_0$.
\endroster
\endproposition

\proof
(1) Let us denote by $C_i$, $i=1,2,3$, the connected components of $A_\R$.
By Proposition \ref{Fano_for_nodal_cubic}, the normalization $\hat F(X_0)$ of $F(X_0)$ is $\Sym^2(A)$. It implies that $\hat F_\R(X)$ consists of
one non-orientable component and
3 torus components $T_k$, $k=1,2,3$,  $T_k=C_i\times C_j$, $1\le i<j\le3$, $k\ne i,j$.
Each torus $T_k$ contains a pair of disjoint curves, $D_{ki}=\phi_i(C_k)$, $i=1,2$, which are identified
to form the double curve of $\hat F(X_0)$
(see Proposition \ref{Fano_for_nodal_cubic}).
Such an identification turns $T_k$
in a component of $F_\R(X_0)$  homeomorphic to the product of a figure-eight with $S^1$.

(2)
 By construction, $A_{0\R}$ splits into three disjoint ellipses $E_i=P_{i\R}^2\cap Q_\R$ traced by the planes $L_i=0, i=1,2,3$.
 Rotation of these planes around the base-line gives us a family of triples of planes $L_{1t}=0, L_{2t}=0, L_{3t}=0$ ending by a cyclic permutation $L_1=L_{10}=L_{31}, L_2=L_{20}=L_{11},L_3=L_{30}=L_{21}$,
 and we get the required family of cubic threefolds by defining it by equations $f_2+L_{1t}L_{2t}L_{3t}=0$.
\endproof

\proposition\label{mergeable-tori-of-C3I}
For a real cubic threefold $X$ of type $C^3_I$, all the tori of $F_\R(X)$ are T-mergeable.
More precisely, there exists a nodal degeneration of $X$ that yields merging of the six toric components of $F_\R(X)$ pairwise, along the three
connected components of the quadrocubic $A_\R$ associated to the node.
\endproposition

\proof
Due to Proposition \ref{Fano-for one-nodal-M-2}, the degeneration ensured by Theorem \ref {quadrocubics_as_edges} is as required.
\endproof

\subsection{Atoric real 6-nodal Segre cubics}\label{real6nodal}
\proposition\label{all_same} If $X$ is an atoric real 6-nodal Segre cubic threefold,
then the central component $R$ of $\hat F(X)=P_1^2\cup P_2^2\cup R$ is a real non-singular cubic $(M-1)$-surface, whereas $P_1^2, P_2^2$
are imaginary complex conjugate projective planes.
\endproposition

\proof Choose one of the nodes  $\ss\in X$. Since it is of signature $(3,1)$, the quadric $Q_\ss$ is an ellipsoid.
By Proposition \ref{simplicity-conditions} the quadrocubic $A$ splits into two nonsingular irreducible curves
$B_1$ and $B_2$ of bidegree $(2,1)$ and $(1,2)$, respectively, intersecting each other transversally at 5 points.
In our case these curves are imaginary complex conjugate, since the complex conjugation transforms each curve of bidegree $(a,b)$
into a curve of bidegree $(b,a)$. Since all the nodes of $X$ are real, all the 5 points of $B_1\cap B_2$ are real (see Proposition \ref{sing-quadricubic}).
Due to Proposition \ref{cubic-surfaces-from quadrics},
this implies that $R$ is an $(M-1)$ real cubic surface. The components  $P_1^2, P_2^2$
are imaginary complex conjugate,
since $P_i^2=\Sym^2(B_i)$ for each $i=1,2$ (see Proposition \ref{quadrocubics-on-quadrics}).
\endproof

Note that $R$ is  an $M$-surface and both $P_1^2, P_2^2$ are real, if the nodes are of signature $(2,2)$.

\subsection{Monodromy of collapsing tori}
\proposition\label{multi-nodal-perturbation}
For any non-singular real cubic threefold $X$
of type $C^k$, $1\le k\le5$, there exists a real degeneration $X_t$, $t\in[0,1]$, of $X=X_1$ to a $(k+1)$-nodal atoric cubic $X_0$.
 The corresponding degeneration of the Fano surfaces $F(X_t), t>0,$ contracts the $\binom{k+1}2$
tori of its real locus to the points  $[l_{ij}]\in F(X_0)$, where $l_{ij}$ are the lines connecting the nodes $\ss_i$ and $\ss_j$,
$0\le i<j\le k$.
\endproposition

\proof
Consider a non-singular real quadric $Q$
whose real locus $Q_\R$ is an ellipsoid
and a real quadrocubic $A$ whose real locus $A_\R$ is formed by $k$ solitary nodes on $Q_\R$ in linearly
general position (one can construct such a quadrocubic as, for example, the image of a real plane quartic with $k-1$ solitary nodes under a standard
birational transformation from real projective plane to a real ellipsoid, choosing the two base points of this transformation to be imaginary complex conjugate points in general position).
Next as usual we pass from such a $k$-nodal quadrocubic $A$ to a $(k+1)$-nodal cubic $X_0\subset P^4$ (see Proposition \ref{cubics-quadrocubics-association})
so that $A$ is associated to one node, $\ss_0$, of $X_0$ and the other nodes $\ss_i$, $i=1,\dots, k$ of $X_0$ represent the nodes of $A$ (viewed as lines $\ss_0\ss_i$).
Lemma \ref{nodes_coherence} guarantees that
all the nodes of $X_0$ are real and have signature $(3,1)$.

Instead of a direct perturbation of all the nodes at once, it will be convenient to perform a perturbation at two steps: first, we keep the distinguished node $\ss_0$ of $X_0$ and perturb the others, which gives one nodal cubic $X_0'$,
and, second, we perturb the remaining node $\ss_0$ and obtain non-singular cubic $X$.
More precisely,
at the first step, we perturb the quadrocubic  so that each solitary node $\ss_0\ss_i$ of $A$ gives birth to an oval $C_i\subset Q_\R$ of the quadrocubic $A'$ associated to $X_0'$,
and then choose and ascending perturbation of $X_0'$,
which implies that the resulting perturbation of $X_0$ is ascending too and thus, according to Theorem \ref{quadrocubics_as_edges}, really yields  $X$ of type $C^k$ as required.
 The required perturbation scheme really exists because in the space of cubic threefolds $\Cal C^{3,3}$
 each stratum of $k$-nodal ones is an intersection of $k$ transversal  branches of the discriminant $\Delta^{3,3}$, for each $k\le 6$.

As it was remarked in Subsection \ref{real-Fano-subsection},
the Fano surface $F_\R(X_0')$
has torus components
$T_{ij}=C_i\times C_j$, $1\le i<j\le k$, as a part of its non-singular locus,
and solitary circle components $C_i\subset A'_\R$, $i=1,\dots,k$, forming its real singular locus.
An ascending perturbation of $X_0'$ leads to smoothing of circles $C_i$ into
toric components $T_{0i}\subset F_\R(X)$, $i=1,\dots, k$.
This means that the tori $T_{0i}$ are collapsing and contract to the points $l_{0i}$ as $X$ contracts to $X_0$.
Since degeneration of $X_0'$ to $X_0$ contracts the oval $C_i$ to points $\ss_0\ss_i$, the other
other toric components of $F_\R(X)$ corresponding to
$T_{ij}=C_i\times C_j$, $1\le i<j\le k$, will contract as the corresponding ovals, namely,
to the points $l_{ij}\in F_\R(X_0)$ representing lines $\ss_i\ss_j$.
\endproof

\proposition\label{toric-monodromy}
Let $X^\tau_0$, ${\tau\in[0,1]}$,  be a continous loop in the space of atoric $k$-nodal cubic threefolds with $k\le 6$,
and let $\sigma$ be the resulting permutation of the nodes $\ss_1,\dots,\ss_k$ of $X^0_0$. Then, there exists a continuous map
$[0,1]^2\to  \Cal C_{3,3}$, $(t,\tau)\mapsto X_t^\tau$,
such that :
\roster
\item $X^0_t=X^1_t$ for any $t$ and $X_t^\tau$ is non-singular for any $\tau\ge 0$ and any $t>0$;
\item $X^0_t$, $t\in[0,1]$, is an ascending perturbation of $X^0_0$;
\item the torus component $T_{ij}\subset F_\R(X^0_1)$, $1\le i<j\le k$,  that contracts
to the point $[\ss_i\ss_j]\in F(X^0_0)$ is mapped by the monodromy along the loop $X^\tau_1$, ${\tau\in[0,1]}$,
to the torus component $T_{\sigma(i)\sigma(j)}$.
\endroster
\endproposition

\proof  Straightforward consequence of Proposition \ref{ascending-coherent}, since for any $k\le 6$ each of the strata of nodal cubic threefolds with $k$ nodes is an intersection of $k$ transversal  branches of the discriminant.
\endproof

\corollary\label{collapsibility-transitivity}
If $X$ is a real non-singular cubic threefold of type $C^k$, $1\le k\le5$, then
all the torus components of $F_\R(X)$ are collapsible.
If, in addition, $k<5$ then the
Fano monodromy group of $X$ acts transitively on the
torus components of $F(X)$.
\endcorollary

\proof
Proposition \ref{multi-nodal-perturbation} implies the first statement.
The second statement follows from
Proposition \ref {toric-monodromy}(3) and
Corollary \ref{monodromy-of-nodes}.
\endproof


\section{Proof of main theorems}\label{summary}

\subsection{Proof of Theorem \ref{N-correspondence}}
If $X$ belongs to the type $C^0$ or $C^1_{I(2)}$, then its Fano surface has only one connected
component, $F_\R(X)=N_5$ (Subsection \ref{real_Fano}, Table 1).
Hence, there exists only one spectral matching $([X],[S])$ involving such a type
of cubic threefolds.
This matching is perfect,
since  due to Lemma \ref{perfect-are-spectral} the perfect matching is spectral.

For $X$ of type  $C^1_I$, Theorem \ref{N-correspondence}
is proved in Subsection \ref{N-correspondence-for exotic}.

If $X$ is of type $C^k$, $1\le k\le 5$, then, by  Proposition \ref{multi-nodal-perturbation} and Corollary \ref{collapsibility-transitivity},
all the tori of $F_\R(X)$ are collapsible.  Thus, by Proposition \ref{collapsible-mergeable},
the spectral matching $([X],[\spec])$ is skew if $l$ is chosen on any of
the tori. Therefore, for the perfect mathching, which is always spectral due to Lemma \ref{perfect-are-spectral},
there remain only the option $l\in N(X)$.

In the remaining case $C^3_I$, Proposition \ref{mergeable-tori-of-C3I} shows that
the tori are mergeable. Thus, by Proposition \ref{collapsible-mergeable}, they all provide
skew matchings.
The remaining component $N(X)$ provides
a perfect matching, as it follows again from Lemma \ref{perfect-are-spectral}.
\qed

\subsection{Proof of Theorem \ref{T-orbits}(1)}\label{T-orbits-proof-1}
For real cubic threefolds of type $C^k$, $2\le k\le 4$, this claim of the
theorem follows from Corollary \ref{collapsibility-transitivity}, while
for other types it is empty:
the types $C^5$ and $C^3_I$ do not make part of the statement,
and in the remaining cases the Fano surface has at most one toric component.
\qed

\subsection{Segre 3-3 bipartition}\label{bipartition}
Claims (1) and (2) of the following proposition go back to L.~Schl\"afli, claim  (3)
is due to B.~Segre \cite{Seg}, \textsection 41--42.

\proposition\label{lines-partition}
Assume that $Y$ is a non-singular real cubic (M-1)-surface. Then:
\roster\item
$Y$ has precisely 15 real and 12 imaginary lines, where the latter ones
form a double six, $\binom{L_1\dots L_6}{\bar L_1\dots \bar L_6}$, so that each line $\bar L_i$ is complex-conjugate to $L_i$.
\item
Each of the real lines intersects precisely two pairs of conjugate imaginary lines, $(L_i,\bar L_i)$ and
$(L_j,\bar L_j)$, so that
this gives a one-to-one correspondence between the 15 real lines and the $\binom62$ choices of  $1\le i<j\le6$.
\item
The set of the six pairs $(L_i,\bar L_i)$
is  partitioned into two triples, so that, with respect to the above correspondence,
the 9 combinations of choices of $(L_i,\bar L_i)$ and $(L_i,\bar L_i)$ from different triples represent hyperbolic lines, while the 6 other combinations
represent elliptic lines.
\qed\endroster
\endproposition

If we consider an atoric real 6-nodal Segre cubic threefold $X$ and the lines on the cubic surface
$Y=R\subset \hat F(X)$, then
Proposition \ref{quadrocubics-on-quadrics}  and Corollary \ref{nodes-as-lines} transform
the partition of the set of 6 pairs $(L_i,\bar L_i)$  given by Proposition \ref{lines-partition}(3)
into a bipartition of the set of the nodes of $X$. We name
the latter bipartition the {\it Segre 3+3 bipartition of nodes}.
A pair of nodes $(\ss_i,\ss_j)$ of nodes of $X$ is called {\it hyperbolic} or {\it elliptic},
in accord with the type of the corresponding line of $Y=R$.

\proposition\label{from3+3to3+2}
Let $\ss_0,\dots,\ss_5$ be the nodes of an atoric real Segre 6-nodal cubic threefold $X$
and  $A$ the quadrocubic associated with $\ss_0$.
Then the Segre  3+3 bipartition of $\{\ss_0,\dots,\ss_5\}$
induces
a 3+2 partition of
$\{\ss_1,\dots,\ss_5\}$ that matches with
the Segre  3+2 partition of the nodes of  $A$ via  the correspondence
$\ss_i\mapsto [\ss_0\ss_i]$,
where $[\ss_0\ss_i]$ denotes the node of $A$ represented by the line $ \ss_0\ss_i$.
\endproposition

\proof
By Propositions \ref{quadrocubics-on-quadrics} and \ref{cubic-surfaces-from quadrics}, the normalized central component $R\subset\hat F(X)$
can be identified with the local quadric $Q_p$
blown up at the points $q_i=[\ss_0\ss_i]$, $i=1,\dots,5$.
According to Proposition \ref{2+3partition},
the Segre 2+3 bipartition of $\{q_0,\dots, q_5\}$  is formed by a pair
$\{q_1,q_2\}$ such that the corresponding exceptional divisors $E_1$ and $E_2$
are elliptic lines  in $R$,  while for the remaining triple $\{q_3,q_4,q_5\}$, the corresponding $E_i$, $i=3,4,5$ are hyperbolic lines.
Each line $E_i$ is incident to the proper image of the imaginary generators of $Q_p$ passing through $q_i$ and to the components
$B_1$ and $B_2$ of the quadrocubic $A=B_1\cup B_2$.
Ellipticity of $E_i$ means, by definition of the Segre 3+3 bipartition,  that the pair
$\ss_i$ and $\ss_0$ is elliptic,  and so such a pair is contained in the same triple. This implies that
$\{\ss_0,\ss_1,\ss_2\}$ form one triple,  and $\{\ss_3,\ss_4,\ss_5\}$ form the other one.
\endproof

\subsection{Monodromy of nodes of atoric real Segre 6-nodal cubics}\label{monodromy-for-M-cubics}
Throughout this Subsection $X$ is assumed to be an atoric real Segre 6-nodal cubic threefold. The nodes of $X$ are denoted  by
$\ss_0,\dots,\ss_5$.

\proposition\label{bipartition-invariance}
The Segre 3+3 bipartition of the set of nodes of
$X$ is preserved under
real equisingular deformations, so that
any monodromy permutation induced on the set of nodes of $X$
either preserves each triple of the bipartition or interchanges them.
\endproposition

\proof
Since $R$ remains  a non-singular $(M-1)$-cubic surface under real equisingular deformations of $X$,
the result follows from
Proposition \ref{lines-partition}(3) and invariance of elliptic and hyperbolic lines under real deformations
of a real cubic surface.
\endproof

The converse is also true.

\proposition\label{bipartition-monodromy}
Any permutation of the nodes of $X$
that preserves or interchanges the triples of the Segre 3+3 bipartition of the set of nodes
can be realized as a monodromy  under  a real equisingular deformation of $X$.
\endproposition

First, we observe transitivity of the monodromy action.

\lemma\label{transitivity-on-nodes}
Real equisingular deformations induce a transitive monodromy action on the set of nodes of $X$.
\endlemma

\proof
We start with a real projective transformation $g\:P^4\to P^4$ that sends $\ss_0$ to any other given node $\ss_i$,
and induces an isomorphism between the local quadric $Q_0\subset P^3_{\ss_0}$ and the local quadric $Q_i\subset P^3_{\ss_0}$.
The latter can be done since any pair of ellipsoids can be connected by a real projective equivalence.
The 5-point configuration of the nodes of the quadrocubic $A_i\subset Q_i$ associated to $\ss_i\in X$
and that of the nodes of
$g(A_0)\subset g(Q_0)=Q_i$ associated to $\ss_i\in X'=g(X)$
can be connected by a deformation
in the class of
5-point subsets of $Q_i$ in linearly general position,
see Corollary \ref{5-point-monodromy}.
Next note that
5 points on $Q_i$ in linearly general position
determine a unique curve of degree $(2,1)$ and a unique curve of degree $(1,2)$
that pass through the given points. These two curves form together a quadrocubic with the 5 given nodes and depend continuously on them.
Now, it remains to apply Lemma \ref{nodes-to-quadrocubics} to obtain an equisingular real deformation connecting the cubic $X$
with itself and moving $\ss_0$ to $\ss_i$.
\endproof

\demo{Proof of Proposition \ref{bipartition-monodromy}}
Lemma \ref{transitivity-on-nodes} reduces the
task to the case of permutations of the nodes $\ss_i$, $i=0,\dots,5$,
that preserve a given node, say $\ss_0$. Let $p_i\in A_0\subset Q_{\ss_0}$ denote the nodes of the quadrocubic
$A_0$
that correspond
to the directions $\ss_0\ss_i$.
By Proposition \ref{2+3partition} and Corollary \ref{5-point-monodromy}
any permutation of the set $\{p_1,\dots,p_5\}$ that preserves the induced Segre 3+2 bipartition
can be realized by deformations in the class of linearly generic 5-point configurations.
Using Proposition \ref{from3+3to3+2},
we can lift such a deformation to a deformation of 6-nodal Segre cubics (like in
the proof of Lemma \ref{transitivity-on-nodes}).
\qed\enddemo

\subsection{Proof of Theorem \ref{T-orbits}(2)}\label{Proof-T-orbits-for-M-cubics}
By Proposition  \ref{multi-nodal-perturbation}, it is sufficient to consider a
cubic $X_1\subset P^4$ of type $C^5$ obtained by a perturbation $X_t$, $t\in[0,1]$, of an atoric real
6-nodal Segre cubic $X_0$.
A monodromy permutation of the nodes $\ss_i$, $i=0,\dots,5$, of $X_0$ realized by a real equisingular deformation of $X_0$ induces
a permutation of the lines $l_{ij}$ connecting the nodes $\ss_i$ and $\ss_j$, $0\le i<j\le 5$, which in its turn provides the same permutation
of the corresponding tori $T_{ij}$,
as it follows from Proposition \ref{toric-monodromy}.
On the other hand,
Propositions \ref{bipartition-invariance} and \ref{bipartition-monodromy} show that
monodromy permutations of the nodes
$\ss_i$ are precisely the ones preserving the Segre 3+3 partition. The action induced by such permutations on the pairs  $(\ss_i, \ss_j)$ produces 2 orbits:
for one orbit $\ss_i$ and $\ss_j$ belong to the same triple of the partition, and
for the other orbit they belong to different triples.

Now it remains to characterise the corresponding two orbits of $T_{ij}\in \TT(X)$ as
$\TT_I(X)$ and $\TT_{II}(X)$, according to their definition in Subsection \ref{monodromy_groups}.
To do it, note that
the spectral curve $S_{0}$ corresponding to the line $l_0=\ss_i\ss_j$ on $X_0$
splits, due to Proposition \ref{simplicity-conditions}(4),
into a pair of conjugate imaginary conics, $C_{ij}$, $C_{ij}'$, intersecting at 4 real points and a real line $L$
(reality of the 4 intersection points is due to their correspondence to the four real planes, $\ss_i\ss_j\ss_k$, $k\in\{0,\dots,5\}\sm\{i,j\}$).
By Proposition \ref{2+3partition}, these 4 points are in a non-convex position in
$\Rp2\sm L$ for the first orbit of $\{s_i,s_j\}$,
and convex for the second orbit.
On the other hand, according to Theorem \ref{N-correspondence}, the
spectral curve $S_1$ of a pair
$(X_1,l_1)$ obtained by a perturbation of $(X_0,l_0)$, is itself a perturbation of $S_{0}$
which gives an $(M-2)$-quintic.
Hence, each of the 4 real nodes of $S_0$ gives birth to an oval, whereas the line $L$ is perturbed into the J-component of $S_1$.
Now it is left to apply Lemma \ref{convexity-criterion}
to distinguish quintics of types $J\+4_I$ and $J\+4_{II}$.
\qed

\subsection{Proof of Theorems \ref{T-orbits}(3)}\label{Case-C3I}
By Theorems \ref{quadrocubics_as_edges} and \ref{DZ}(3), a cubic $X\subset P^4$ of type $C^3_I$ can be obtained by an ascending perturbation
from a real one-nodal cubic $X_0$ whose associated quadrocubic $A$
has real locus formed by three non-contractible components.
Recall that $F_\R(X_0)$ contains 6 tori that split into 3 pairs: in each pair the tori are merged together along the corresponding component
of $A_\R$, see Proposition \ref{mergeable-tori-of-C3I}. The induced monodromy action cyclically permutes these 3 merged pairs and after an ascending perturbation
we obtain an induced action of $\Z/3$ on the set of the six tori of $F_\R(X)$.
This implies that there are at most 2 monodromy orbits on the set of tori, each orbit formed by 3 cyclically permuted tori.
On the other hand, since the (M-4)-quintic $S$ can be of Klein's type I or II, and by Lemma \ref{skew-matchings-are-spectral} there are spectral matchings
$([X],[S])$ with the both types, and they are realised by lines $l$ on the toric components of $F(X)$ according to Theorem \ref{N-correspondence},
we conclude that
there exist two monodromy orbits on the set
of tori:
one orbit of type I, and another of type II.
\qed

\subsection{Summary of the proof of Theorem \ref{main} }
Part (1) of Theorem \ref{main} is derived in Subsection \ref{scheme-main} and
part (3) in Subsection \ref{surjectivity of SMC}.
Theorem \ref{N-correspondence} together with uniqueness of the odd component
 $N(X)\subset F_\R(X)$ implies (2) in the case of perfect matchings.
By Theorem \ref{T-orbits}(1), pairs $(X,l)$ that
represent a given skew matching $([X],[S])$, are deformation equivalent unless
$[X]$ is $C^5$ or $C^3_I$. In the case of $X$ having type $C^5$ (respectively, $C^3_I$) we apply Theorem \ref{T-orbits}(2)
(respectively, \ref{T-orbits}(3)) showing that the deformation classes of pairs $(X,l)$ (within the given class $[X]$) are distinguished by
$[S]$.


\section{Concluding remarks}

\subsection{Expanded spectral correspondence}\label{full_correspondence}
In this paper it was sufficient for us to stay within the most classical setting not involving
discussion of theta-characteristics on singular curves and to use only the most basic properties of the spectral correspondence
understood as a regular bijective morphism:
$\Cal C^*/PGL(5; \C)$ $ \to {\Cal S}/PGL(3; \C)$ when working over $\C$ and $\Cal C_\R^*/PGL(5; \R)\to {\Cal S_\R}/PGL(3; \R)$
over $\R$ (see Subsection \ref{spectral-cor}).
Both morphisms are indeed isomorphisms. Furthermore,
one can extend $\Cal S$
to a larger space $\Cal {\hat S}$ by
adding
nodal quintics and
define a theta-characteristic on a nodal quintic not as a line bundle but as a rank-1 torsion free sheaf that is not locally free at each of the nodes
({\it cf.}, \cite{Co}).
Then, these isomorphisms
extend to isomorphisms $\Cal C/PGL(5; \C)\cong {\Cal {\hat S}}/PGL(3; \C)$
and  $\Cal C_\R/PGL(5; \R)=(\Cal C/PGL(5; \C))_\R\cong ({\Cal {\hat S}}/PGL(3; \C))_\R= {\Cal {\hat S_\R}}/PGL(3; \R)$, respectively.

\subsection{Monodromy groups} Using the approach and the results of this paper it is possible to evaluate the Fano real component monodromy group $G_{\roman{mon}}$
(see Subsection \ref{monodromy_groups})
for each of the deformation classes.

For example, in the case of maximal real cubic threefolds this group is isomorphic to the
semi-direct product of $S_3\times S_3$ and $\Z/2$, namely, to the
subgroup $S_{123, 456}$ of the symmetric group $S_6$ on 6 elements $1, 2,\dots, 6$
formed by permutations leaving invariant
the partition $\{1,\dots,6\}=\{1,2,3\}\cup\{4,5,6\}$ (with possibility to interchange 3-element subsets).
Indeed, as we have seen in Section \ref{monodromy-for-M-cubics} via degenerating maximal cubic threefolds to a 6-nodal Segre cubic,
one can enumerate
the real tori, $T_{ij}=T_{ji}$, $1\le i\ne j\le 6$, of the Fano surface so that
each element
of
$G_{\roman{mon}}$ preserves the 6+9 partition of the set of tori in $\{T_{ij} \,\text{with}\,  i,j\le 3\, \text{or} \, i,j\ge 4\}$ and
 $\{T_{ij} \, \text{with} \, i\le 3\, \text{and} \, j\ge 4\}$, and, in addition,  $G_{\roman{mon}}$
 contains as a subgroup the group $S_{123, 456}$ naturally acting on the set of tori. Furthermore, one can check
that
 $H_1(F; \Z/2)$ contains 6 independent elements, $c_1,\dots, c_6$, such that
 the intersection of the images of $H_1(T_{ij}; \Z/2)$ and $H_1(T_{kl}; \Z/2)$ in $H_1(F; \Z/2)$ is zero, if $i,j$ and $k,l$ have no common element, and generated by $c_i$, if $k$ or $l$ is equal to $i$.
Clearly, the action of $G_{\roman{mon}}$ preserves these relations.
Hence, each element of $G_{\roman{mon}}$ sends tori having a common index to tori having a common index, and the result follows.

Passing through one-nodal degeneration from cubics in class $C^5$ (M-cubics)  to cubics of class  $C^4$, from  $C^4$  to $C^3$, {\it etc.}, and noticing that the real structure of the Fano surfaces is modified only
along a loci not containing the surviving real tori,
one concludes that the above exchange property in $H_1$ is preserved for the remaining tori and, applying Proposition \ref{k-point-monodromy}, obtains that for cubics of class  $\Cal C^{5-k}$ with $k\ge 1$
the group $G_{\roman{mon}}$ is isomorphic to $S_{6-k}$.

%
\subsection{Conormal projection}\label{conormal-remark}
The conormal projection used in Section \ref{conormal-section}
for studying Fano surfaces can be also used
for studying real plane
quintic curves. Here is one example that concerns M-quintics.

First, we recall a simple {\it convexity property} of the 6 ovals of an M-quintic: if we pick
one point inside each oval,
they will form in some affine plane $\R^2\subset P^2_\R$ a convex hexagon disjoint
from the one-sided component of the quintics
(see Figure \ref{hexagonal-bitangents}c)).
Our observation is that
for each generic
M-quintic $S$ there exist 6 and only 6 real bitangents whose tangency points are either imaginary,
or both real and belonging to the same connected component of $S_\R$;
these bitangents are in a natural correspondence with the six ovals of $S$:
each bitangent separates the corresponding oval from the other five by cutting an angle from the
hexagon described
(see Figure \ref{hexagonal-bitangents}(d)).

\midinsert
\hskip3mm\epsfbox{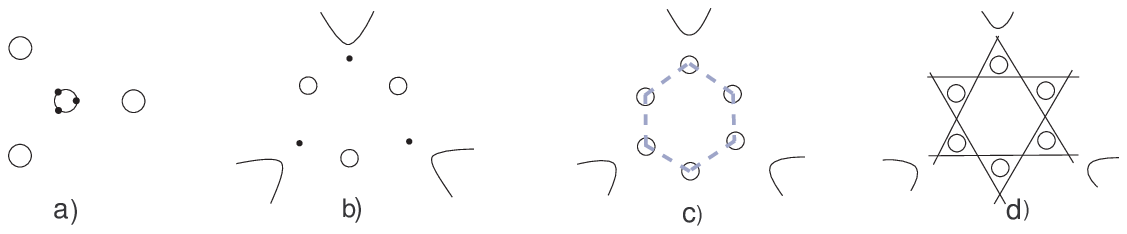}
\vskip-4mm
\figure{a) Quartic $A$ with 3 fundamental points of the quadratic transformation; b) quintic $S_0$; c) M-quintic $S$ with a hexagonal configuration of ovals;
d) six bitangents at imaginary points.}\endfigure
\endinsert\label{hexagonal-bitangents}

For the proof, we treat, first, some special example.
We start from considering
the real quartic $A_0$ defined (for better visibility) in polar coordinates as
$r^2-2r^3\cos 3\phi + r^4=0$.
Its real locus consists of
4 solitary nodes: at the origin ($r=0$) and at
 three points $s_k$, $k=0,1,2$, with $r=1$ lying on the rays $R_k=\{\phi=\frac{2k}3\pi, r\ge0\}$.
A small real perturbation $r^2-2r^3\cos 3\phi + r^4=\epsilon^2$
is a quartic $A$ (shown on Figure \ref{hexagonal-bitangents}a)) that has 4 small ovals around the nodes of $A_0$.
Next, we consider the intersection points of the rays $R_k$, $k=0,1,2$, with
the central oval of $A$ and perform the standard quadratic Cremona transformation with these three points as fundamental points.
The image of the quartic $A$
will be a quintic $S_0$ with 3 ovals that are the images of the 3 non-central ovals of $A$,
a one-sided component $J$ that is the image of the central oval, and 3 solitary points that
are the images of the pairs of imaginary points of $A$ that
lie on the sides of the fundamental triangle (see Figure \ref{hexagonal-bitangents}b)).
Each of the conics $\Theta_0$ that pass through these 3 solitary points and intersect two of these 3 ovals contains the remaining oval inside. Finally, we consider an M-quintic $S$ obtained by a small perturbation of $S_0$ and the theta-conic $\Theta$ obtained by a small perturbation of $\Theta_0$ whose real locus
contains inside all the 6 ovals of $S$ and is tangent to 5 of them.
The curve dual to $S$ is shown on Figure \ref{dual-projection}.

Let $X$ and $l$ be a real cubic threefold and a real line on it defined by $(S,\Theta)$ as above, and let $F$ be the Fano surface of $X$.
 By Proposition \ref{folds-near-ovals} the image $\l_\R(F_{l\R})$
of the conormal projection $\l_{\R}\:F_{l\R}\to\widehat{\P_\R^2}$ is the complement of the six bands on $\widehat{\P^2_\R}$
formed by the points dual to the lines in $\P^2_\R$ crossing the six ovals of $S$.
\midinsert
\hskip40mm\epsfbox{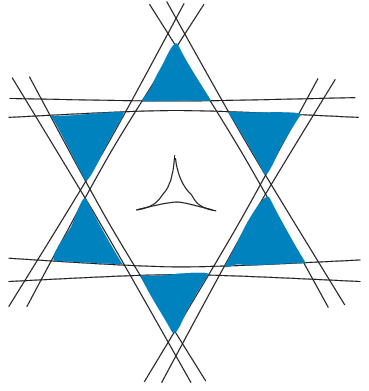}
\vskip-4mm
\figure{
The dual curve, $\hat S$, and the 6 triangle components of $\l_\R(F_{l\R})$.}\endfigure\label{dual-projection}
\endinsert
This complement consists of
16 curvilinear polygons: 6 triangles, 9 quadrilaterals and a hexagon, as it is shown on Figure \ref{dual-projection}.
So, the 16 connected components of $F_{l\R}$ are mapped into the above 16 polygons, each in its own.
The projection to the triangles and quadrilaterals have 4 sheets folded at the boundary and
the Euler characteristic argument shows that
the corresponding components of $F_{l\R}$ must be tori, while the non-orientable component of $F_{l\R}$ is mapped to the hexagon.
Moreover, the projection
of a torus to a triangle must have a branch point. It can be easily seen that
such a point is dual to a real bitangent to $S$ with imaginary tangency points,
and the projective duality implies that
this bitangent is located in $\P_\R^2$ as is stated.

To deduce the general claim, it is sufficient to notice that all
M-quintics are deformation equivalent and that a bitangent tangent to two distinct components of the quintic can not disappear
during a deformation.

\subsection{Mutual position of $S_\R$ and $\Theta_\R$}
The methods developed allows to obtain some restrictions on the position of the ovals of a real quintic $S$ with respect to a contact conic $\Theta$.
The fact that the associated theta-characteristic differs from Rokhlin's
one implies for instance that
the ovals of an M-quintic $S$ cannot lie all outside the interior of $\Theta_\R$ and have each even contact with $\Theta$
(for instance, it can not happen that each of the 5 tangency points $S\cap \Theta$ is imaginary or lies on the $J$-component
of $S_\R$).

Indeed, Corollary \ref{Skew-matching-criterion} together with the fact that every contact pair $(S,\Theta)$ is spectral (Theorem \ref{tang-conic})
shows that such a mutual position of $S_\R$ and $\Theta_\R$ implies
that for the corresponding cubic threefold $X$ we have $d_X=d_S-2$ (the matching $(X,S)$ is skew)
 which is impossible in the case $d_S=0$.
The latter argument shows also impossibility of a similar mutual position for $(M-1)$-quintics too.

Using that the spectral theta-characteristic must be odd, one can obtain some additional restrictions.
For example, if all the ovals of an M-quintic $S$ lie outside
$\Theta_\R$
and only one of them is tangent to $\Theta_\R$,
then the latter cannot be one of the three
ovals joined by a vanishing bridge-cycle with the one-sided component of $S_\R$.



\Refs\widestnumber\key{ABCD}

\ref{Ab} \paper Cubic surfaces with a double line
\by Sh.~Abhyankar
 \jour Memoirs of the College of Sciences, Univ. Kyoto, series A
 \vol 32
 \yr 1960
\pages 455 -- 511
\endref\label{Ab}

\ref{ACT} \paper Hyperbolic geometry and moduli of real cubic surfaces
\by D.~Allcock, J.A.~Carlson, D.~Toledo
 \jour Ann. Sci. Ãc. Norm. SupÃ©r. (4)
 \vol 43
 \yr 2010
\pages 69 --115
\endref\label{ACT}



\ref{At} \paper Riemann surfaces and spin structures
\by M.Atiyah
 \jour Ann. Sci. ENS
 \vol 4
 \yr 1971
\pages 47--62
\endref\label{Atiyah}

\ref{B1} \paper Vari\'et\'es de Prymes et jacobiennes interm\'ediares
\by A.~Beauville
 \jour Ann. Sci. ENS
 \vol 10
 \yr 1977
\pages 309 - 391
\endref\label{Beau}

\ref{B2} \paper Determinantal hypersurfaces
\by A.~Beauville
 \jour Michigan Math. J.
 \vol 48
 \yr 2000
\pages 39 - 64
\endref\label{Beau-Det}



\ref{CG} \paper The intermediate Jacobian of the cubic threefold
\by C. Clemens,   P. Griffiths
 \jour Ann. of Math.
 \vol 95
 \yr 1972
\pages 281-356
\endref\label{Clemens-Griffiths}

\ref{Co}
\by M.Cornalba.
\paper Moduli of curves and theta-characteristics.
\book
Lectures  on  Riemann  surfaces (Trieste, 1987)
\bookinfo M. Cornalba,
ed. X.Gomez-Mont, and A.Verjovsky,  World Sci. Publ.
\pages  560-589
\yr 1989
\endref\label{Co}

\ref{DK}
\paper  Topological properties of real
algebraic varieties : du c\^{o}t\`{e} de chez Rokhlin
\by A.~Degtyarev and V.~Kharlamov
\jour
Uspekhi Mat. Nauk.
\vol 55
\yr 2000
\issue 4
\pages  129--212
\endref\label{DK}

\ref{DZ}
\paper  Rigid isotopy classification of real algebraic curves of bidegree $(3,3)$
on quadrics
\by A.~Degtyarev and V.~Zvonilov
\jour
Mat. Zametki
\vol 66
\yr 1999
\issue 6
\pages  810-815
\endref\label{Degt-Zvon}



\ref{Do}
\by I.~Dolgachev
\paper Corrado Segre and nodal cubic threefolds
\inbook  From classical to modern algebraic geometry, ,
Trends Hist. Sci., Birkhäuser/Springer, Cham
 \yr2016
 \pages429--450
\endref\label{Dolgachev-Segre}

\ref{Fa} \paper  Sul sistema $\infty^2$ di rette contenuto in une
variet\`a cubica generale dello spazio a quattro dimensioni
 \by  G.~Fano
\jour Atti R. Accad. Sci. Torino
\vol 39
 \yr 1904
\pages 778Â-792
\endref\label{Fano}


\ref{Fi}
\by T.~Fiedler
\paper Pencils of lines and the topology of real algebraic curves
\jour Math. USSR-Izv
\vol 21
\yr 1983
\pages  161 -- 170
\endref\label{Fiedler}




%

\ref{FK}
 \by S. Finashin, V. Kharlamov
 \paper Topology of real cubic
fourfolds
 \jour J. Topol.
 \vol 3
 \issue 1
 \yr 2010
 \pages 1 - 28
\endref\label{FK}


\ref{GS}
\by D.~Gudkov, E.~Shustin
\paper
On the intersection of close algebraic curves
\inbook Lecture Notes Math. 1060
\yr 1984
\pages  278 - 289
\endref\label{GS}

\ref{GLS}
\by G-M.~Greuel, C.~Lossen, E.~Shustin
\book Introduction to Singularities and Deformations
\yr 2007
\bookinfo Springer-Verlag Berlin Heidelberg
\endref\label{GLS}

%





\ref{Ha}
 \paper Tangent Conics at Quartic Surfaces and Conics in Quartic Double Solids
 \by I.~Hadan
 \jour  Math. Nachr.
 \vol 210
 \yr 2000
 \pages 127 -- 162
\endref\label{Hadan}

\ref{HT}
 \paper Flops on holomorphic symplectic fourfolds and determinantal cubic hypersurfaces
 \by B.~Hasset, Y.~Tschinkel
 \jour   J. Inst. Math. Jussieu
 \vol 6
 \yr 2010
 \pages 125 -- 153
\endref\label{Hasset-T}




\ref{J}
\by D. Johnson
\paper Spin structures and quadratic forms on surfaces,
\jour J. London Math Soc.
\vol 2 (22)
\yr 1980
\pages 365-373
\endref\label{Johnson}



\ref{Kh}
\by V.~Kharlamov
\jour Funct. Anal. Appl.
\paper Rigid isotopy classification of real plane curves of degree 5
\vol 15
\issue 1
\yr 1981
\pages 73 - 74
\endref\label{Kh-quintics}


\ref{Kr1}
 \by V.~Krasnov
 \paper Rigid isotopy classification of real three-dimensional
 cubics
 \jour Izvestiya: Mathematics
 \vol 70
 \issue 4
 \yr2006
 \pages 731--768
\endref\label{Kr-rigid}

\ref{Kr2}
 \by V.~Krasnov
\paper
Topological classification of real three-dimensional cubics
\jour Mathematical Notes
\vol 85
\issue 5--6
\yr 2009
\pages 	841--847
\endref\label{Kr-top}

\ref{Kr3}
\by V.~Krasnov
\paper The topological classification of Fano surfaces of real three-dimensional cubics
\jour Izv. RAN. Ser. Mat.
\yr 2007
\vol 71
\issue 5
\pages 3--36
\transl
\jour Izv. Math.
\yr 2007
\vol 71
\issue 5
\pages 863--894
\endref\label{Kr-Fano}

%

\ref{Mum}
\by
D.~Mumford
\paper Theta characteristics of an algebraic curve
\jour Ann. Sci. ENS SÃ©r. 4
\vol 4 (2)
\pages 181--192
\yr 1971
\endref\label{Mumford}

\ref{Mur}
\by
J.P.~Murre
\paper Algebraic equivalence modulo rational equivalence on a cubic threefold
\jour Compositio Mathematica
\vol 25
\pages 161--206
\yr 1972
\endref\label{Murre}


\ref{Ro}
 \by V.~A.~Rokhlin
 \paper Proof of a conjecture of Gudkov
 \jour Functional Analysis and Its Applications
 \vol 6
  \yr 1972
 \pages 136--138
\endref\label{Ro}



\ref{Sch}
 \by L.~Schl\"afli
 \paper An attempt to determine the twenty-seven lines upon a surface of the third order,
 and to divide such surfaces  into species in reference to the reality of the lines upon
 a surface
  \jour Quart. J. Pure Applied Math.
 \vol 2
 \yr 1858
 \pages 110 - 120
\endref\label{Sch}


\ref{SeB}
 \by B.~Segre
\book The Non-Singular Cubic Surfaces. A new method of Investigation with Special Reference to
Questions of Reality
\bookinfo Oxford Univ. Press, London
 \yr 1942
\endref\label{Seg}

\ref{SeC}
 \by C.~Segre
 \paper Sulle variet\'a cubiche dello spazio a quattro dimensioni e su certi sistemi di rette e certe superficie dello
spazio ordinario
\jour Mem. Acad. Sci. Torino (2)
 \yr 1887
 \vol 39
 \pages 2 -- 48
\endref\label{Corrado}

\ref{T}
\paper Five lectures on three-dimensional varieties
 \by A.N.~Tyurin
 \jour Russian Math. Surveys
 \vol 27
 \yr 1972
 \issue 5
 \pages 1 - 53
\endref\label{Tyu}



\ref{Wa}
\by C.T.C.~Wall
\paper Sextic curves and quartic surfaces with higher singularities
\jour preprint
\yr 1999
\pages 1--32
\endref\label{Wall}

\ref{Wh}
\by F.~P.~White
\paper On the 5-tangent conics of a plane quintic curve
\jour Proc. Lond. Math. Soc.
\yr 1930
\pages 347--358
\endref\label{White}


\ref{Zv}
\by V.~Zvonilov
\paper Complex topological invariants of real algebraic curves on a hyperboloid and an ellipsoid
\jour St. Petersburg Math. J.
\vol 3
\yr 1992
\pages 1023--1042
\endref\label{Zvon}



\endRefs
\enddocument